# ON THE GLOBAL LINEAR ZARANKIEWICZ PROBLEM

PANTELIS E. ELEFTHERIOU AND ARIS PAPADOPOULOS

Abstract. The 'global' Zarankiewicz problem for hypergraphs asks for an upper bound on the number of edges of a finite $r$-hypergraph $V$ in terms of the number $|V|$ of its vertices, assuming the edge relation is induced by a fixed $K_{k,\dots,k}$-free $r$-hypergraph $E$, for some $k \in \mathbb{N}$. In [4], such bounds of size $O(|V|^{r-1})$ were achieved for a *semilinear* $E$, namely, definable in a linear o-minimal structure. We establish the same bounds in five new settings: when $E$ is definable in (a) a semibounded o-minimal structure and the vertex set of $V$ is 'sufficiently distant', (b) a model of Presburger arithmetic, (c) the expansion $\langle \mathbb{R}, <, +, \mathbb{Z} \rangle$ of the real ordered group by the set of integers, (d) a stable 1-based structure without the finite cover property, and (e) a locally modular regular type in a stable theory, such as the generic type of the solution set of the Heat differential equation.

Our methods include techniques for reducing Zarankiewicz's problem to the setting of arbitrary subgroups of powers of groups, used in geometric cases (a)–(c). They also include an abstract version of Zarankiewicz's problem for general saturated 'linear structures' that yields the desired bounds in the model-theoretic settings (d)–(e), as well as a parametric version in (b). Furthermore, the bounds in (a) characterise those o-minimal structures that do not recover a global field, and in (c) they yield new versions of Zarankiewicz's problem for certain ordered abelian groups.

*Date*: January 4, 2026.

The first author was partially supported by an EPSRC Early Career Fellowship (EP/V003291/1), and the second author by a Leeds Doctoral Scholarship, University of Leeds.

## 1. Introduction

The Zarankiewicz problem for hypergraphs is a central problem in extremal combinatorics. It asks for upper bounds on the number of edges of a hypergraph in terms of the number of its vertices, assuming it has no complete sub-hypergraphs of a given size. The *global* Zarankiewicz problem asks for the same upper bounds, assuming moreover that the edge relation of the hypergraphs we consider is induced by the edge relation of a *fixed* hypergraph with no complete sub-hypergraphs. A nice blend of recent model-theoretic and combinatorial techniques around structures without the independence property (NIP) yielded bounds in both problems in a variety of geometrically tame settings, including semialgebraic [27], o-minimal and distal [10, 53], and semilinear [4]. In particular, in the last reference a link was made between the global Zarankiewicz problem and a classical theme in model theory, that of recognising the existence of a definable field in a given structure from combinatorial data (see, for example, [33, 46]). A relevant combinatorial problem, that of Elekes-Szabó, with different Zarankiewicz-type bounds and recognising algebraic groups instead of fields, has also been considered in the o-minimal and other model-theoretic settings in [12]. The semilinear setting has seen further combinatorial interest in Basit-Tran [5], as well as interest from the perspective of valued fields and non-archimedean tame topology in Hrushovski-Loeser [32] (see also [14]). For further historical background on Zarankiewicz's problem for hypergraphs, the reader may consult the introductions in [4, 10, 27].

In this paper, we focus on the global Zarankiewicz problem and establish the same 'linear' bounds as in [4] in five new settings. Let $M$ be a set and $d_1, \dots, d_r \in \mathbb{N}_{>0}$, with $d_1 + \cdots + d_r = n$. By a $(d_1, \dots, d_r)$-*grid* $B$ on $M$, we mean a Cartesian product $B = B_1 \times \cdots \times B_r$ of sets $B_i \subseteq M^{d_i}$, $i \in [r]$. We omit '$(d_1, \dots, d_r)$-' when it is clear from the context. Given $k \in \mathbb{N}$, we call $B$ a $k$-*grid* if each $|B_i| = k$. We call $B$ *infinite* if each $B_i$ is infinite. An $r$-ary relation

$$E \subseteq \prod_{i\in[r]} M^{d_i} = M^n$$

induces on each grid $B$ an $r$-hypergraph $\mathcal{G}_B^E$; namely,

$$\mathcal{G}_B^E = \big(B_1, \dots, B_r; E \cap (B_1 \times \cdots \times B_r)\big).$$

If $B$ is finite, we write $n_B := \max_i |B_i|$. Clearly, $|E \cap B| \leq n_B^r$. We are interested in relations $E$ such that $|E \cap B| \leq n_B^{r-1}$, under the further condition that $E$ does not contain certain grids. To provide the definitions that apply to the entire paper, let us also write, following [4], for a finite grid $B$,

$$\delta(B) = \sum_{1\leq i_1 < \cdots < i_{r-1} \leq r} \left( \prod_{i\in[r-1]} |B_{i_j}| \right).$$

Note that $\delta(B) = O(n_B^{r-1})$, and hence one may think of $\delta(B)$ as $n_B^{r-1}$. By *Zarankiewicz bounds* for $E$, we mean bounds for the quantity $|E \cap B|$, as $B$ ranges over all finite grids from a class $\mathcal{C}$ of grids. A novelty of the present account is that $\mathcal{C}$ need not contain all grids. This allows for 'Zarankiewicz statements' in a variety of settings, such as in Theorem A below, where they would otherwise fail. Let us introduce these key notions, relativising the ones from [4] to a class $\mathcal{C}$ of grids.

**Definition 1.1.** Let $E \subseteq \prod_{i\in[r]} M^{d_i}$ be an $r$-ary relation, and $\mathcal{C}$ a class of grids.

(1) We say that $E$ has *linear Zarankiewicz bounds (linear $Z$-bounds) for $\mathcal{C}$* if there is $\alpha \in \mathbb{R}_{>0}$ such that for every finite $B \in \mathcal{C}$,

$$|E \cap B| \leq \alpha\delta(B).$$

In other words, the number of edges of a finite hypergraph $\mathcal{G}_B^E$, where $B \in \mathcal{C}$, is $O(n_B^{r-1})$. In this case, we say that these bounds are *witnessed* by $\alpha$. (If $E$ is binary, $O(n_B)$ is linear, supporting the terminology.)

(2) We say that $E$ is *$\mathcal{C}$-$\infty$-free* if it contains no infinite grid $B \in \mathcal{C}$. For $k \in \mathbb{N}$, we say that $E$ is *$\mathcal{C}$-$k$-free* if it contains no $k$-grid $B \in \mathcal{C}$.

(3) We abbreviate the following two *Zarankiewicz statements*:

$$Zar^{\infty}(E, \mathcal{C})\colon \text{ if } E \text{ is } \mathcal{C}\text{-}\infty\text{-free, then } E \text{ has linear } Z\text{-bounds for } \mathcal{C}.$$

$$Zar(E, \mathcal{C})\colon \text{ if there is } k \in \mathbb{N}_{>0}, \text{ such that } E \text{ is } \mathcal{C}\text{-}k\text{-free, then } E \text{ has linear } Z\text{-bounds for } \mathcal{C}.$$

We write $Zar_{\alpha}^{\infty}(E, \mathcal{C})$ or $Zar_{\alpha}(E, \mathcal{C})$, respectively, if $\alpha$ witnesses the conclusion of the corresponding statement. Of course, $Zar_{\alpha}^{\infty}(E, \mathcal{C})$ implies $Zar_{\alpha}(E, \mathcal{C})$. If $\mathcal{C}$ is the class of all grids, we omit $\mathcal{C}$ from all notation above, and say, for example that $E$ is $k$-free or that $Zar(E)$ holds. We also abbreviate the following *uniform* (that is, *parametric*) Zarankiewicz statements, for a family $\mathcal{E} = \{E_b\}_{b \in I}$ of relations

$$E_b \subseteq \prod_{i \in [r]} M^{d_i}.$$

$$Zar^{\infty}(\mathcal{E}, \mathcal{C})\colon \text{ there is } \alpha \in \mathbb{R}_{>0}, \text{ such that for every } b \in I,\ Zar_{\alpha}^{\infty}(E_b, \mathcal{C}) \text{ holds.}$$

$$Zar(\mathcal{E}, \mathcal{C})\colon \text{ there is } \alpha \in \mathbb{R}_{>0}, \text{ such that for every } b \in I,\ Zar_{\alpha}(E_b, \mathcal{C}) \text{ holds.}$$

If $\mathcal{C}$ is the class of all grids, we again omit $\mathcal{C}$ and write $Zar(\mathcal{E})$, $Zar^{\infty}(\mathcal{E})$.

Note that $E$ is $k$-free if and only if the hypergraph $\mathcal{G}_B^E$, with $B = \prod_{i \in [r]} M^{d_i}$, is $K_{k,\ldots,k}$-free in the usual graph-theoretic terminology.

For succinctness, we often refer to the Zarankiewicz problem in various settings by use of an adjective, such as linear/semibounded/Presburger/abstract Zarankiewicz. This is the problem of whether relations, or families thereof, definable in a structure from the corresponding setting, satisfy some of the Zarankiewicz statements above.

**In the rest of this paper, $\mathcal{M} = \langle M, \ldots \rangle$ denotes a structure, 'definable' means 'definable in $\mathcal{M}$ with parameters', and '$A$-definable' indicates the parameters are from $A \subseteq M$.** Definable sets in each of the first four settings below enjoy purely geometric characterisations (mentioned in this introduction and explored further in Section 2).

1.1. **Linear Zarankiewicz.** We restate the relevant results from [4]. Let $M = \mathbb{R}$. By a *semilinear* set $X \subseteq \mathbb{R}^n$, we mean a Boolean combination of solution sets in $\mathbb{R}^n$ to linear equalities and inequalities over $\mathbb{R}$.

**Fact 1.2** ([4, Corollary 5.12])**.** *Let $E \subseteq \prod_{i \in [r]} \mathbb{R}^{d_i}$ be a semilinear relation. Then $Zar^{\infty}(E)$ holds.*

In model-theoretic terminology, a set $X \subseteq \mathbb{R}^n$ is semilinear if and only it is definable in the real vector space $\mathbb{R}_{vec}$ over $\mathbb{R}$. In fact, $\mathbb{R}_{vec}$ is just one example of a *linear* o-minimal structure (Definition 2.11), yielding the above semilinear sets. Moreover, the full version of Fact 1.2 characterises linear o-minimal structures via a Zarankiewicz statement, as stated next. An advantage of employing model-theoretic methods is that we further obtain the parametric versions of the statements at hand, such as in (2) below, akin to many of the results of this paper and other known results in the literature (for example, [27]).

**Fact 1.3** ([4, Corollary 5.11])**.** *Let $\mathcal{M} = \langle M, <, \ldots \rangle$ be an o-minimal structure. Then the following are equivalent:*

*(1) $\mathcal{M}$ is linear.*
*(2) Let $\mathcal{E} = \{E_b\}_{b \in I}$ be a definable family of relations $E_b \subseteq \prod_{i\in[r]} M^{d_i}$. Then $Zar^\infty(\mathcal{E})$ holds.*
*(3) For every binary definable $E \subseteq M^{d_1} \times M^{d_2}$, $Zar(E)$ holds.*

*Remark* 1.4. As pointed out in [4], the above clauses are also equivalent to the fact that there is no infinite field definable in $\mathcal{M}$. This is an alternative characterisation of linearity for o-minimal structures ([35, 45]). The equivalence with (3) is due to the Szemerédi-Trotter lower bound $\Omega(n^{4/3})$ ([51]) for the point-line incidence relation, which is definable in the presence of an infinite field. Finally, the uniformity in (2) indicates that Zarankiewicz bounds only depend on the formula that defines each $E_b$, and not on the parameter $b$ (which is also implicit, say, in [27, Theorem 1.1]).

We now proceed to describe the results of this paper, in two fairly independent parts: Sections 1.2 – 1.4 (geometric settings) and Section 1.5 (model-theoretic settings).

1.2. **Semibounded Zarankiewicz.** Motivated by the themes in [4] and recent literature from the so-called *semibounded* o-minimal geometry, we prove a version of the global Zarankiewicz problem for *eventually linear* sets $E$, obtaining $Zar^\infty(E, \mathcal{C})$, for a class $\mathcal{C}$ of 'sufficiently distant' grids. We postpone the precise definition of a semibounded o-minimal structure $\mathcal{M}$ to Definition 2.11, but the reader can think of an eventually linear set $X \subseteq \mathbb{R}^n$ as a finite union of sets $V + D$, where $D \subseteq \mathbb{R}^n$ is bounded and $V$ the sum of unbounded segments of 1-dimensional $\mathbb{R}_{vec}$-subspaces. For example, the vertical cylinder $S^1 + V$, where $S^1$ is the unit circle on the $xy$-plane and $V$ is the segment $\{(0, 0, z) : 0 < z < \infty\}$ of the $z$-axis, is eventually linear. In model-theoretic terminology, a set $X \subseteq \mathbb{R}^n$ is eventually linear if and only if it is definable in an o-minimal expansion of $\mathbb{R}_{vec}$ by bounded sets. By [20], a semibounded $\mathcal{M}$ may define bounded infinite fields (but not unbounded ones), and hence by Remark 1.4, we cannot expect $Zar^\infty(E, \mathcal{C})$ to hold with $\mathcal{C}$ the class of all grids. One could try to remedy this problem by varying the Zarankiewicz bounds, but in this work we preserve the linear $Z$-bounds by restricting $\mathcal{C}$ to the class of all 'sufficiently distant' grids.

**Definition 1.5.** Let $\mathcal{M} = \langle M, <, +, \ldots \rangle$ be an expansion of an ordered group, and $m \in M_{\geq 0}$. Two elements $x, y \in M^n$ are called *$m$-distant* if $|x - y| > m$. A set $X \subseteq M^n$ is called *$m$-distant* if any two distinct $x, y \in X$ are $m$-distant. A grid $B = B_1 \times \cdots \times B_r$ is called *$m$-distant* if every $B_i$ is an $m$-distant set. We denote by $\mathcal{C}_m$ the class of all $m$-distant grids.

Using semibounded geometry, we are able to prove the following results.

**Theorem A.** *Let $E \subseteq \prod_{i\in[r]} \mathbb{R}^{d_i}$ be an eventually linear relation. Then there is $m \in \mathbb{R}_{>0}$, such that $Zar^\infty(E, \mathcal{C}_m)$ holds.*

As in the semilinear setting, our full theorem includes a parametric version of the Zarankiewicz statement. Note that a semibounded structure need not be a reduct of a real closed field, for example $\langle \mathbb{R}, <, +, \sin_{\restriction(0,1)}\rangle$.

**Theorem A′** (3.12, 3.13)**.** *Let $\mathcal{M} = \langle M, <, +, \dots\rangle$ be a semibounded expansion of an ordered group. Let $\{E_b\}_{b\in I}$ be a definable family of relations $E_b \subseteq \prod_{i\in[r]} M^{d_i}$. Then there are $\alpha \in \mathbb{R}_{>0}$ and $\{m_b\}_{b\in I}$ with $m_b \in M_{\geq 0}$, such that for every $b \in I$, $Zar^{\infty}_{\alpha}(E_b, \mathcal{C}_{m_b})$ holds.*

Coming to reducts of real closed fields, we characterise those that do not recover the full multiplication via a Zarankiewicz statement; this is perhaps surprising, since we manage to capture the unbounded behavior of a structure using finite grids.

**Theorem A″** (3.16)**.** *Let $\mathcal{M} = \langle M, <, +, \dots\rangle$ be a reduct of a real closed field $\mathcal{R} = \langle M, <, +, \cdot\rangle$. Then the following are equivalent:*

1. *The multiplication $\cdot$ is not definable in $\mathcal{M}$.*
2. *For every binary definable $E \subseteq M^{d_1} \times M^{d_2}$, there is $m \in M_{\geq 0}$, such that $Zar(E, \mathcal{C}_m)$ holds.*

*Remark* 1.6. The condition of being semibounded in Theorem A′ is equivalent to the fact that there is no definable ordered field with domain $M$ whose ordered agrees with $<$ ([20]). In fact, Corollary 3.12 yields Zarankiewicz statements sensitive to the (non-)existence of definable fields on intervals of any given length, bounded or unbounded. Theorem A″ is again due to an 'unbounded version' of Szemerédi-Trotter's lower bound theorem for any ordered field (Proposition 3.15). We do not know if Theorem A″ is true without assuming $\mathcal{M}$ is a reduct of the field $\mathcal{R}$ (and recovering in (1) *some* unbounded field); on the other hand, we also do not know of any semibounded structure $\mathcal{M}$ that does not satisfy this assumption (see also Remark 3.17 and Question 9.1).

1.3. **Presburger Zarankiewicz.** By a *Presburger set* $X \subseteq \mathbb{Z}^n$, we mean a Boolean combination of solution sets in $\mathbb{Z}^n$ to linear equalities and inequalities over $\mathbb{Q}$. In model-theoretic terminology, $X$ is Presburger if and only if it is definable in the structure $\langle \mathbb{Z}, <, +\rangle$. Presburger arithmetic, denoted by $\mathsf{Pres}$, is the theory of that structure. If $\mathcal{M} = \langle M, <, +\rangle \models \mathsf{Pres}$ is a model of Presburger arithmetic, then $\mathcal{M}$ is an 'elementary extension' of $\langle \mathbb{Z}, <, +\rangle$ (Section 2.1). We also call *Presburger sets* the sets definable in $\mathcal{M}$.

**Theorem B** (4.31)**.** *Let $\mathcal{M} = \langle M, <, +\rangle \models \mathsf{Pres}$, and $E \subseteq \prod_{i\in[r]} M^{d_i}$ a definable set. Then $Zar(E)$ holds.*

We observe that $Zar^{\infty}(E)$ need not be true (Remark 4.28), due to the fact that Presburger arithmetic does not 'eliminate $\exists^{\infty}$'. Likewise, and more importantly, if $\mathcal{E}$ is a definable family, $Zar(\mathcal{E})$ need not hold. For example, let $E_t = [0, t) \subseteq M$, $t \in M_{>0}$. Then for every $t \in \mathbb{Z}$, there is $k \in \mathbb{N}$ such that $E_t$ is $k$-free, but $E_t$ itself is a grid of size $t$. Hence there is no $\alpha \in \mathbb{R}_{>0}$, such that for every $t \in \mathbb{Z}$, $|E_t| \leq \alpha t^{1-1} = \alpha$. We remedy this problem in the second part of this paper, where we prove a parametric Zarankiewicz statement for the class of all '$\mathbb{Z}$-distant grids' in the 'saturated' setting.[1] See Section 6.1 for the definitions of a saturated structure and a type-definable set.

[1] By 'saturated', in this paper, we mean 'sufficiently saturated'.

**Definition 1.7.** Let $\mathcal{M} \models \mathsf{Pres}$. A set $X$ or grid $B$ is called $\mathbb{Z}$-*distant* if it is $m$-distant for every $m \in \mathbb{N}$.

**Theorem B′** (8.4)**.** *Let $\mathcal{M} = \langle M, <, +\rangle \models \mathsf{Pres}$ be saturated, $\mathcal{E} = \{E_b\}_{b\in I}$ a type-definable family of relations $E_b \subseteq \prod_{i\in[r]} M^{d_i}$, and $\mathcal{C}_\mathbb{Z}$ the class of all $\mathbb{Z}$-distant grids. Then $Zar^\infty(\mathcal{E}, \mathcal{C}_\mathbb{Z})$ holds.*

Note that in the example above, if $\mathcal{M}$ is saturated, then $E_t = [0, t)$ is $\mathcal{C}_\mathbb{Z}$-$k$-free for some $k \in \mathbb{N}$ if and only if $t \in \mathbb{N}_{>0}$. In that case $|E_t \cap B| \leq 1$, for any $\mathbb{Z}$-distant grid $B$.

1.4. **Mixed real-integer Zarankiewicz.** Our third setting concerns the expansion $\langle \mathbb{R}, <, +, \mathbb{Z}\rangle$ of the real ordered group by the set of integers. This structure combines the semilinear and Presburger geometries, and has been studied from at least two different perspectives: as a natural environment to host applications to linear and integer programming (see [34, 55, 56], and [38] for more general 'real-integer systems'), and as a prototype of 'tame expansions of o-minimal structures' – a model-theoretic environment originated by A. Robinson [49] and grown substantially ever since (for example, in [6, 7, 15, 16, 18, 22, 39]). Definable sets in this setting are unions of semilinear families of sets where the parameters can vary over any Presburger set (Fact 5.9). Thus, our motivation for establishing a Zarankiewicz statement in $\langle \mathbb{R}, <, +, \mathbb{Z}\rangle$ is also two-fold: (a) to yield new parametric versions of linear Zarankiewicz (Fact 1.3) (Corollary 5.23), and (b) to demonstrate how combinatorial considerations can extend to more general tame settings. The latter aspect has potentially a nice range of applications, already illustrated here by new Zarankiewicz statements for some further ordered abelian groups (Corollary 5.24). It also triggers some curious questions in the tame setting (Section 9.4).

**Theorem C** (5.21)**.** *Let $\mathcal{M} = \langle \mathbb{R}, <, +, \mathbb{Z}\rangle$, and $E \subseteq \prod_{i\in[r]} \mathbb{R}^{d_i}$ a definable relation. Then $Zar(E)$ holds.*

**On the methods so far.** Our proofs in the above three geometric settings involve a reduction of the Zarankiewicz statements to the setting of arbitrary subgroups of a power $\langle M^n, +\rangle$ of a group $\langle M, +\rangle$ (literally, subgroups of $\prod_{i\in[r]}\langle M^{d_i}, +\rangle$), which are independently handled in Theorem 2.7. This reduction involves the notion of a 'shell' of a set $E$, namely the subgroup of $\langle M^n, +\rangle$ generated by $E$, and essentially consists of proving Property $(*)$ from Section 2.7: if $E$ is $k$-free, or $\infty$-free, then so is its shell. In our three geometric settings, the reduction is manifested in various ways, which we describe further in the Reduction Strategy of Section 2.7. As a by-product of this machinery, we obtain an alternative, model theory-free proof of Fact 1.3, after establishing $(*)$ for linear cells containing 0 (Corollaries 2.21, 2.22). We also note that for every $\infty$-free semilinear set $E$, the Zarankiewicz bounds are $O(n_B^{\mathsf{dim} E})$, answering a question by S. Starchenko (Remark 2.23(2)). Notably, a major part of our work towards Theorems B′ and C is carried out at the general level of definably complete ordered groups (Section 4), where we introduce the notion of an abstract cell and prove a Zarankiewicz statement for itself (Proposition 4.29).

1.5. **Abstract Zarankiewicz.** In this second part, we prove an abstract version of Zarankiewicz's problem for an arbitrary saturated structure $\mathcal{M} = \langle M, \dots\rangle$, an operator $\mathsf{cl}$ (that is, a map $\mathsf{cl} : \mathcal{P}(M) \to \mathcal{P}(M)$), and a class $\mathcal{C}$ of grids satisfying some key properties (Theorem D). As an application, we obtain Zarankiewicz statements in the settings of stable 1-based theories without the finite cover property (Theorem E), and locally modular regular types in a stable theory (Theorem F). We further deduce Theorem B′ and strengthen Theorem A′ in the saturated setting (Theorem 8.9).

The key properties consist of (DEF), (UB), and (TIGHT), which we introduce in Section 6. The operator $\mathsf{cl}$ induces a notion of independence, and a grid $B = B_1 \times \cdots \times B_r$ is called $\mathsf{cl}$-*independent over* $A$ if each $B_i$ is $\mathsf{cl}$-independent over $A$. Property (TIGHT) captures the 'linearity' or 'modularity' nature of $\mathcal{M}$, and is versatile enough so that the following theorem applies to different settings.

**Theorem D** (6.7)**.** *Let $\mathcal{M}$ be a saturated structure, $\mathsf{cl}$ an operator, and $\mathcal{C}$ a class of grids. Suppose* (DEF), (UB) *and* (TIGHT)*. Let $\mathcal{E} = \{E_b\}_{b\in I}$ be a type-definable family of relations $E_b \subseteq \prod_{i\in[r]} M^{d_i}$. Then $Zar^\infty(\mathcal{E},\mathcal{C})$ holds.*

When moreover $(M, \mathsf{cl})$ is a pregeometry, we obtain Theorem 6.17, which has the same assumptions with Theorem D, except that $\mathcal{C}$ is required to contain all infinite $\mathsf{cl}$-independent grids, and $\mathsf{cl}$ is 'weakly locally modular' (WLM) instead of (TIGHT).

Theorem D yields the following application, with $\mathsf{cl} = \mathsf{acl}^{eq}$.

**Theorem E** (7.2)**.** *Let $T$ be a stable 1-based theory without the finite cover property, and $\mathcal{M}$ a model of $T$. Let $\mathcal{E} = \{E_b\}_{b\in I}$ be a type-definable family of relations $E_b \subseteq \prod_{i\in[r]} M^{d_i}$. Then $Zar(\mathcal{E})$ (and $Zar^\infty(\mathcal{E})$, if $\mathcal{M}$ is saturated) holds.*

As a corollary (7.10), we obtain that the 'ab initio Hrushovski constructions' enjoy the Zarankiewicz statements.

Theorem 6.17 yields the following applications. First, it implies a known result ([4, Theorem 5.6]), if $\mathcal{M}$ eliminates $\exists^\infty$, $\mathsf{cl} = \mathsf{acl}$ if a pregeometry and $\mathcal{C}$ is the class of all grids. However, in its generality, Theorem 6.17 applies to settings without elimination of $\exists^\infty$, such as Presburger arithmetic (Theorem B′), and to settings where $\mathsf{acl}$ is not locally modular, such as semibounded structures (Theorem 8.9), and the 'locus' of locally modular regular types in stable theories (Theorem F).

**Theorem F** (7.16)**.** *Let $T$ be a stable theory, $\mathcal{M}$ a model of $T$, $p$ a locally modular regular type, $U = p(\mathcal{M})$ its set of realisations in $\mathcal{M}$, and $\mathsf{cl}_p$ the forking closure operator on $U$. Let $\mathcal{U}$ be the induced structure on $U$ by $\mathcal{M}$, and $\mathcal{C}$ the class of all Morley grids in $p$. Let $\mathcal{E} = \{E_b\}_{b\in I}$ be a type-definable (in $\mathcal{U}$) family of relations $E_b \subseteq \prod_{i\in[r]} U^{d_i}$. Then $Zar(\mathcal{E},\mathcal{C})$ (and $Zar^\infty(\mathcal{E},\mathcal{C})$, if $\mathcal{M}$ is saturated) holds.*

Theorem F applies to the theory $\mathsf{DCF}_{0,2}$ of differentially closed fields of characteristic 0 with two commuting derivations, and $p$ the generic type of the Heat variety $\delta_1 y = \delta_2^2 y$.

*Remark* 1.8. Property $(*)$ from Section 2.7 and (TIGHT) can be viewed as two new notions of 'linearity' in our geometric and model-theoretic settings, respectively (see also Examples 2.16 and 6.5).

**Structure of the paper.** In Section 2, we introduce the key notions of this paper, such as shells, linear and semibounded o-minimal structures, Presburger arithmetic, and our Reduction Strategy to subgroups of arbitrary groups. In Section 3, using the Reduction Strategy, we settle Theorems A′ and A″. In Section 4, we carry out some analysis of topological nature in arbitrary definably complete ordered groups with a fixed copy $\mathbf{Z}$ of the integers, after introducing the notion of $N$-internal points and arbitrary $\mathbf{Z}$-cells. We use this analysis to implement the Reduction Strategy and settle Theorem B in Section 4.5, and Theorem C in Section 5. For the latter, we also use the work in Appendix A, which is carried out for arbitrary ordered vector spaces over ordered division rings. In Section 6, we turn to the abstract model-theoretic setting, proving Theorems D and 6.17, which we use to settle Theorems E and F in Section 7. In Section 8, we return to our geometric settings, and use Theorem 6.17 to settle Theorem B′, as well as obtain a stronger version of Theorem A′, namely, Theorem 8.9. Section 9 contains a list of open questions.

**Acknowledgments.** We thank Pablo Andújar Guerrero, Artem Chernikov, Joshua Losh, Dugald Macpherson, Vincenzo Mantova, Mervyn Tong and Yilong Zhang for multiple discussions on the topics of this paper. We thank Amador Martin-Pizarro and Rosario Mennuni for helping with the proof of Proposition 7.7, Mario Edmundo, Ya'acov Peterzil and Sergei Starchenko for several discussions on semibounded o-minimal structures, Alexander Berenstein and Evgueni Vassiliev on weak local modularity, and Joel Nagloo, Rémi Jaoui, Anand Pillay on regular types. We thank Mariana Vicaria for suggesting the Zarankiewicz problem in Presburger arithmetic.

## 2. Preliminaries and first results

In this section, we introduce some key notions that are used throughout this paper, such as grid-configurations and shells. We prove some preliminary results around Zarankiewicz's problem, moving gradually from arbitrary structures, to groups, then linear and semibounded o-minimal structures, and finally Presburger arithmetic. We present (Section 2.7) a Reduction Strategy for Zarankiewicz's problem for our three geometric settings from the introduction to Theorem 2.7 for arbitrary subgroups of powers of groups. We employ this reduction in Sections 3–5. As a by-product, we obtain an alternative proof of Fact 1.3 (Corollary 2.22).

### 2.1. **Notation and terminology.**

*For the rest of this paper, and unless stated otherwise, the following notation and terminology will be fixed.*

- $M$ is a set, usually the domain of a structure
- $r, n, d_1, \dots, d_r \in \mathbb{N}_{>0}$, with $d_1 + \cdots + d_r = n$.
- $E$ is an $r$-ary relation, that is a subset of $\prod_{i\in[r]} M^{d_i} = M^n$.
- By a grid, we mean a $(d_1, \dots, d_r)$-grid.
- For $j = 1, \dots, r$, we denote by
$$p_j : \prod_{i\in[r]} M^{d_i} \to \prod_{i\neq j} M^{d_i}$$
the projection that omits the '$j$-th block of coordinates', namely the $d_j$-many coordinates that occur at the $j$-th component of $\prod_{i\in[r]} M^{d_i}$.
- For $x \in M^n$, we denote by $x^i \in M^{d_i}$ the tuple of its $i$-th block of coordinates.
- For $x \in M^{d_j}$, we denote by
$$\hat{x}^j = (0, \dots, 0, x, 0, \dots, 0) \in \prod_{i\in[r]} M^{d_i}$$
the tuple consisting of zeros everywhere except for the $j$-th block of coordinates, which equals $x$.

**Further standard notation.** For $q \in \mathbb{N}_{>0}$, we write $[q] = \{1, \dots, q\}$, and $[0] = \emptyset$. A tuple of elements is denoted just by one element, and we write $b \subseteq B$ if $b$ is a tuple with coordinates from $B$. If $A, B \subseteq M$, we may write $AB$ for $A \cup B$, and $Ab$ for $A\{b\}$. By $|A|$ we denote the cardinality of $A$.

Given a set $X \subseteq M^m \times M^n$ and $a \in M^m$, we write $X_a$ for

$$\{b \in M^n \ : \ (a, b) \in X\}.$$

We deal with many coordinate projections, but unless stated otherwise, $\pi : M^n \to M^{n-1}$ denotes the projection onto the first $n-1$ coordinates. We routinely identify $X^k \times X^l$ with $X^{k+l}$. A family of sets $\mathcal{X} = \{X_t\}_{t \in I}$, $I \subseteq M^m$, $X_t \subseteq M^n$, is called (type-)definable in a given structure $\mathcal{M} = \langle M, \dots \rangle$ if the set

$$\bigcup_{t \in I} \{t\} \times X_t \subseteq M^{m+n}$$

is (type-)definable.

If $\langle M, +, \dots \rangle$ expands a group, we denote by 0 the group identity element. We also write 0 for the origin $(0, \dots, 0) \in M^n$, for a given $n$. By convention, $M^0 = \{0\}$, and we identify $(x, y) \in M^0 \times M^n$ with $y$. Perhaps a feature of this paper is that we often start recursive definitions and inductive proofs with $n = 0$. For $n \in \mathbb{N}$ and $t \in M$, we denote by $nt$ the sum of $t$ with itself $n$ times. For $n \in \mathbb{N}$ and $A \subseteq M$, we denote $nA = \{nt : t \in A\}$.

By an *order*, we mean a total order. If $\langle M, <, \dots \rangle$ is an ordered structure, we extend $<$ to $M \cup \{\pm\infty\}$ in the standard way. For $x \in M$, $|x|$ denotes its absolute value. For $x = (x_1, \dots, x_n) \in M^n$, we let $|x| = \max_i |x_i|$. Given $f, g : X \subseteq M^n \to M$ two functions, or $f = -\infty$ or $g = +\infty$ (viewed as constant functions), with $f(x) < g(x)$ for all $x \in X$, we write $f < g$ and define the cylinders

$$(f, g) = \{(x, y) : x \in X, f(x) < y < g(x)\}$$

and

$$[f, g] = \{(x, y) : x \in X, f(x) \leq y \leq g(x)\},$$

where, by convention, $[-\infty, \infty] = (-\infty, \infty)$, $[-\infty, g] = (-\infty, g]$ and $[f, \infty] = [f, \infty)$. We denote by $\Gamma(f)$ the graph of $f$, and by $Im(f)$ its image. If $Y \subseteq X$, we write $\Gamma(f)_Y$ for $\Gamma(f_{\restriction Y})$, $(f, g)_Y$ for $(f_{\restriction Y}, g_{\restriction Y})$ and similarly for the rest of the cylinders. We call those sets the *restrictions* of the corresponding set to $Y$. We write $C = \Gamma(f)$ for $C = \Gamma(f)_{\pi(C)}$ and $C = (f, g)$ for $C = (f, g)_{\pi(C)}$. Given $x \in M$ and $A \subseteq M$, we write $x < A$ if $x < y$ for every $y \in A$, and similarly for $A < x$. We also denote $A_{>x} = \{y \in A : x < y\}$ and similarly for $A_{\geq x}, A_{<x}, A_{\leq x}$. We endow $M$ with the order topology, whose basic open sets are the open intervals, and each $M^n$ with the product topology. Given $A \subseteq M^n$, we denote by $\overline{A}$ its topological closure.

If $\langle M, <, +, \dots \rangle$ expands an ordered group, $x \in M^n$ and $r \in M$, we denote by $B_r(x) = x + [-r, r]^n$ the box around $x$ of radius $r$. By an *open box in* $M^n$, we mean a set of the form $B_1 \times \cdots \times B_n$ where each $B_i$ is an open interval. We extend addition to $M \cup \{\pm\infty\}$ in the standard way, for example $\infty + a = \infty$ and $a - \infty = -\infty$.

For general background in logic and model theory, the reader can consult [36, 46, 52]. Given a language $\mathcal{L}$, an $\mathcal{L}$*-theory* is a consistent set of $\mathcal{L}$-sentences, and it is *complete* if it is maximal such. The theory of an $\mathcal{L}$-structure $\mathcal{M}$ is the set of all $\mathcal{L}$-sentences that are true in $\mathcal{M}$, and it is complete. An $\mathcal{L}$-structure $\mathcal{M}$ *eliminates* $\exists^\infty$ if for every $\mathcal{L}$-formula $\phi(x, y)$ $(|x| = 1)$, there is $N \in \mathbb{N}$, such that for every $b \in M^{|y|}$, if the realisations of $\phi(x, b)$ in $\mathcal{M}$ are finitely many, then they are at most $N$. For two $\mathcal{L}$-structures $\mathcal{M}, \mathcal{N}$, we write $\mathcal{M} \preccurlyeq \mathcal{N}$ if $\mathcal{M}$ is an *elementary substructure* of $\mathcal{N}$; that is, every $\mathcal{L}$-sentence with parameters from $M$ is true in $\mathcal{M}$ if and only if it is true in $\mathcal{N}$. A structure $\mathcal{M} = \langle M, \dots \rangle$ is a *reduct* of a structure $\mathcal{N} = \langle N, \dots \rangle$ if $M = N$ and every $\emptyset$-definable subset of $\mathcal{M}$ is $\emptyset$-definable in $\mathcal{N}$.

For the basics of o-minimality, such as dimension and cell decomposition, the reader is referred to [17], although most relevant background is introduced as we move along. Below, we define 'cells' in various contexts. We keep the terminology *linear cell*, *Presburger cell*, and $\mathbf{Z}$-*cell*, for the settings of linear o-minimal structures, models of $\mathsf{Pres}$, and arbitrary definably complete ordered groups in Sections 2.5, 2.6 and 4.2, respectively.

2.2. **Grid-configurations and basic lemmas.** The starting point for obtaining linear $Z$-bounds for a set $E$ is the following.

**Lemma 2.1.** *Let $l \in [r]$. Suppose there is $N \in \mathbb{N}$, such that for every $B \in \mathcal{C}$, the map $p_{l \restriction E \cap B}$ is $N$-to-$1$. Then $E$ has linear $Z$-bounds for $\mathcal{C}$, witnessed by $N$.*

*Proof.* Denote $\pi = p_l$. For every finite grid $B \in \mathcal{C}$, we have

$$|E \cap B| \leq N|\pi(E) \cap \pi(B)| \leq N|\pi(B)| = N \prod_{i \in [r] \setminus \{l\}} |B_i| \leq N\delta(B),$$

as needed. $\square$

The next remark becomes handy when dealing with unions of sets.

*Remark* 2.2. Suppose $E = E_1 \cup \dots \cup E_l$, and $\mathcal{C}$ is a class of grids. If $Zar(E_i, \mathcal{C})$ holds for every $i$, then so does $Zar(E, \mathcal{C})$. Indeed, if $E$ is $\mathcal{C}$-$k$-free for some $k \in \mathbb{N}_{>0}$, then so is each $E_i$, and hence it has linear $Z$-bounds for $\mathcal{C}$, say witnessed by $\alpha_i \in \mathbb{R}_{>0}$. Then $\alpha = \sum_{i=1}^{l} \alpha_i$ witnesses that $E$ has linear $Z$-bounds for $\mathcal{C}$.

In various reductions, we replace grids by others of the same 'configuration'.

**Definition 2.3.** Let $M, N$ be two sets, $A \subseteq \prod_{i \in [r]} M^{d_i}$ and $B \subseteq \prod_{i \in [r]} N^{d_i}$. We say that $A = \{x_1, \dots, x_l\}$ and $B$ *have the same grid-configuration* if there is a bijection $f : A \to B$, such that for every $j, k \in [l]$ and $i \in [r]$, we have

$$x_j^i = x_k^i \iff f(x_j)^i = f(x_k)^i.$$

We denote this fact by $A \sim B$, and say that $f$ *witnesses* it.

Clearly, $\sim$ is an equivalence relation. Note that if $A \sim B$ is witnessed by $f$, and $C \subseteq A$, then $C \sim f(C)$.

2.3. **Group-theoretic Zarankiewicz.** We prove a key tool for handling the geometric settings in Sections 3–5, namely Theorem 2.7, which establishes $Zar(V)$ for an arbitrary subgroup $V$ of a power of a group. More precisely, let $\langle M, + \rangle$ be a group (not necessarily abelian), and denote by 0 its identity element. Let $V \subseteq \prod_{i \in [r]} M^{d_i}$ be a subgroup $V \leqslant \langle M^n, + \rangle$ (literally, a subgroup of $\prod_{i \in [r]} \langle M^{d_i}, + \rangle$). Clearly, for $l \in [r]$, $p_{l \restriction V}$ is finite-to-1 if and only if it is uniformly finite-to-1 (there is $N \in \mathbb{N}$ such that $p_{l \restriction V}$ is $N$-to-1).

*Remark* 2.4. In case there is an order $<$ such that $\langle M, <, + \rangle$ is an o-minimal ordered group or a model of $\mathsf{Pres}$, the above statements are moreover equivalent to $p_{l \restriction V}$ being injective.

For $l = 1, \ldots, r$, we define

$$U_l^V = \{x \in M^{d_l} : \hat{x}^l \in V\} \leqslant M^{d_i}.$$

We omit the superscript '$V$' if it is clear from the context. Observe that

$$\{0\} \times U_l \times \{0\} = \ker p_{l\restriction V} \leqslant V,$$

where the first $\{0\} \subseteq M^{d_1} \times \cdots \times M^{d_{l-1}}$, and the second $\{0\} \subseteq M^{d_{l+1}} \times \cdots \times M^{d_r}$. The following fact is then immediate.

**Fact 2.5.** *Let $V \leqslant \langle M^n, +\rangle$, and $l \in [r]$. Then*

$$p_{l\restriction V} \text{ is finite-to-1} \;\Leftrightarrow\; U_l \text{ is finite.}$$

**Lemma 2.6.** $U_1 \times \cdots \times U_r \subseteq V$.

*Proof.* We just need to observe that

$$U_1 \times \cdots \times U_r = (U_1 \times \{0\}) + \cdots + (\{0\} \times U_l \times \{0\}) + \cdots + (\{0\} \times U_r)$$

is a direct sum of subgroups of $V$. □

**Theorem 2.7.** *Let $V \subseteq \prod_{i\in[r]} M^{d_i} \subseteq M^n$ be a subgroup of $\langle M^n, +\rangle$. Then the following are equivalent:*

*(1) $V$ is $\infty$-free.*
*(2) there is $l \in \{1, \ldots, r\}$, such that $p_{l\restriction V}$ is finite-to-1.*
*(3) $V$ has linear $Z$-bounds for the class of all grids.*

*In particular, $Zar^\infty(V)$ holds. If, moreover, there is an order $<$ such that $\langle M, <, +\rangle$ is an o-minimal ordered group or a model of* Pres*, we can replace 'finite-to-1' by 'injective', and $Zar^\infty(V)$ by $Zar_1^\infty(V)$.*

*Proof.* (1)⇒(2). If each $p_{l\restriction V}$ is not finite-to-1, then by Fact 2.5 every $U_l$ is infinite. It follows from Lemma 2.6 that $V$ is not $\infty$-free.

(2)⇒(3). By Lemma 2.1.

(3)⇒(1). Immediate from the definitions.

The 'moreover' clause is by Remark 2.4. □

2.4. **Shells and linear maps.** Let $\langle M, +, \ldots\rangle$ be an expansion of a group. A key notion of this paper is the following.

**Definition 2.8.** Let $C \subseteq M^n$ be a set. The *shell of $C$*, denoted by $Sh(C)$, is the subgroup of $\langle M^n, +\rangle$ generated by $C$.

**Note:** Formally, the terminology should include a reference to +, but we omit this for simplicity.

Since the projection $\pi : \langle M^n, +\rangle \to \langle M^{n-1}, +\rangle$ is a group homomorphism, we obtain

$$\pi(Sh(C)) = Sh(\pi(C)).$$

Indeed, an element of the above sets is of the form

$$\pi(x_1 \pm \cdots \pm x_l) = \pi(x_1) \pm \cdots \pm \pi(x_l),$$

for $x_i \in C$.

We will be using throughout the notion of a linear map.

**Definition 2.9.** A map $f : X \subseteq M^m \to M^k$ is called *linear* (or, rather, *affine*) if for every $x, y, x+t, y+t \in X$,

$$f(x+t) - f(x) = f(y+t) - f(y).$$

We denote by $L(X)$ the set of all linear maps $f : X \subseteq M^m \to M$ and let $L_\infty(X) = L(X) \cup \{\pm\infty\}$, where we regard $-\infty$ and $+\infty$ as constant functions on $X$.

If $f : X \subseteq M^m \to M^k$ is a linear map with $0 \in X$, then clearly $f(0) = 0$ if and only if for every $x, y, x+y \in X$, $f(x+y) = f(x) + f(y)$. In this case, $f$ can be extended (uniquely) to a linear map $\hat{f} : Sh(X) \to M^k$ via

$$\hat{f}(x_1 \pm \cdots \pm x_l) = f(x_1) \pm \cdots \pm f(x_l).$$

We keep writing $f$ for $\hat{f}$.

**Lemma 2.10.** *Let $f : C \subseteq M^m \to M^k$ be a linear map with $f(0) = 0$. Then*

$$Sh(\Gamma(f)_C) = \Gamma(f)_{Sh(C)}.$$

*Proof.* Clear, since an element of the above sets has form

$$(x_1, f(x_1)) \pm \cdots \pm (x_l, f(x_l)) = (x_1 \pm \cdots \pm x_l, f(x_1) \pm \cdots \pm f(x_l)),$$

for $x_i \in C$. □

2.5. **Linear and semibounded o-minimal structures.** In this subsection we introduce the settings of linear and semibounded o-minimal structures, and prove some lemmas for the former setting that will be used in the sequel. The primary example of a linear o-minimal structure is that of an ordered vector space over an ordered division ring, such as $\mathbb{R}_{vec} = \langle \mathbb{R}, <, +, 0, \{x \mapsto \lambda x\}_{\lambda \in \mathbb{R}} \rangle$. The primary example of a semibounded non-linear structure is the expansion $\mathbb{R}_{sbd}$ of $\mathbb{R}_{vec}$ by all bounded semialgebraic sets. In both examples, the collection of all maps $\{x \mapsto \lambda x\}_{\lambda \in \mathbb{R}}$ forms a ring, that of all definable endomorphisms (with addition and composition).

Let $\mathcal{M} = \langle M, <, +, \dots \rangle$ be an o-minimal expansion of an ordered group. By a *partial endomorphism* we mean a linear map $f : (-c, c) \to M$ with $f(0) = 0$. We denote by $\Lambda$ the set of all $\emptyset$-definable partial endomorphisms, and by $\mathcal{B}$ the collection of all bounded definable sets.

**Definition 2.11.** We say that:

(a) $\mathcal{M}$ is *linear* ([35]) if every definable set in $\mathcal{M}$ is already definable in the structure $\langle M, <, +, 0, \{\lambda\}_{\lambda \in \Lambda} \rangle$.
(b) $\mathcal{M}$ is *semibounded* ([20, 42]) if every definable set in $\mathcal{M}$ is already definable in the structure $\langle M, <, +, 0, \{\lambda\}_{\lambda \in \Lambda}, \{B\}_{B \in \mathcal{B}} \rangle$.

Definable sets in (a) are also called *semilinear*, and in (b) *eventually linear*.

Obviously, if $\mathcal{M}$ is linear, then it is semibounded. By [45], $\mathcal{M}$ is not linear if and only if there is an infinite definable field. By [20], $\mathcal{M}$ is not semibounded if and only if there is an unbounded definable field whose order agrees with $<$ if and only if $\mathcal{M}$ expands a real closed field if and only if every definable map $f : M \to M$ is linear on some $(m, \infty) \subseteq M$. By [37, 42, 48], $\mathbb{R}_{sbd}$ is the unique structure that lies strictly between $\mathbb{R}_{vec}$ and the real field (in terms of their classes of definable sets). We postpone any further background on semibounded structures until Section 3, and prove here some lemmas for linear o-minimal structures.

*In the rest of this subsection, $\mathcal{M} = \langle M, <, +, \dots \rangle$ is a linear o-minimal structure.*

Our next goal is Lemma 2.13, which along with Fact 2.20, will be used in the proof of Corollary 2.22 (linear Zarankiewicz), and along with Fact 3.9, in that of Proposition 3.10 (semibounded Zarankiewicz). As a byproduct, we obtain a direct proof of Fact 1.3 (Corollary 2.22). A definable set $X \subseteq M^n$ is called *definably connected*, if for any $x, y \in X$, there is a definable continuous map $\gamma : [0, q] \subseteq M \to X$ with $\gamma(0) = x$ and $\gamma(q) = y$.

Consider the following property for a semilinear set $D \subseteq \prod_{i\in[r]} M^{d_i} = M^n$ containing 0:

(†) for every $l \in [r]$,
  (a) the set $D \cap \{\hat{x}^l : x \in M^{d_l}\}$ is definably connected.
  (b) If $p_{l\restriction D}$ is not injective, the set $D \cap \{\hat{x}^l : x \in M^{d_l}\}$ is infinite.
  (c) If $p_{l\restriction D}$ is injective, then, after permuting blocks of coordinates, $D = \Gamma(f)_C$ for some linear map $f : C \subseteq \prod_{i\in[r-1]} M^{d_i} \to M^{d_r}$.

Examples of semilinear sets satisfying (†) include 'linear cells' (Fact 2.20) and sets $\widetilde{V}$ as in Fact 3.9.

**Lemma 2.12.** *Let $D \subseteq \prod_{i\in[r]} M^{d_i} = M^n$ be a semilinear set with $0 \in D$ satisfying (†), and $l \in [r]$.*

*(1) If $p_{l\restriction D}$ is not injective, then for every $\varepsilon \in M_{>0}$, $B_\varepsilon(0) \cap U_l$ is infinite, where $U_l = U_l^{Sh(D)}$ (as in Section 2.3).*
*(2) If $p_{l\restriction D}$ is injective, then so is $p_{l\restriction Sh(D)}$.*
*(3) $D$ is open in $Sh(D)$.*

*Proof.* (1). By (†)(b), $D \cap \{\hat{x}^l : x \in M^{d_l}\}$ is infinite. By (†)(a), it is also definably connected, hence since $0 \in D$, it must contain elements arbitrarily close to the origin. That is, for every $\varepsilon \in M_{>0}$, there is $x \in M^{d_l} \setminus \{0\}$, with $|x| < \varepsilon$, such that $\hat{x}^l \in D$. It then suffices to observe further that $\{x \in M^{d_l} : \hat{x}^l \in D\} \subseteq U_l$.

(2). By (†)(c), after permuting coordinates, $D = \Gamma(f)$, where $f : \prod_{i\in[r-1]} M^{d_r} \to M^{d_r}$ is a linear map. Since applying permutation of coordinates on $D$, taking the shell of the resulting set, and then permuting coordinates back results to the shell of $D$, we may assume $D = \Gamma(f)$, where $f : \prod_{i\in[r-1]} M^{d_i} \to M^{d_r}$ (injectivity of $p_l$ on a set is also not affected). Now, by Lemma 2.10, $Sh(D) = \Gamma(f)_{Sh(p_r(D))}$. Thus, $p_{r\restriction Sh(D)}$ is injective, showing (2).

(3). Let $x \in C$. We want to find an open box $B$ containing $x$ such that $B \cap V \subseteq C$. We work by induction on $n$. For $n = 0$, if $D \subseteq M^0$ then $D = Sh(D)$ and we are done. Let $n > 0$. If $\mathsf{dim} D = n$, then $D$ is open, hence we are done. Let $\mathsf{dim} D < n$, then by (†)(c), we may again assume that $D = \Gamma(f)$, where $f : M^{n-1} \to M$ is a linear map. Let $\pi : M^n \to M^{n-1}$ be the projection onto the first $n - 1$ coordinates. By Inductive Hypothesis, since $\pi(Sh(C)) = Sh(\pi(C))$, we obtain that $\pi(C)$ is open in $\pi(Sh(C))$. So let $B_1 \subseteq M^{n-1}$ be an open box containing $\pi(x)$ with $B_1 \cap \pi(C) \subseteq \pi(V)$. Let $B = B_1 \times M$. If $C = \Gamma(f)_{\pi(C)}$ is the graph of a linear map $f$, then $Sh(C) = \Gamma(f)_{\pi(Sh(C))}$ (Lemma 2.10), and we have $B \cap V = \Gamma(f)_{B_1\cap\pi(C)} \subseteq C$, as needed. □

**Lemma 2.13.** *Let $D \subseteq M^n$ be a semilinear set satisfying (†) and containing 0, and $q \in \{0, \dots, r\}$. Then the following are equivalent:*

*(1) for every $l \in [r] \setminus [q]$, $p_{l \restriction Sh(D)}$ is non-injective*
*(2) there are definably connected $D_l \subseteq M^{d_l}$ containing $0$, $l \in [r] \setminus [q]$, such that*

$$\{0\} \times D_{q+1} \times \cdots \times D_r \subseteq D.$$

*Proof.* (2)⇒(1) is clear. For (1)⇒(2), let $V = Sh(D)$. By Lemma 2.12(3), $D$ is open in $V$. Now, for every $l \in [r] \setminus [q]$, since $p_{l \restriction Sh(D)}$ is non-injective, by Fact 2.5 $U_l := U_l^{Sh(D)}$ is infinite, and by Lemma 2.6, $U_1 \times \cdots \times U_r \subseteq V$. Since $D$ is open in $V$, there is an open box $B = B_1 \times \cdots \times B_r$ containing 0, such that $B \cap V \subseteq D$. Therefore

$$B \cap (\{0\} \times U_{q+1} \times \cdots \times U_r) \subseteq D,$$

and hence

$$\{0\} \times (B_{q+1} \cap U_{q+1}) \times \cdots \times (B_r \cap U_r) \subseteq D.$$

Moreover, for $l \in [r] \setminus [q]$, since $p_{l \restriction Sh(D)}$ is non-injective, by Lemma 2.12(2) we obtain $p_{l \restriction D}$ is non-injective. Hence, by Lemma 2.12(1), $B_l \cap U_l$ is infinite. Finally, note that $B \cap V = B \cap D$ is definable, and hence so is $\pi_l(B \cap V)$. Since

$$B_l \cap U_l = \pi_l(B \cap (\{0\} \times U_l \times \{0\})) \subseteq \pi_l(B \cap V) \subseteq \pi_l(D),$$

we may set $D_l = \pi_l(B \cap V)$. It is easy to see that $B \cap V$ (and hence $D_l$) is definably connected. □

We finish this subsection by recalling some further background for a concrete o-minimal linear setting. While $\langle \mathbb{R}, <, + \rangle$ in Section 5 is our main targeted structure (as an ordered vector space over $\mathbb{Q}$), many of our results hold in the following $\mathcal{M}$. The key point is that a 'linear cell decomposition theorem' is known for this $\mathcal{M}$.

*Let $\mathcal{M} = \langle M, <, +, \{x \mapsto \lambda x\}_{\lambda \in \Lambda} \rangle$ be an ordered vector space over an ordered division ring $\Lambda$.*

In this setting, a definable function $f : X \subseteq M^n \to M$, $n > 0$, is in $L(X)$ if and only if it has form

$$f(x_1, \ldots, x_n) = \lambda_1 x_1 + \cdots + \lambda_n x_n + a,$$

where $\lambda_i \in \Lambda$ and $a \in M$ (this can be seen either directly, or using the linear cell decomposition theorem from [25, Section 3], also stated below). Clearly, if $f \in L(X)$ then it extends uniquely to a map in $L(\overline{X})$, where is $\overline{X}$ is the topological closure of $X$, and hence we write $f(a)$ for its value at $a$, even if $a \in \overline{X}$. It also extends to a (non-necessarily unique) linear map in $L(M^n)$.

A *linear cell* in $M^n$ is defined similarly to [17, Chapter 3, (2.3)], recursively, as follows. A linear cell in $M^0$ is just $M^0$. A linear cell in $M^{n+1}$ with $n \geq 0$ is either a graph $\Gamma(\alpha)$, for some $\alpha \in L(X)$, or a cylinder $(\alpha, \beta)$, for some $\alpha, \beta \in L_\infty(X)$, and $X \subseteq M^n$ a linear cell. By the *linear cell decomposition theorem* ([25, Section 3]), every semilinear set is a finite union of linear cells, and moreover, definable functions are piecewise linear (although we will not use the latter statement in this paper).

2.6. **Presburger arithmetic.** Let $\mathcal{M} = \langle M, <, + \rangle \models \mathsf{Pres}$ be a model of Presburger arithmetic. The equivalence relation $\equiv \bmod n$, for $n \in \mathbb{N}_{>0}$, is definable. Recall that $\mathcal{M}$ is an elementary extension of $\langle \mathbb{Z}, <, + \rangle$, and that definable sets are also called Presburger sets.

In this setting, a definable function $f : X \subseteq M^n \to M$, $n > 0$, is in $L(X)$ if and only if it has form

$$f(x) = \sum_{i=1}^{n} a_i \left( \frac{x_i - c_i}{r_i} \right) + \gamma,$$

where $\gamma \in M$, $a_i$ and $0 \leq c_i < r_i$ are integers, for $i = 1, \ldots, m$, and for every $x \in X$, each $x_i \equiv c_i \operatorname{mod} r_i$. (This again can be seen either directly, or using the Presburger cell decomposition theorem from Fact 2.15 below). We call such a map Pres-linear.

Given Pres-linear maps $\alpha, \beta : X \subseteq M^n \to M$ and integers $0 \leq c < r$, we denote by $[\alpha, \beta]_r^c$ the Pres-*cylinder*

$$[\alpha, \beta]_r^c = \{(x, t) \in M^{n+1} : x \in D, \alpha(x) \leq t \leq \beta(x), t \equiv c \operatorname{mod} r\}.$$

We define $[\alpha, \infty)_r^c$ by removing condition '$\alpha(x) \leq t$' in the above definition, and similarly for $(-\infty, \beta]_r^c$ and $(-\infty, \infty)_r^c$. If $n = 0$, we call a Pres-cylinder a Pres-*interval.*

**Definition 2.14** ([13])**.** We define *Presburger cells* recursively as follows. A Presburger cell in $M^0$ is just $M^0$. A Presburger cell $C$ in $M^{n+1}$ with $n > 0$ is either

(1) a graph $\Gamma(\alpha)$, where $\alpha : X \subseteq M^n \to M$ is a linear map, or
(2) a Pres-cylinder $[\alpha, \beta]_r^c, [\alpha, \infty)_r^c, (-\infty, \beta]_r^c, (-\infty, \infty)_r^c$, for some linear maps $\alpha, \beta : X \subseteq M^n \to M$, such that there is no $N \in \mathbb{N}$ satisfying that for all $x \in X$, $|C_x| < N$, where $X$ is a Presburger cell. By *dimension* of $X$ we mean the number of times we apply the second step in its construction.

**Fact 2.15** (Presburger cell decomposition [13])**.** *Every $A$-definable set can be partitioned into finitely many $A$-definable cells, and $A$-definable functions are piecewise $A$-definable and linear.*

Given a definable set $X$, we define

$$\mathsf{dim} X = \max\{\mathsf{dim} C : C \subseteq X \text{ is a Presburger cell}\}.$$

In this paper we use 'dim' both in the o-minimal and Presburger settings, and it is clear from the context which one it is.

2.7. **Property $(*)$ and Reduction Strategy to Theorem 2.7.** We now return to Theorem 2.7 and explain how it can be utilised in our geometric settings. Let $\langle M, + \rangle$ be a group. Consider the following properties for a set $E \subseteq \prod_{i \in [r]} M^{d_i}$:

$(*)$ for every $k \in \mathbb{N}$, if $E$ is $k$-free, then $Sh(E)$ is $k$-free.

$(*)_\infty$ if $E$ is $\infty$-free, then $Sh(E)$ is $\infty$-free.

These properties capture some notion of 'linearity', as the following example shows.

**Example 2.16.** Consider the binary relations $\Delta = \{(x, x) \in \mathbb{R}^2 : x \in (0, 1)\}$ and $E = \{(x, x^2) \in \mathbb{R}^2 : x \in (0, 1)\}$. They are both 2-free, $Sh(\Delta) = \{x, x) \in \mathbb{R}^2 : x \in \mathbb{R}\}$ is 2-free, but $Sh(E) = \mathbb{R}^2$ is not.

**Proposition 2.17.** *Let $E \subseteq \prod_{i \in [r]} M^{d_i}$. Suppose $(*)$ (respectively, $(*)_\infty$) holds. Then $Zar(E)$ (respectively, $Zar^\infty(E)$) holds. If moreover there is an order $<$ such that $\langle M, <, + \rangle$ is an o-minimal ordered group or a model of* Pres*, we may replace $Zar(E)$ by $Zar_1(E)$, and $Zar^\infty(E)$ by $Zar_1^\infty(E)$.*

*Proof.* Suppose $E$ is $k$-free ($\infty$-free). By $(*)$ (respectively, $(*)_\infty$), so is its shell $Sh(E)$. By Theorem 2.7(1)⇒(3), $Sh(E)$ has linear $Z$-bounds for the class of all grids (witnessed by 1 in the 'moreover' cases). Hence so does $E \subseteq Sh(E)$, as needed. □

**Corollary 2.18.** *Let $\mathcal{M} = \langle M, +, \dots \rangle$ be an expansion of a group, such that every definable set is a finite union of sets each satisfying $(*)$ (respectively, $(*)_\infty$) after perhaps a translation. Then for every definable set $E$, $Zar(E)$ (respectively, $Zar^\infty(E)$).*

*Proof.* By Remark 2.2 and Proposition 2.17. A translation $f$ does not affect the validity of the Zarankiewicz statements since for any finite set $B$, $B \sim f(B)$. □

As mentioned earlier, we will reduce the Zarankiewicz statements in each of our geometric settings to Theorem 2.7 by virtue of Property $(*)$. The reduction is manifested in diverse ways, through variants of Corollary 2.18.

**Reduction Strategy:**

(1) In Presburger arithmetic and ordered vector spaces, $(*)$ or $(*)_\infty$ hold: if a Presburger/linear cell $C$ containing 0 is $k$-free ($k \in \mathbb{N}$) or $\infty$-free, respectively, then so is its shell $Sh(C)$. (Proposition 4.34 and Corollary 2.21.)

(2) In semibounded structures, if a 'cone' $E = V + D$ is $\mathcal{C}$-$\infty$-free for the class $\mathcal{C}$ of 'sufficiently distant' grids, then the shell $Sh(V)$ of its 'linear skeleton' $V$ is $\infty$-free ($D$ is bounded). Note that we cannot expect the same for $Sh(E)$, see Remark 3.11. This implies that for some $l$, the projection $p_{l\restriction E\cap B}$ is injective for every $B \in \mathcal{C}$ (Proposition 3.10).

(3) In $\langle \mathbb{R}, <, +, \mathbb{Z} \rangle$, if the sum $E = S + J + D$ of a Presburger cell $S$ with a 'product linear cell' $J + D$ is $k$-free, and $0 \in S \cap J$, then so is $Sh(S + J) + D$, where $S + J$ is the 'purely unbounded part' of $E$ ($D$ is bounded) (Lemma 5.5). This implies that some projection $p_{l\restriction E}$ is uniformly finite-to-1 (Proposition 5.6). In fact, we prove the above statement with $f(S)$ in place of $S$, where $f$ is a linear map.

2.8. **Linear Zarankiewicz.** Here, we provide a direct proof of Fact 1.3 for linear o-minimal structures (Proposition 2.24). Let first $\mathcal{M} = \langle M, <, +, \{x \mapsto \lambda x\}_{\lambda \in \Lambda} \rangle$ be an ordered vector space over an ordered division ring $\Lambda$. Our goal is to prove that every linear cell containing the origin satisfies $(*)_\infty$. We begin by observing this statement in a slightly more general setting.

**Proposition 2.19.** *Let $C$ be a semilinear set that contains 0 and satisfies (†) (from Section 2.5). Then it satisfies $(*)_\infty$. In fact, if $C$ is $\infty$-free, then $Sh(C)$ is 2-free. In particular, $Zar_1^\infty(C)$ holds.*

*Proof.* If $Sh(C)$ is not 2-free, then all $p_{l\restriction Sh(C)}$ are not injective. It follows from Lemma 2.13 that $C$ is not $\infty$-free. We conclude by Proposition 2.17. □

We now turn to linear cells.

**Fact 2.20.** *Let $C \subseteq \prod_{i=1}^r M^{d_i} = M^n$ be a linear cell containing 0. Then (†) holds.*

*Proof.* By a straightforward induction on $n$, left to the reader. □

**Corollary 2.21.** *Every linear cell $C$ that contains 0 satisfies $(*)_\infty$. In fact, if $C$ is $\infty$-free, then $Sh(C)$ is 2-free. In particular, $Zar_1^\infty(C)$ holds.*

*Proof.* By Proposition 2.19 and Fact 2.20. □

**Corollary 2.22.** *Let $E \subseteq \prod_{i=1}^r M^{d_i}$ be a semilinear set. Then $Zar^\infty(E)$.*

*Proof.* By linear cell decomposition, and Corollaries 2.21 and 2.18. □

*Remark* 2.23. We note the following:

(1) The parametric version of linear Zarankiewicz (Fact 1.3(2)) also follows. Indeed, by linear cell decomposition, there is $N \in \mathbb{N}$ such that each $E_b$ is a union of at most $N$ linear cells $D$, and for each of them $Zar_1^\infty(D)$ holds, by Corollary 2.21. By Remark 2.2, we are done.

(2) If $\mathsf{dim} E = s$ and $E$ is $\infty$-free, then $|E \cap B| = O(n^s)$, for finite grids $B$.[2] Indeed, by Remark 2.2, we may assume $E$ is a linear cell. i If $s > r-1$, then since $|E \cap B| = O(n^{r-1})$, the result follows. If $s \leq r-1$, let $p : M^n \to M^s$ be some coordinate projection that is injective on $C$. Then $p(C)$ is contained in some product of $q$ many blocks $M^{d_i}$, with $q \leq s$ ($q = s$ if all $d_i = 1$). But then $|E \cap B| = O(n^q)$, as needed.

We can now conclude the main result of this subsection.

**Proposition 2.24.** *Let $\mathcal{M} = \langle M, <, +, \dots \rangle$ be a linear o-minimal structure. Let $\mathcal{E} = \{E_b\}_{b \in I}$ be a definable family of relations $E_b \subseteq \prod_{i \in [r]} M^{d_i}$. Then $Zar(\mathcal{E})$ holds.*

*Proof.* By [35], $\mathcal{M}$ can be elementarily embedded into a reduct of an ordered vector space $\mathcal{N} = \langle N, +, <, 0, \{\lambda\}_{\lambda \in \Lambda} \rangle$ over an ordered division ring $\Lambda$. Since for a fixed $k \in \mathbb{N}$, $k$-freeness is a first-order property, $Zar(E)$ is preserved under elementary equivalence and taking reducts. Hence it holds for $\mathcal{M}$. □

Note that Fact 1.3 is still more general than Proposition 2.24, since it does not assume an ambient group structure (and we would also need [4, Corollary 5.8] to obtain $Zar^\infty(\mathcal{E})$ from it).

## 3. Semibounded Zarankiewicz

*In the section, $\mathcal{M} = \langle M, <, +, 0, \dots \rangle$ denotes a semibounded o-minimal expansion of an ordered group which is not linear. We fix an element $1 > 0$ such that a real closed field, with universe $(0, 1)$ and whose order agrees with $<$, is definable in $\mathcal{M}$.*

Our goal is to prove Theorems A′ and A″ via Reduction Strategy (2) of Section 2.7 (Proposition 3.10). We begin with some basics of semibounded structures. Following [44], an interval $I \subseteq M$ is called *short* if there is a definable bijection between $I$ and $(0, 1)$; otherwise, it is called *long*. Equivalently, an interval $I \subseteq M$ is short if a real closed field whose domain is $I$ is definable. An element $a \in M$ is called *short* if either $a = 0$ or $(0, |a|)$ is a short interval; otherwise, it is called *tall*. A tuple $a \in M^n$ is called *short* if $|a| := |a_1| + \cdots + |a_n|$ is short, and *tall* otherwise. A definable set $X \subseteq M^n$ (or its defining formula) is called *short* if it is in definable bijection with a subset of $(0, 1)^n$; otherwise, it is called *long*. Notice that this is compatible, for $n = 1$, with the notion of a short interval.

[2] We thank S. Starchenko for suggesting this statement.

3.1. **Cones.** We adopt the definition of a cone from [21] which is a refinement of that from [20]. We first fix some standard terminology and notation.

Recall from Section 2.5 that $\Lambda$ is the set of all $\emptyset$-definable partial endomorphisms of $\langle M, <, +, 0\rangle$. For $v \in \Lambda$, we denote by $\mathrm{dom}(v)$ and $\mathrm{ran}(v)$ the domain and range of $v$, respectively. We write $vt$ for $v(t)$. It is a standard practice in this section that whenever we write an expression of the form '$vt$', with $v \in \Lambda$ and $t \in M$, we mean in particular that $t \in \mathrm{dom}(v)$. Sometimes, however, we say explicitly that $t \in \mathrm{dom}(v)$. As in Section 2.4, a partial endomorphism $v : (0, a) \to M$, $a \in M$, can be (uniquely) extended to a partial endomorphism $f : \bigcup_{n\in\mathbb{N}}(-na, na) \to M$. We keep the notation $v, \mathrm{dom}(v), \mathrm{ran}(v)$ for this extension.

Following [44], we say that two $\lambda, \mu \in \Lambda$ *have the same germ at* 0, denoted $\sim_g$ if there is $\epsilon > 0$, such that the restrictions of $\lambda, \mu$ on $(-\epsilon, \epsilon)$ are the same. It is observed in [44, Section 6], that $\Lambda$ modulo $\sim_g$ can be given the structure of an ordered field with multiplication given by composition. For a matrix $A = (a_{ij})$ with entries from $\Lambda$, the *rank* of $A$ is the rank of the matrix $\overline{A} = (\bar{a}_{ij})$, where $\bar{a}_{ij}$ is the $\sim_g$-equivalence class of $a_{ij}$. As explained in [21, Section 2.1], it is a routine to check that various classical results from linear algebra hold for matrices with entries from $\Lambda$. For example, a $l \times k$ linear system with coefficients from $\Lambda$ has a unique solution if and only if the coefficient matrix has rank $k$. We freely use such results.

We now proceed to the notion of $\Lambda$-independence, which is needed to define cones.

**Definition 3.1.** If $v = (\lambda_1, \dots, \lambda_n) \in \Lambda^n$ and $\mu \in \Lambda$, we denote $\mu v := (\mu\lambda_1, \dots, \mu\lambda_n)$. We say that $v_1, \dots, v_k \in \Lambda^n$ are $\Lambda$-*independent* if for all $\mu_1, \dots, \mu_k$ in $\Lambda$, with $\mathrm{ran}(v_i) \subseteq \mathrm{dom}(\mu_i)$,

$$\mu_1 v_1 + \cdots + \mu_k v_k = 0 \text{ implies } \mu_1 = \cdots = \mu_k = 0.$$

If $v = (\lambda_1, \dots, \lambda_n) \in \Lambda^n$ and $t \in M$, we denote $vt := (\lambda_1 t, \dots, \lambda_n t)$ and $\mathrm{dom}(v) := \bigcap_{i=1}^n \mathrm{dom}(\lambda_i)$.

We are now ready to present the notion of a cone.

**Definition 3.2.** Let $k \in \mathbb{N}$. A *k-cone* $C \subseteq M^n$ is a definable set of the form

$$C = \left\{ b + \sum_{i=1}^{k} v_i t_i : b \in S,\ t_i \in J_i \right\},$$

where $S \subseteq M^n$ is a short cell, $v_1, \dots, v_k \in \Lambda^n$ are $\Lambda$-independent and $J_1, \dots, J_k$ are long intervals each of the form $(0, a_i)$, $a_i \in M_{>0} \cup \{\infty\}$, with $J_i \subseteq \mathrm{dom}(v_i)$. So a 0-cone is just a short cell. A *cone* is a $k$-cone, for some $k \in \mathbb{N}$. We say that the long cone $C$ is *normalised* if for each $x \in C$ there are unique $b \in B$ and $t_1 \in J_1, \dots, t_k \in J_k$ such that $x = b + \sum_{i=1}^k v_i t_i$. In this case, we write:

$$C = S + \sum_{i=1}^{k} v_i t_i | J_i \subseteq M^n.$$

In what follows, all long cones are assumed to be normalised, and we thus drop the word 'normalised'. The *linear skeleton* of $C$ is the set

$$\langle C\rangle = \left\{ \sum_{i=1}^{k} v_i t_i : t_i \in \bigcup_{n\in\mathbb{N}} (-na_i, na_i) \right\}.$$

*Remark* 3.3. We note the following:

(1) Our notion of a $k$-cone is the same with that of a $k$-long cone from [21, Definition 2.9], since '$M$-independence' mentioned there is equivalent to $\Lambda$-independence ([21, Lemma 2.4]).
(2) $\langle C\rangle$ is a subgroup of $\langle M^n,+\rangle$ that need not be definable. Our notation differs from that in [21], where $\langle C\rangle = \left\{\sum_{i=1}^k v_i t_i : t_i \in (-a_i, a_i)\right\}$, as well as from that in [20], where $\langle C\rangle = \left\{\sum_{i=1}^k v_i t_i : t_i \in M\right\}$ (which are both definable).
(3) In case $M$ does not contain any tall elements (such as when $M=\mathbb{R}$), in our definition of a cone $C$ we have each $a_i=\infty$, and hence our notion of a cone $C$ and the notation $\langle C\rangle$ coincide with those in [20].

3.2. **Projections of cones intersected with grids.** The goal of this subsection is to prove Proposition 3.7 below, saying that projection maps restricted to intersections of a cone $E$ with sufficiently distant grids behave similarly to their restriction on the linear skeleton $\langle E\rangle$ of $E$. 'Sufficiently distant' here means precisely $m$-distant, where $m=m(E)$ is defined right before Lemma 3.6.

Let us fix a $k$-cone

$$E = S + \sum_{i=1}^k v_i t_i | J_i \subseteq M^n,$$

where $J_i=(0,a_i)$ is a long interval, for some $a_i \in M_{>0}\cup\{\infty\}$. Since the vectors $v_i = \begin{pmatrix} v_i^1 \\ \vdots \\ v_i^n \end{pmatrix}$, $i=1,\dots,k$, are $\Lambda$-independent, the matrix

$$A = \begin{pmatrix} v_1^1 & \dots & v_k^1 \\ \vdots & \dots & \vdots \\ v_1^n & \dots & v_k^n \end{pmatrix}$$

has rank $k$. For $j=1,\dots,n$, denote by $A_j$ the $j$-th row of $A$.

The following lemma is straightforward from elementary linear algebra, but we prove it for completeness.

**Lemma 3.4.** *Fix any $l\le n$ coordinates of $M^n$. Let $\pi : M^n \to M^l$ be the projection map onto those coordinates. Let also $D$ be the $l\times k$ submatrix of $A$ consisting of the corresponding $l$ rows. Then the following are equivalent:*

*(1) $\pi_{\restriction \langle E\rangle}$ is injective.*
*(2) the system $Dx=0$ has a unique solution.*

*Proof.* Without loss of generality, we may assume that the fixed $l$ coordinates are the first $l$ ones.

(2)⇒(1): Let $c, c' \in \langle E\rangle$ be distinct with $\pi(c)=\pi(c')$. Suppose

$$c = \sum_{i=1}^k v_i t_i \quad \text{and} \quad c' = \sum_{i=1}^k v_i t_i',$$

for some $t_i, t_i' \in \mathrm{dom}(v_i)$. So we have

$$\pi\left(\sum_{i=1}^k v_i t_i\right) = \pi\left(\sum_{i=1}^k v_i t_i'\right).$$

Then $D\begin{pmatrix} t_1-t'_1 \\ \vdots \\ t_k-t'_k \end{pmatrix}=0$, and hence $Dx=0$ has non-trivial solutions.

(1)⇒(2) Let $s=(s_1,\dots,s_k)\neq 0$ be a solution to $Dx=0$. Take any $c=\sum_{i=1}^k v_i t_i \in \langle E\rangle$. Then for $t'_i=t_i+s_i$, $i=1,\dots,k$, and $c'=\sum_{i=1}^k v_i t'_i$, we have $c\neq c'$ and $\pi(c)=\pi(c')$. Hence, $\pi_{\restriction\langle E\rangle}$ is not injective. □

We now show how condition (1) of the previous lemma applied to a cone $E$, further implies this: for a suitable $m=m(E)\in M_{\geq 0}$ the restriction of $\pi$ on $m$-distant sets intersected with $E$ is injective.

**Definition of $m(E)$.** Let $m(E)\in M_{\geq 0}$ be the supremum of the absolute values of all elements of the form

$$(b_j-b'_j)+A_jD^{-1}\begin{pmatrix} b_{i_1}-b'_{i_1} \\ \vdots \\ b_{i_k}-b'_{i_k} \end{pmatrix},$$

where

- $D$ ranges over all $k\times k$ submatrices of $A$, $j,i_1,\dots,i_k\in\{1,\dots,n\}$, and
- $b_l,b'_l$ denote projections of elements $b,b'$ onto the $l$-th coordinate, as $b,b'$ range over $S$.

*Remark* 3.5. The element $m(E)$ is short. Indeed, since $S$ is a short connected set, all expressions of the form $b-b'$ with $b,b'\in S$, are short. Moreover, $\lambda s$ is short whenever $\lambda\in\Lambda$ and $s\in M$ is short. It follows that $m(E)$ is short.

**Lemma 3.6.** *Let $l,\pi,D$ be as in Lemma 3.4, and $m=m(E)$. Let also $B\subseteq M^n$ be an $m$-distant set. Suppose $\pi_{\restriction\langle E\rangle}$ is injective. Then $\pi_{\restriction E\cap B}$ is injective.*

*Proof.* Without loss of generality, we may again assume that the fixed $l$ coordinates are the first $l$ ones. By Lemma 3.4, if $D$ is the matrix consisting of the first $l$ rows of $A$, then $Dx=0$ has a unique solution. Equivalently, $\mathsf{rank}(D)=k$. In particular, $l\geq k$. Let $D'$ be a $k\times k$ submatrix of $D$ of rank $k$. Without loss of generality, $D'$ consists of the first $k$ rows of $D$. It is enough to prove that for the projection $\pi':M^n\to M^k$ onto the first $k$ coordinates, and every $m$-distant set $B$, $\pi'_{\restriction E\cap B}$ is injective. In other words, we may assume that $l=k$.

Let $B$ be an $m$-distant grid and denote $C=E\cap B$. We prove that $\pi_{\restriction C}$ is injective. Let $c,c'\in C$ with $\pi(c)=\pi(c')$. Suppose

$$c=b+\sum_{i=1}^k v_i t_i \quad\text{and}\quad c'=b'+\sum_{i=1}^k v_i t'_i,$$

for some $b,b'\in S$ and $t_i,t'_i\in\mathrm{dom}(v_i)$. So we have

$$\pi\left(b+\sum_{i=1}^k v_i t_i\right)=\pi\left(b'+\sum_{i=1}^k v_i t'_i\right).$$

Then $D\begin{pmatrix} t_1-t'_1 \\ \vdots \\ t_k-t'_k \end{pmatrix}=\begin{pmatrix} b'_1-b_1 \\ \vdots \\ b'_k-b_k \end{pmatrix}$. Hence $\begin{pmatrix} t_1-t'_1 \\ \vdots \\ t_k-t'_k \end{pmatrix}=D^{-1}\begin{pmatrix} b'_1-b_1 \\ \vdots \\ b'_k-b_k \end{pmatrix}$ and each row of the last matrix consists of an element in $S-S'$. It follows that for every $j=1,\dots,n$,

we have

$$c_j - c'_j = (b_j - b'_j) + A_j \begin{pmatrix} t_1 - t'_1 \\ \vdots \\ t_k - t'_k \end{pmatrix} = (b_j - b'_j) + A_j D^{-1} \begin{pmatrix} b'_1 - b_1 \\ \vdots \\ b'_k - b_k \end{pmatrix}.$$

By the choice of $m$, we have that $|c_j - c'_j| \leq m$. By the assumption on $B$, it follows that $c = c'$. ☐

The converse of Lemma 3.6 also holds, but we will not be using it here.

**Proposition 3.7.** *Let $E \subseteq \prod_{i=1}^{r} M^{d_i} = M^n$ be a $k$-cone and $m = m(E)$. Let $\pi : M^n \to M^l$ be a projection onto some $l$ coordinates. If $\pi_{\restriction \langle E \rangle}$ is injective, then for every $m$-distant grid $B$, we have that $\pi_{\restriction E \cap B}$ is injective.*

*Proof.* Note that an $m$-distant grid is in particular an $m$-distant set, and hence it satisfies the assumption of Lemma 3.6. ☐

3.3. **Reduction strategy.** In this section, we complete Reduction Strategy (2) from Section 2.7 (Proposition 3.10). We first point out a minor correction from [21]. As stated, [21, Corollary 2.14] is not correct. There are two cases: (I) $M$ contains tall elements (such as in a saturated setting), (II) $M$ does not contain any tall elements (such as when $M = \mathbb{R}$). The proof of [21, Corollary 2.14] assumes we are in Case (I), and this assumption should be added to the statement. In Case (II), a version of [21, Corollary 2.14] remains true (namely, [20, Lemma 3.2]). These two cases are now incorporated correctly in Fact 3.8 below. We note that the rest of the results in [21] in Case (II) remain unaffected. In fact, in Case (II) the main results of [21] were already known from [20], since in this case, the notions of cones in [20] and [21] (and here) coincide (Remark 3.3(3)).

Let us denote

$$\widetilde{V} = \left\{ \sum_{i=1}^{k} v_i t_i : t_i \in (-a_i, a_i) \right\}.$$

**Fact 3.8.** *Let $\lambda \in \Lambda^d$ and $t \in M_{>0}$ with $\lambda t \in \widetilde{V}$. Then either*

- *there is a tall $s \in M_{>0}$ with $\lambda s \in \widetilde{V}$, or*
- *$\{\lambda s : s \in M_{>0}\} \subseteq \widetilde{V}$.*

*Proof.* Case (I): $M$ contains tall elements. Then this is by [21, Corollary 2.14].

Case (II): $M$ does not contain any tall elements. In this case, $\widetilde{V} = \langle v_1, \ldots, v_k \rangle$ in the notation of [20, Lemma 3.2](1). Moreover, our assumption implies the assumption of that lemma, namely that $\lambda t \in \langle v_1, \ldots, v_k \rangle$ with $t$ positive. Hence its conclusion holds, giving that

$$\{\lambda s : s \in M_{>0}\} \subseteq \langle v_1, \ldots, v_k \rangle = \widetilde{V},$$

as needed. ☐

We will also need the following fact.

**Fact 3.9.** *The set $\widetilde{V} \subseteq \prod_{i=1}^{r} M^{d_i} = M^n$ contains 0 and satisfies (†) from Section 2.5. Hence the equivalence of Lemma 2.13 holds for $D = \widetilde{V}$.*

*Proof.* Easy, left to the reader. ☐

Recall that for $m \in M_{\geq 0}$, $\mathcal{C}_m$ denotes the class of all $m$-distant grids.

**Proposition 3.10.** *Let $E$ be a cone, $V = \langle E\rangle$, and $m = m(E)$. Then the following are equivalent:*

*(1) $E$ is $\mathcal{C}_m$-$\infty$-free.*
*(2) $V$ is 2-free.*
*(3) there is $l \in [r]$, such that $p_{l\restriction V}$ is injective.*
*(4) there is $l \in [r]$, such that for all $B \in \mathcal{C}_m$, $p_{l\restriction E\cap B}$ is injective.*
*(5) $E$ has linear $Z$-bounds for $\mathcal{C}_m$, witnessed by $\alpha = 1$.*

*In particular, $Zar_1^\infty(E, \mathcal{C}_m)$.*

*Proof.* (1)⇒(2). Suppose

$$E = S + \sum_{i=1}^{k} v_i t_i | J_i \subseteq M^n,$$

and observe that $V = Sh(\widetilde{V})$. If $V$ is not 2-free, then some projection $p_{l\restriction V}$ is not injective, and hence by Fact 3.9 and Lemma 2.13, there are infinite definable $D_l \subseteq M^{d_l}$ containing 0, for $l \in [r]$, such that

$$D_1 \times \cdots \times D_r \subseteq \widetilde{V}.$$

Since each $D_l$ is definably connected, there is a definable curve inside it approaching 0; namely, there is a definable map $\gamma_l : (0,p) \to D_l$ with $\lim_{t\to 0^+} \gamma_l(t) = 0$. By piecewise linearity, $\gamma_l$ is eventually linear, that is of the form $\lambda_l t$, where $\lambda_l \in \Lambda^{d_l} \setminus \{0\}$. We obtain that there is some $t \in M_{>0}$ such that $\lambda_l t \in D_l$. Thus $\{0\} \times \{\lambda t\} \times \{0\} \in \widetilde{V}$. Hence, Fact 3.8 applies, with $\lambda = (0,\dots,0,\lambda_l,0,\dots,0) \in \Lambda^n \setminus \{0\}$.

Suppose first that there is a tall element $s \in M_{>0}$ such that $\lambda s \in \widetilde{V}$. Since $\lambda_l \in \Lambda^n \setminus \{0\}$, one of the coordinates of $g_l = \lambda_l s$ is tall, and hence $|g_l| > m$ (recall that $m$ is short).

Suppose now that

$$\{\lambda s : s \in M_{>0}\} \subseteq \widetilde{V}.$$

Then clearly we can find large enough $s \in M_{>0}$ so that for $g_l = \lambda_l s$, we have $|g_l| > m$.

In both cases, for every $l \in [r]$, we found $g_l \in D_l$ with $|g_l| > m$. Therefore, since $V = Sh(\widetilde{V})$,

$$\mathbb{Z}g_1 \times \cdots \times \mathbb{Z}g_r \subseteq V,$$

showing that $V$ (and hence $E$) is not $\mathcal{C}_m$-$\infty$-free.

(2)⇒(3). By Theorem 2.7(1)⇒(2).

(3)⇒(4). By Proposition 3.7.

(4)⇒(5). By Lemma 2.1.

(5)⇒(1). Immediate from the definitions.

The 'in particular' clause is by the fact that (1)⇒(5). □

*Remark* 3.11. It is not possible to replace (2) by '$Sh(E)$ is $\infty$-free'. For example, for the cylinder $E = S^1 + V$ from Section 1.2 with all $d_i = 1$, $E$ is $\mathcal{C}_1$-$\infty$-free, but $Sh(E) = M^3$ is not.

We can now conclude the proof of Theorem A′.

**Corollary 3.12.** *Let $R = (-\gamma, \gamma) \subseteq M$ be an open interval, $\gamma \in R \cup \{\infty\}$, such that there is no definable field with domain $R$. Then for every definable $E \subseteq \prod_{i\in[r]} M^{d_i}$, there is $m \in R$, such that $Zar(E, \mathcal{C}_m)$ holds.*

*Proof.* By Remark 2.2, we may assume that $E$ is a cone. Let $m = m(E)$. By Proposition 3.10, $Zar^\infty(E, \mathcal{C}_m)$ holds. By Remark 3.5, $m$ is a short element; that is, there is a definable real closed field with domain $(0, m)$. Hence $m \in R$. ☐

*Remark* 3.13. A uniform version of Corollary 3.12 also holds. More precisely:

*Let $\{E_b\}_{b\in I}$ be a definable family of relations $E_b \subseteq \prod_{i\in[r]} M^{d_i}$. For each $b \in I$, let $m_b = m(E_b) \in R$ and $\mathcal{C}_{m_b}$ be the class of all finite $m_b$-distant grids. Then there is $\alpha \in \mathbb{R}_{>0}$, such that for every $b \in I$, $Zar_\alpha^\infty(E_b, \mathcal{C}_{m_b})$.*

Indeed, suppose $\{E_b\}_{b\in I}$ is $A$-definable. By the cone decomposition theorem, each $E_b$ is a finite union of $Ab$-definable cones, and by [1, Lemma 4.10], each of $S$ and $J_i$ from Definition 3.2 appearing in each such cone is also $Ab$-definable. By a standard compactness argument (working in a saturated elementary extension of $\mathcal{M}$), we obtain the following uniform cone decomposition for $\{E_b\}_{b\in I}$: there are definable families $\mathcal{D}_1, \dots, \mathcal{D}_l$, with $\mathcal{D}_j = \{C_{j,b}\}_{b\in I}$, where each $C_{j,b}$ is either a cone or an empty set, and for each $b \in I$,

$$E_b = \bigcup_{j=1}^{l} C_{j,b}.$$

For each $b \in I$, let $m_b = \max_j\{m(C_{j,b})\}$. Since, by Proposition 3.10, we have $Zar_1(C_{j,b}, \mathcal{C}_{m_b})$, it follows that $Zar_l(E_b, m_b)$, as needed.

*Remark* 3.14. We note that our semibounded Zarankiewicz (Theorem A′) yields also a stronger linear Zarankiewicz, since there are semilinear sets that are not $\infty$-free, but are $\mathcal{C}_m$-$\infty$-free for some $m \in M_{\geq 0}$.

3.4. **Recognising the full multiplication.** We begin with an 'unbounded version' of the Szemerédi-Trotter theorem, which is useful in the semibounded context.

**Proposition 3.15** (Unbounded Szemerédi-Trotter)**.** *Let $\langle R, <, +, \cdot\rangle$ be an ordered field. Let $E \subseteq R^2 \times R^2$ be the binary relation given by*

$$E(x, y, a, b) \;\Leftrightarrow\; y = ax + b,$$

*saying that $(x, y)$ lies on the line $y = ax + b$. Let $m \in R_{\geq 0}$. Then for finite $m$-distant $n$-grids $B \subseteq R^2 \times R^2$, we have*

$$|E \cap B| = \Omega(n^{4/3}).$$

*Proof.* The classical Szemerédi-Trotter lower bound theorem guarantees (say with the proof appearing in [38, Proposition 4.2.1]) that for arbitrarily large $n$ there are $n$ points and $n$ lines on the plane, with at least $\frac{1}{4}n^{4/3}$ incidences. That is, there are $n$-grids $B \subseteq R^2 \times R^2$ of arbitrarily large size $n$, such that

$$|E \cap B| \geq \frac{1}{4}n^{4/3}.$$

We check that with a small tweak in the proof we can always find such grids that are moreover $m$-distant. Indeed, let $r, s \in R_{>0}$ with $r, s > m$.

The construction from [38, Proposition 4.2.1] is this: let

$$P = \{(x, y) : x = 0, 1, \dots, k-1, \;\; y = 0, 1, \dots, 4k^2 - 1\}$$

and

$$L = \{(a,b) : a = 0, 1, \dots, 2k-1,\ \ b = 0, 1, \dots, 2k^2 - 1\}.$$

Then for $n = 4k^3$ and $B$ the $n$-grid $B = P \times L$, we have

$$|E \cap B| \geq \frac{1}{4} n^{4/3}.$$

Here, we define

$$P' = rP = \{(x,y) : x = 0, r, \dots, (k-1)r,\ \ y = 0, rs, \dots, (4k^2-1)rs\}$$

and

$$L' = \{(a,b) : a = 0, s, \dots, (2k-1)s,\ \ b = 0, r, \dots, (2k^2-1)r\}.$$

We have, for any $x, y, a, b \in R$,

$$y = ax + b \Leftrightarrow ysr = asxr + bsr$$

that is,

$$E(x,y,a,b) \Leftrightarrow E(xr, ysr, as, br),$$

and since the map $(x,y,a,b) \mapsto (xr, ysr, as, br)$ is a bijection, we get

$$|E \cap (P \times L)| = |E \cap (P' \times L')|.$$

It remains to check that each of $P'$ and $L'$ is an $m$-distant set. The two sets handled analogously, we only look at $P'$: for any two distinct $(x,y), (x',y')$, either $x \neq x'$ or $y \neq y$. Any two such distinct $x, x'$ differ by $\geq r > m$, and any two such distinct $y \neq y'$ differ by $\geq sr > m$, so in all cases $(x,y), (x',y')$ are $m$-distant. □

We can now conclude the proof of Theorem A$''$.

**Theorem 3.16.** *Let $\mathcal{M} = \langle M, <, +, \dots\rangle$ be a reduct of a real closed field $\mathcal{R} = \langle M, <, +, \cdot\rangle$. Then the following are equivalent:*

*(1) The multiplication $\cdot$ is not definable in $\mathcal{M}$.*
*(2) Let $1 \leq \beta < \frac{4}{3}$ in $\mathbb{R}_{>0}$. For every binary definable $E \subseteq M^{d_1} \times M^{d_2}$, there are $\alpha \in \mathbb{R}_{>0}$ and $m \in R_{\geq 0}$, such that if $E$ is $k$-free, for some $k \in \mathbb{N}$, then for every $m$-distant $n$-grid $B \subseteq M^{d_1} \times M^{d_2}$, we have*

$$|E \cap B| \leq \alpha n^{\beta}.$$

*Proof.* (1) $\Rightarrow$(2): Since $\mathcal{M}$ is a reduct of a real closed field, and is not semibounded, it defines the field's multiplication. This statement follows from [41], where it is shown for reducts of the real field, and the fact that $\mathcal{R}$ is elementarily equivalent to the real field. Since $\mathcal{M}$ is not semibounded, the result follows from Corollary 3.12, for $\beta = 1$.

(2) $\Rightarrow$(1): By Proposition 3.15. □

*Remark* 3.17. The above proof does not go through without assuming that $\mathcal{M}$ is a reduct of a field $\mathcal{R}$. For example, there are non-semibounded o-minimal expansions of $\langle \mathbb{R}, <, +\rangle$ where the real multiplication is not defined ([43, Section 3]).

## 4. $N$-internal points and Presburger Zarankiewicz

In this section, we carry out some analysis of topological nature in any 'definably complete' ordered group $\mathcal{M} = \langle M, <, +\rangle$ that will be employed when implementing Reduction Strategy (1)&(3) from Section 2.7 in Presburger arithmetic and $\langle \mathbb{R}, <, +, \mathbb{Z}\rangle$ (Proposition 4.29 and Lemma 5.5, respectively). An ordered structure is *definably complete* if every bounded definable subset of its universe has a supremum. As a result of our analysis, we obtain Theorem B (Corollary 4.35).

*In the rest of this section, $\mathcal{M} = \langle M, <, +\rangle$ is any definably complete ordered group. We fix an element $\mathbf{1} \in M_{>0}$, and write $\mathbf{Z}$ for the subgroup of $\langle M, +\rangle$ generated by $\mathbf{1}$, and $\mathbf{N} = \mathbf{Z}_{\geq 0}$.* Of course, for every $n \in \mathbb{N}$, $n\mathbf{Z} \subseteq \mathbf{Z}$. Whenever $T = T'\mathbf{1} \in \mathbf{N}$, with $T' \in \mathbb{N}$, and $x \in M$, we write $Tx$ for $T'x$. Clearly, if $x < y$, then $Tx < Ty$.

We next define the notions of a purely $\mathbf{Z}$-unbounded set, an abstract $\mathbf{Z}$-cell, and an $N$-internal point, for $N \in \mathbf{N}$ (Definitions 4.10, 4.1, 4.15). Our goal is to prove that a purely $\mathbf{Z}$-unbounded $\mathbf{Z}$-cell $C$ contains, for every $N \in \mathbf{N}$, an $N$-internal point (Proposition 4.22); that is, a point whose $N$-neighborhood in $Sh(C)$ is contained in $C$. This guarantees 'enough space' to accommodate grids inside $C$ given they exist in $Sh(C)$ (Corollaries 4.27 and 4.33), which is the key to obtaining $(*)$ in Reduction Strategy (1) for $C$ (Proposition 4.34) and hence Presburger Zarankiewicz for a single set (Corollary 4.35). We also achieve Reduction Strategy (3), by applying Corollary 4.27 to purely unbounded linear cells in $\langle \mathbb{R}, <, +\rangle$ (Lemma 5.5).

### 4.1. Purely ($\mathbf{Z}$-)unbounded sets.

The class of purely $\mathbf{Z}$-unbounded sets is defined recursively.

**Definition 4.1.** Let $C \subseteq M^n$. We say that $C$ is *purely* $\mathbf{Z}$*-unbounded* if

(1) $n = 0$, $C = M^0$.
(2) $n > 0$, $C$ is the graph of a function, and $\pi(C)$ is purely $\mathbf{Z}$-unbounded.
(3) $n > 0$, $C$ is not the graph of a function, $\pi(C)$ is purely $\mathbf{Z}$-unbounded, and for every $N \in \mathbf{N}$, there is $x \in C$, such that $C_{\pi(x)} \not\subseteq B_N(x)$.

We also define the class of purely unbounded sets, which will be used in the sequel.

**Definition 4.2.** The class of *purely unbounded* sets is defined in exactly the same way as in Definition 4.1 after replacing $\mathbf{N}$ by $M_{\geq 0}$. Alternatively, if we were to vary $\mathbf{Z}$, a set is purely unbounded if and only if it is purely $\mathbf{Z}$-unbounded for any $\mathbf{Z}$.

*Remark* 4.3. We note the following:

(1) A purely ($\mathbf{Z}$-)unbounded set $C$ can even be a singleton. But if $C$ is not the graph of a map, then its fibers are not uniformly bounded (in $\mathbf{N}$).
(2) A finite union of non-purely ($\mathbf{Z}$-)unbounded sets is not purely ($\mathbf{Z}$-)unbounded.
(3) Every Presburger cell is purely $\mathbb{Z}$-unbounded. We build further on this in Proposition 5.15.
(4) If $\mathbf{Z}$ is cofinal in $M$, then a set is purely $\mathbf{Z}$-unbounded if and only if it is purely unbounded. This is, for example, the case in any o-minimal structure over the reals, as well as in the standard model of $\mathsf{Pres}$ with $\mathbf{Z} = \mathbb{Z}$ (and it is not the case in a saturated model in either setting and any $\mathbf{Z}$).

4.2. **Abstract Z-cells.** We proceed towards the definition of an abstract $\mathbf{Z}$-cell.

**Definition 4.4.** Two elements $x, y \in M^n$ are called $\mathbf{Z}$-*close* if there is $T \in \mathbf{N}$ with $|x - y| < T$. A set $A \subseteq M$ is called $\mathbf{Z}$-*dense* if either it contains no $x, y$ with $|x - y| > \mathbf{1}$, or it contains $\mathbf{Z}$-close $x, y$ with $|x - y| > \mathbf{1}$.

**Example 4.5.** In a model $\mathcal{M}$ of $\mathsf{Pres}$, every cell in $M$ is $\mathbb{Z}$-dense, but if $a > \mathbb{Z}$, then $\mathbb{Z}a$ is not. In a linear o-minimal structure, and our fixed but arbitrary $\mathbf{Z}$, every interval is $\mathbf{Z}$-dense (if infinitesimal with respect to $\mathbf{1}$, it satisfies the former condition of $\mathbf{Z}$-density, otherwise the latter).

**Definition 4.6.** Let $f : D \subseteq M^n \to M$ be a map. We call $f$ *bi-Lipschitz* if there is $n \in \mathbb{N}$, such that for every $x, y \in D$,

$$|f(x) - f(y)| < n|x - y| \quad \text{and} \quad |x - y| < n|f(x) - f(y)|.$$

**Fact 4.7.** *The following hold:*

(1) *A composition of bi-Lipschitz (respectively, linear) maps is bi-Lipschitz (respectively, linear).*
(2) *Let $f : D \subseteq M^n \to M$ be a bi-Lipschitz map. Then for every $x, y \in D$, if $x, y$ are $\mathbf{Z}$-close, then so are $f(x), f(y)$.*

*Proof.* (1) is well-known and straightforward to prove, whereas (2) follows from the definition of a bi-Lipschitz map and the fact that $n\mathbf{Z} \subseteq \mathbf{Z}$ for any $n \in \mathbb{N}$. □

**Fact 4.8.** *Let $\mathcal{M} = \langle M, <, +\rangle$ be an o-minimal ordered group, or $\mathcal{M} \models \mathsf{Pres}$. Then a linear map $f : X \subseteq M^n \to M$ is bi-Lipschitz.*

*Proof.* In the former case, $f$ has form $f(x_1, \dots, x_n) = \lambda_1 x_1 + \dots + \lambda_n x_n + a$, where each $\lambda_i \in \mathbb{Q}^n$ and $a \in M$. In the latter case, $f$ has form $f(x) = \sum_{i=1}^n a_i \left(\frac{x_i - c_i}{r_i}\right) + \gamma$, where $\gamma \in M$, $a_i$ and $0 \leq c_i < r_i$ are integers. In both cases, the result follows. □

*Remark* 4.9. If $\mathcal{M}$ is an ordered vector space over an arbitrary ordered division ring $\Lambda$, linear maps need no longer be bi-Lipschitz. For example, $\Lambda$ may contain $\lambda$ with $\lambda(\mathbf{1}) > \mathbf{N}$.

We can now define recursively the notion of an abstract $\mathbf{Z}$-cell $C$. Along with it, we define the cylinder $\widehat{C}$. The notion of a shell was introduced in Definition 2.8.

**Definition 4.10.** An *abstract open* $\mathbf{Z}$-*cell* or simply *open* $\mathbf{Z}$-*cell* in $M^n$ is defined as follows. An open $\mathbf{Z}$-cell in $M^0$ is just $M^0$. An open $\mathbf{Z}$-cell $C$ in $M^{n+1}$ with $n > 0$ is either:

- a graph $\Gamma(\alpha)$, where $\alpha : X \subseteq M^n \to M$ is a definable bi-Lipschitz linear map, and $\widehat{C} = C$, or
- infinite, contained in a cylinder $\widehat{C} = (\alpha, \beta)$, where $\alpha, \beta : X \subseteq M^n \to M$ are definable bi-Lipschitz linear maps, or $\alpha = -\infty$, $\beta = \infty$, such that for every $x \in X$, $\inf C_x = \alpha(x)$ and $\sup C_x = \beta(x)$, and the following uniformity and density conditions hold, for every $z \in \pi(C)$:

(UC) $C_z = \widehat{C}_z \cap Sh(C)_z$,
($\mathbf{Z}$-DENSE) $C_z$ is $\mathbf{Z}$-dense,
where $X$ is an open $\mathbf{Z}$-cell in $M^n$.

We define an *abstract closed* $\mathbf{Z}$-*cell* or simply a *closed* $\mathbf{Z}$-*cell* in $M^n$ in exactly the same way after replacing $(\alpha, \beta)$ by $[\alpha, \beta]$, and 'open' by 'closed'. In particular, $\widehat{C} = [\alpha, \beta]$. By a $\mathbf{Z}$-*cell* we mean an open or a closed $\mathbf{Z}$-cell. To lighten notation, we write *purely* $\mathbf{Z}$-*unbounded cell* for purely $\mathbf{Z}$-unbounded $\mathbf{Z}$-cell.

*Remark* 4.11. We note the following:

(1) (UC) is independent from perhaps another interesting uniformity condition:

$$\text{for every } z, w \in \pi(C) \text{ and } a < b \in M,\ C_z \cap [a, b] = C_w \cap [a, b].$$

(2) We have $\pi(\widehat{C}) \subseteq \widehat{\pi(C)}$ but the other inclusion need not be true.

Our notion of a $\mathbf{Z}$-cell is more meaningful when $0 \in C$, as indicated for example by the second item of the next lemma.

**Lemma 4.12.** *Let $C \subseteq M^n$.*

*(1) If $\mathcal{M}$ is an o-minimal ordered group, and $C$ is a linear cell, then it is an open $\mathbf{Z}$-cell (for our fixed but arbitrary $\mathbf{Z}$).*

*(2) If $\mathcal{M} \models \mathsf{Pres}$, and $C$ is a Presburger cell with $0 \in C$, then $C$ is a purely $\mathbb{Z}$-unbounded closed cell.*

*Proof.* (1) is clear, by Example 4.5, Fact 4.8, and since for every $z \in \pi(C)$, $C_z = \widehat{C}_z$.

(2) By induction on $n$. For $n = 0$, it is clear. For $n > 0$, the graph case is again clear, by Inductive Hypothesis. Let now $C = [\alpha, \beta]^0_r$. Then the last coordinate of an element in $Sh(C)$ is still congruent to $0 \bmod r$. Together with Example 4.5, it follows that $C$ is a closed $\mathbb{Z}$-cell. By Definition 2.14(2), it is purely $\mathbb{Z}$-unbounded. □

We now combine the two freshly defined notions into the following lemma.

**Lemma 4.13.** *Let $D \subseteq M^n$ be a purely $\mathbf{Z}$-unbounded cell, $n > 0$. Then for every $T \in \mathbf{N}$, there are $x, y \in D$, such that $x, y$ are $\mathbf{Z}$-close and $|x - y| > T$.*

*Proof.* By induction on $n$. For $n = 1$, since $D$ is purely unbounded there are $x, y \in D$ with $|x - y| > \mathbf{1}$. By ($\mathbf{Z}$-DENSE), we can find such $x, y$ which are also $\mathbf{Z}$-close. By (UC), $y' = x + T(x - y) \in D$, and hence $x, y'$ satisfy the conclusion.

Now let $n > 1$. Suppose first that $D = \Gamma(\alpha)$. Take any $x, y \in \pi(D)$ that are $\mathbf{Z}$-close with $|x - y| > T$. Then clearly $|(x, \alpha(x)) - (y, \alpha(y))| > T$, and since by Fact 4.7(2), $\alpha(x), \alpha(y)$ are $\mathbf{Z}$-close, so are $(x, \alpha(x))$, $(y, \alpha(y))$.

Suppose now that $\widehat{D} = (\alpha, \beta)$ or $[\alpha, \beta]$. Take again $x, y \in \pi(D)$ that are $\mathbf{Z}$-close with $|x - y| > T$, and let $t, s \in M$ with $0 < t, s < \mathbf{1}$, so that both $(x, \alpha(x) + t)$, $(y, \alpha(y) + s)$ are in $D$. It follows from the graph case that those two elements satisfy the conclusion. □

**Corollary 4.14.** *Let $D \subseteq M^n$ be a purely $\mathbf{Z}$-unbounded cell, $n > 0$, and $\alpha : D \subseteq M^n \to M$ a non-constant bi-Lipschitz function. Then for every $T \in \mathbf{N}$, there are $x, y \in D$, such that $x, y$ are $\mathbf{Z}$-close and $|\alpha(x) - \alpha(y)| > T$.*

*Proof.* Immediate from Lemma 4.13, the definition of a bi-Lipschitz function and the fact that $n\mathbf{Z} \subseteq \mathbf{Z}$ for any $n \in \mathbb{N}$. □

4.3. ***N*-internal points.** We now turn to the notion of $N$-internality. Given an interval $S = (a, b)$ or $[a, b]$, we denote by $\|S\| = b - a$ its length, possibly $\infty$. Observe that for every $\mathbf{Z}$-cell $C$ and $x \in \pi(C)$, $\|C_x\| = \|\widehat{C}_x\|$.

**Definition 4.15.** Let $C \subseteq M^n$ and $N \in \mathbf{N}$. An element $x \in C$ is called *$N$-internal (in $C$)* if

$$B_N(x) \cap Sh(C) \subseteq C.$$

The next definition will play an important role in the proofs that follow.

**Definition 4.16.** Let $C \subseteq M^n$ be a $\mathbf{Z}$-cell, $L \in \mathbf{N}$ and $x \in \pi(C)$. Denote:

$$S_L(x, C) = \bigcap_{y \in B_L(x) \cap \pi(C)} \widehat{C}_y.$$

The *$L$-cylinder above $x$ in $C$* is the set $B_L(x) \times S_L(x, C)$, and its *height* is the quantity $\|S_L(x, C)\|$.

4.4. **Existence of $N$-internal points.** The main result of this subsection is Proposition 4.22, and we begin with a sketch of its proof (Case II: $\widehat{C} = (\alpha, \beta)_D$ or $[\alpha, \beta]_D$). There are two ingredients. First, we prove that the height of an $L$-cylinder above a point $x$ in $C$ is bounded from below by $\|C_x\|$, uniformly in $x$ (Lemma 4.18). The proof is a straightforward application of linearity. Second, we ensure that arbitrarily large (in $\mathbf{N}$) fibers $C_x$ in $C$ exist, within some fixed radius of $x$ (in $\mathbf{N}$), again uniformly in $x$ (Lemma 4.20). This step is a bit more involved as it can only be established when $x$ is sufficiently internal. The two steps combined guarantee that $C$ has enough space 'vertically', namely there is $y \in D$ with $\|S_N(y, C)\| > 2N$, within a fixed radius of any $x$ (Corollary 4.21). Therefore, in the inductive step of the proof of Proposition 4.22(Case II) we can find $N$-internal points $y \in D$ with $\|S_N(y, C)\| > 2N$, and by Claim 4.17(3), $C$ contains an $N$-internal point.

Before establishing our two ingredients, we collect some basic facts. In the following claim, the only assumption on $C$ from Definition 4.10 that is used is (UC).

**Claim 4.17.** *Let $C \subseteq M^n$ be a $\mathbf{Z}$-cell, $n > 0$, and $L, N \in \mathbf{N}$. We have:*

1. *If $x$ is $L$-internal in $C$ and $N \leq L$, then $x$ is also $N$-internal in $C$.*
2. *If $x$ is $2L$-internal in $C$ and $y \in B_L(x) \cap C$, then $y$ is $L$-internal in $C$.*
3. *If $y$ is $N$-internal in $\pi(C)$ and $\|S_N(y, C)\| > 2N$, then there is an $N$-internal point in $C$.*

*Proof.* (1) and (2) are straightforward from the definitions. For (3), pick an element $t$ in $S_L(y, C)$ such that $B_L(t) \subseteq S_L(y, C)$. We claim that the element $z = (y, t)$ is $L$-internal. Take $x = (x_1, x_2) \in B_L(z) \cap Sh(C)$. Then

$$x_1 \in B_L(y) \cap \pi(Sh(C)) = B_L(y) \cap Sh(\pi(C)) \subseteq \pi(C),$$

since $\pi(Sh(C)) = Sh(\pi(C))$ and $y$ is $L$-internal in $\pi(C)$. On the other hand,

$$x_2 \in B_L(t) \cap Sh(C)_{x_1} \subseteq S_L(y, C) \cap Sh(C)_{x_1} \subseteq \widehat{C}_{x_1} \cap Sh(C)_{x_1} = C_{x_1},$$

by choice of $t$, definition of $S_L(y, C)$ and (UC). Hence $x \in C$, as needed. □

We now proceed to our first ingredient.

**Lemma 4.18.** *Let $C \subseteq M^n$ be a $\mathbf{Z}$-cell with $\widehat{C} = (\alpha, \beta)_D$ or $[\alpha, \beta]_D$. Then for all $N \in \mathbf{N}$, there is $Q \in \mathbf{Z}$, such that for all $x \in \pi(C)$,*

$$\|S_N(x, C)\| > \|C_x\| + Q.$$

*Proof.* We only handle the case $\widehat{C} = (\alpha, \beta)_D$, as the other one is almost identical. Let $x \in \pi(C)$ and denote $U_x = B_N(x) \cap \pi(C)$. Since for each $y \in \pi(C)$, $\widehat{C}_y = (\alpha(y), \beta(y))$, it follows easily that $S_N(x, C)$ is either empty, if $\sup_{U_x} \alpha \geq \inf_{U_x} \beta$, or

$$S_N(x, C) = \left( \sup_{U_x} \alpha, \inf_{U_x} \beta \right),$$

otherwise. In either case, $\|S_N(x, C)\| \geq \inf_{U_x} \beta - \sup_{U_x} \alpha$. We can write

$$\begin{aligned} \inf_{U_x} \beta - \sup_{U_x} \alpha &= \inf_{U_x} \beta - \beta(x) + \alpha(x) - \sup_{U_x} \alpha + (\beta(x) - \alpha(x)) \\ &\geq \inf_{U_x} \beta - \beta(x) + \alpha(x) - \sup_{U_x} \alpha + \mathbf{1} + \|C_x\|, \end{aligned}$$

since $\|C_x\| \leq \beta(x) - \alpha(x) + \mathbf{1}$. Therefore, it suffices to prove that there is $Q \in \mathbf{Z}$, such that each of $\inf_{U_x} \beta - \beta(x)$ and $\alpha(x) - \sup_{U_x} \alpha$ is greater than $N$, as $x$ ranges over $\pi(C)$. We prove this for the former as for the latter it is similar.

Since $U_x = B_N(x) \cap \pi(C)$, we have

$$\inf_{U_x} \beta - \beta(x) = \inf\{\beta(x + t) - \beta(x) \,:\, |t| \leq N,\, x + t \in \pi(C)\}.$$

By linearity of $\beta$, for every $x, x' \in \pi(C)$ with $x + t, x' + t \in \pi(C)$, we have

$$\beta(x + t) - \beta(x) = \beta(x' + t) - \beta(x').$$

Hence a lower bound $Q \in \mathbf{Z}$ for $\inf_{U_x} \beta - \beta(x)$, for a fixed $x \in \pi(C)$, is also a lower bound for the same expression for any $x \in \pi(C)$. Such a lower bound exists, by Fact 4.7(2). $\square$

We next proceed towards our second ingredient (Lemma 4.20).

**Claim 4.19.** *Let $D \subseteq M^n$ be a purely $\mathbf{Z}$-unbounded cell, and $\alpha : D \subseteq M^n \to M$ a non-constant definable bi-Lipschitz linear function, with $\alpha(x) \geq 0$, for every $x \in D$. Then for all $T \in \mathbf{N}$, there is $K \in \mathbf{N}$, such that for all $K$-internal points $x \in D$,*

$$\textit{there is } y \in B_K(x) \cap D \textit{ with } \alpha(y) > T.$$

*Proof.* Since $\alpha$ is non-constant, $n > 0$. Suppose $\inf Im(\alpha) = y_0$. Choose any $x_0 \in D$ so that $|\alpha(x_0) - y_0| < \mathbf{1}$. By translating $\Gamma(\alpha)$ so that $(x_0, \alpha(x_0))$ is mapped to the origin, we may assume that $\alpha(0) = 0$ and $\alpha \geq -\mathbf{1}$ instead. Indeed, a $K \in \mathbf{N}$ that satisfies the conclusion of the statement for this translated $\alpha$ would also satisfy it for the original $\alpha$.

Let $T \in \mathbf{N}$. By Corollary 4.14, there are $K \in \mathbf{N}$ and $x_0, z_0 \in D$ with $|x_0 - z_0| < K$ and $\alpha(z_0) - \alpha(x_0) > T + \mathbf{1}$. Let $t = |x_0 - z_0|$. Then for every $K$-internal $x \in D$, we have $x + t \in B_K(x) \cap Sh(D) \subseteq D$, and by linearity of $\alpha$, we obtain

$$\alpha(x + t) - \alpha(x) = \alpha(x_0 + t) - \alpha(x_0) > T + \mathbf{1}.$$

Since $\alpha(x) \geq -\mathbf{1}$, we obtain $\alpha(x + t) > T$. We have shown that for every $K$-internal $x \in D$, there is $y = x + t \in B_K(x) \cap D$ with $\alpha(y) > T$, as needed. $\square$

**Lemma 4.20.** *Let $C \subseteq M^n$ be a purely $\mathbf{Z}$-unbounded cell, with $\widehat{C} = (\alpha, \beta)_D$ or $[\alpha, \beta]_D$. Then for all $T \in \mathbf{N}$, there is $K \in \mathbf{N}$, such that for all $K$-internal $x \in D$,*

$$\textit{there is } y \in B_K(x) \cap D \textit{ with } \|C_y\| > T.$$

*Proof.* If $\beta - \alpha$ is constant, then by definition of purely $\mathbf{Z}$-unbounded cells, for every $x \in D$ and $T \in \mathbf{N}$, $T < \beta(x) - \alpha(x)$, and hence the conclusion of the lemma holds with $y = x$. If $\beta - \alpha$ is not constant, apply Claim 4.19 to the map $\beta - \alpha$. $\square$

Combining our two ingredients, we obtain the following corollary.

**Corollary 4.21.** *Let $C \subseteq M^n$ be as in Lemma 4.20. Then for all $N \in \mathbf{N}$, there is $K \in \mathbf{N}$, so that for every $K$-internal $x \in D$,*

$$\textit{there is } y \in B_K(x) \cap D \textit{ with } \|S_N(x, C)\| > 2N.$$

*Proof.* By Lemma 4.18, there is $Q \in \mathbf{Z}$ such that for all $x \in \pi(C)$,

$$\|S_N(x, C)\| > \|C_x\| + Q.$$

By Lemma 4.20, there is $K \in \mathbf{N}$ such that for every $K$-internal $x \in D$,

$$\text{there is } y \in B_K(x) \cap D \text{ with } \|C_y\| > 2N - Q.$$

For such $y$, we have

$$\|S_N(y, C)\| > \|C_y\| + Q > 2N,$$

as needed. □

We are now ready to prove the key proposition of this subsection, which suggests some notion of relative 'openness' for purely $\mathbf{Z}$-unbounded cells $C$ containing 0.

**Proposition 4.22** (Existence of $N$-internal points)**.** *Let $C \subseteq M^n$ be a purely $\mathbf{Z}$-unbounded cell with $0 \in C$, and $N \in \mathbf{N}$. Then there is an $N$-internal point in $C$.*

*Proof.* We perform induction on $n$. For $n = 0$, $C = Sh(C)$ and the statement is trivial. For the inductive step, let $n > 0$. We split into three cases.

**Case I.** $C = \Gamma(\alpha)_D$. By Inductive Hypothesis, there is an $N$-internal point $x_1$ of $D$. We claim that $x = (x_1, \alpha(x_1))$ is an $N$-internal point of $C$. By Lemma 2.10

$$B_N(x) \cap Sh(C) = B_N(x) \cap \Gamma(\alpha)_{Sh(D)},$$

so take $y = (y_1, \alpha(y_1)) \in B_N(x) \cap \Gamma(\alpha)_{Sh(D)}$. Then $y_1 \in B_N(x_1) \cap Sh(D) \subseteq D$, and hence $y \in \Gamma(\alpha)$, as needed.

**Case II.** $\widehat{C} = (\alpha, \beta)_D$ or $[\alpha, \beta]_D$ with $\alpha, \beta : D \to M$. By Corollary 4.21, there is $K \in \mathbf{N}$, so that for every $K$-internal $x \in D$,

$$\text{there is } y \in B_K(x) \cap D \text{ with } \|S_N(x, C)\| > 2N.$$

Now, let $L = \max\{N, K\}$. By Inductive Hypothesis, there is a $2L$-internal point $x_0$ in $D$. In particular, $x_0$ is $K$-internal in $D$ (Claim 4.17(1)).

By our choice of $K$,

$$\text{there is } y \in B_K(x_0) \cap D \text{ with } \|S_N(y, C)\| > 2N.$$

Since $x_0$ is $2L$-internal in $D$ and $y \in B_L(x_0) \cap D$ (since $K \leq L$), we obtain that $y$ is a $L$-internal point of $D$ (Claim 4.17(2)). By Claim 4.17(1) again, $y$ is $N$-internal in $D$. Since also $\|S_N(y, C)\| > 2N$, we obtain by Claim 4.17(3), that there is an $N$-internal point in $C$.

**Case III.** $\widehat{C} = (\alpha, \beta)$ or $[\alpha, \beta]_D$, with either $\alpha = -\infty$ or $\beta = \infty$. Let $D = \pi(C)$ and take an $N$-internal point $y_1 \in D$ (by Inductive Hypothesis). Observe that in this case $||S_N(y_1, C)|| = \infty$, and hence we are again done by Claim 4.17(3). □

We record the following corollary, although not explicitly used in the sequel.

**Corollary 4.23.** *Let $C \subseteq M^n$ with $0 \in C$.*

(1) *If $\mathcal{M}$ is an o-minimal ordered group, and $C$ a purely $\mathbf{Z}$-unbounded linear cell (for our fixed but arbitrary $\mathbf{Z}$), then for every $N \in \mathbf{N}$, there is an $N$-internal point in $C$.*
(2) *If $\mathcal{M} \models \mathsf{Pres}$ and $C$ is a Presburger cell, then for every $N \in \mathbb{N}$, there is an $N$-internal point in $C$.*

*Proof.* By Lemma 4.12, a linear cell is an open $\mathbf{Z}$-cell, and a Presburger cell is a purely $\mathbb{Z}$-unbounded closed cell. By Proposition 4.22, we are done. □

*Remark* 4.24. Although not used in this paper, it can also be proved that in an ordered vector space $\mathcal{M}$ over an arbitrary ordered division ring, purely unbounded cells still contain $N$-internal points, $N \in M_{\geq 0}$, despite the fact that linear maps are not bi-Lipschitz (Remark 4.9) (and hence the conclusion of Lemma 4.12(1) fails).

4.5. **Corollaries to Proposition 4.22.** We start with a general lemma.

**Lemma 4.25.** *Let $C \subseteq M^n$. For every $x, y \in Sh(C)$ and $N \in \mathbf{N}$,*

$$B_N(x) \cap Sh(C) \sim B_N(y) \cap Sh(C),$$

*witnessed by the map $f : t \mapsto t + (y - x)$.*

*Proof.* Since $Sh(C)$ is a subgroup of $\langle M^n, +\rangle$, we have $f(Sh(C)) = Sh(C)$. Moreover, $f$ yields a bijection between $B_N(x) \cap Sh(C)$ and $B_N(y) \cap Sh(C)$. Since it is a translation, it clearly witnesses the $\sim$-equivalence of the last two sets. □

**Definition 4.26.** By a *copy of* $\mathbf{Z}^n$, we mean a set of the form $a + \mathbf{Z}^n$, for some $a \in M^n$. We say $X \subseteq M^n$ is of $\mathbf{Z}$-*volume* if it is contained in a single copy of $\mathbf{Z}^n$.

**Corollary 4.27.** *Let $C \subseteq M^n$ be a purely $\mathbf{Z}$-unbounded cell with $0 \in C$ and $A \subseteq Sh(C)$ a finite set of $\mathbf{Z}$-volume. Then there is $B \subseteq C$ with $B \sim A$. In fact, $B = a + A$, for some $a \in M^n$.*

*In particular, this statement applies to purely unbounded linear cells in $\langle \mathbb{R}, <, +\rangle$ and to Presburger cells in $\langle \mathbb{Z}, <, +\rangle$, with $\mathbf{Z} = \mathbb{Z}$.*

*Proof.* Since $B$ is of $\mathbf{Z}$-volume, there are $x \in C$ and $N \in \mathbf{N}$ such that

$$A \subseteq B_N(x).$$

(For example, let $x$ be any element in $A$ and $N$ the maximum among all $x_i - y_i$, where $y \in A$, $i \in [n]$.) Let $y \in C$ be an $N$-internal point of $C$, provided by Proposition 4.22. By Lemma 4.25,

$$B_N(x) \cap Sh(C) \sim B_N(y) \cap Sh(C),$$

witnessed by the map $f : t \mapsto t + (y - x)$. Let $B = f(A) \sim A$. Since $y$ is $N$-internal, we have

$$B \subseteq B_N(y) \cap Sh(C) \subseteq C,$$

as needed. The 'in particular' clause is by Lemma 4.12. □

*Remark* 4.28. Corollary 4.27 fails if we allowed $A$ to be infinite. Indeed, let $\mathcal{M} = \langle \mathbb{Z}, <, +\rangle$, $\mathbf{Z} = \mathbb{Z}$, and

$$C = \{(x, y) \in \mathbb{Z}^2 : 0 < y < x\}.$$

Then $Sh(C) = \mathbb{Z}^2$, and the statement fails for $A = \mathbb{Z}^2$ (even without requiring $B'$ to be a translate of $A$). This example also shows the failure of $Zar^{\infty}(C)$ in $\mathcal{M}$. (For saturated models of $\mathsf{Pres}$, however, it holds, Remark 6.19.)

Although not needed in what follows, we can now conclude a '$\mathbf{Z}$-version' of Zarankiewicz's problem for our arbitrary definably complete ordered group.

**Proposition 4.29.** *Assume that $\mathbf{Z}$ is cofinal in $M$. Let $C \subseteq M^n$ be a purely $\mathbf{Z}$-unbounded cell that contains $0$. The following are equivalent:*

(1) *For every $l \in [r]$ and $k \in \mathbb{N}$, $p_{l\restriction Sh(C)}$ is not $k$-to-$1$.*
(2) *For every $k \in \mathbb{N}$, $C$ contains a $k$-grid.*

*Hence, $C$ satisfies $(*)$: if $C$ is $k$-free, for some $k \in \mathbb{N}$, then so is $Sh(C)$. In particular, $Zar(C)$ holds.*

*Proof.* (2)⇒(1) is clear. For (1)⇒(2), by Lemma 2.6 and Fact 2.5, there is an infinite grid $G = G_1 \times \cdots \times G_r \subseteq Sh(C)$, with all $G_i$ infinite subgroups of $M^{d_i}$. Since $\mathbf{Z}$ is cofinal in $M$, $G$ is of $\mathbf{Z}$-volume. It follows that we can find a $k$-grid inside $Sh(C)$ of $\mathbf{Z}$-volume. By Corollary 4.27, $C$ contains a $k$-grid, as needed. We conclude by Proposition 2.17. ☐

The above proof could not yield that $C$ satisfies $(*)_\infty$ because Corollary 4.27 does not apply to infinite grids (Remark 4.28).

**Corollary 4.30.** *Assume that $\mathbf{Z}$ is cofinal in $M$, and every definable set is a finite union of purely $\mathbf{Z}$-unbounded cells. Let $E \subseteq \prod_{i\in[r]} M^{d_i}$ be definable. Then $Zar(E)$.*

*Proof.* By Proposition 4.29 and Corollary 2.18. ☐

We can also obtain Theorem B in the standard model.

**Corollary 4.31.** *Let $\mathcal{M} = \langle \mathbb{Z}, <, + \rangle$ and $\mathbf{Z} = \mathbb{Z}$. Let $E \subseteq \prod_{i\in[r]} M^{d_i}$ be definable. Then $Zar(E)$.*

*Proof.* The assumptions of Corollary 4.30 hold. ☐

4.6. **Proof of Theorem B (Corollary 4.35).** As with the standard model (Proposition 4.29 and Corollary 4.31), our strategy consists in proving Property $(*)$ from Section 2.7 for a Presburger cell $C$ containing 0. The key point in doing so is to obtain Corollary 4.27 without assuming that $A$ is of $\mathbf{Z}$-volume and Proposition 4.29 without $\mathbf{Z}$ being cofinal in $M$. We do this for the case of a model of Presburger arithmetic and $\mathbf{Z} = \mathbb{Z}$ (Corollary 4.33 and Proposition 4.34, respectively).

*In this subsection, we fix a model $\mathcal{M} = \langle M, <, + \rangle \models \mathsf{Pres}$.* Recall $\langle \mathbb{Z}, <, + \rangle \preccurlyeq \mathcal{M}$.

The following proposition in a sense suggests that the $\mathbb{Z}$-density of a Presburger cell (Lemma 4.12) 'passes' to its shell.

**Proposition 4.32.** *Let $C \subseteq M^n$ be a Presburger cell with $0 \in C$, and $A \subseteq Sh(C)$ a finite set. Then there is $B \subseteq Sh(C) \cap \mathbb{Z}^n$ with $A \sim B$.*

Before proving Proposition 4.32, we illustrate how it implies the Zarankiewicz statements for any Presburger set. First, it implies the desired strengthening of Corollary 4.27 in the Presburger setting (although $B$ need not be a translate of $A$).

**Corollary 4.33.** *Let $C \subseteq M^n$ be a Presburger cell with $0 \in C$, and $A \subseteq Sh(C)$ a finite set. Then there is $B \subseteq C$ with $A \sim B$.*

*Proof.* By Proposition 4.32, there is $B' \subseteq Sh(C)$ of $\mathbb{Z}$-volume with $B' \sim A$. By Corollary 4.27, there is $B \subseteq C$ with $B \sim B' \sim A$, as needed. ☐

The following proposition is an analogue of Lemma 2.13 from the linear case.

**Proposition 4.34.** *Let $C \subseteq M^n$ be a Presburger cell that contains $0$. The following are equivalent:*

*(1) For every $l \in [r]$, $p_{l \restriction Sh(C)}$ is non-injective.*
*(2) For every $k \in \mathbb{N}$, $C$ contains a $k$-grid.*

*Hence, $C$ satisfies $(*)$: if $C$ is $k$-free, for some $k \in \mathbb{N}$, then so is $Sh(C)$. In particular, $Zar_1(C)$ holds.*

*Proof.* (2)⇒(1) is clear. For (1)⇒(2), by Lemma 2.6, Fact 2.5 and the beginning of Section 2.3, there is an infinite grid $G = G_1 \times \cdots \times G_r \subseteq Sh(C)$, with all $G_i$ infinite subgroups of $M^{d_i}$. It follows that we can find a $k$-grid inside $Sh(C)$. By Corollary 4.33, $C$ contains a $k$-grid, as needed. We conclude by Proposition 2.17. □

**Corollary 4.35.** *Let $\mathcal{M} \models \mathsf{Pres}$. Let $E \subseteq \prod_{i \in [r]} M^{d_i}$ be definable. Then $Zar(E)$.*

*Proof.* By Presburger cell decomposition, Proposition 4.34 and Corollary 2.18. □

We now turn to the proof of Proposition 4.32. Our strategy is this: given a Presburger cell $C \subseteq M^n$, we introduce in Definition 4.36 the '$\mathbb{Z}$-cover' $\mathcal{Z}(C)$ of $C$, another Presburger cell that is $\mathbb{Z}$-definable, contains $C$, and satisfies $Sh(C) \cap \mathbb{Z}^n = Sh(\mathcal{Z}(C)) \cap \mathbb{Z}^n$. Since the finite set $A \subseteq Sh(C) \subseteq \mathcal{Z}(C)$ is only contained in a $\mathbb{Z}$-definable 'fragment' of $Sh(\mathcal{Z}(C))$, Proposition 4.32 reduces to the same statement but in the standard model $\langle \mathbb{Z}, <, + \rangle$, where it trivially holds.

We begin with some useful notation. Given a $\mathsf{Pres}$-linear map $f : D \subseteq M^n \to M$, we let $\zeta(f) : D \to M$ with

$$\zeta(f) = \begin{cases} \infty & \text{if } f(0) > \mathbb{Z} \\ f & \text{if } f(0) \in \mathbb{Z} \\ -\infty & \text{if } f(0) < \mathbb{Z}, \end{cases}$$

Note that if $f(0) \in \mathbb{Z}$, then $f$ is $\mathbb{Z}$-definable (by the form of $\mathsf{Pres}$-linear maps, Section 2.6). In the next definition, if $f, g : X \subseteq M^n \to M$ are two $\mathsf{Pres}$-linear maps, or $f = -\infty$ or $g = \infty$, by $[f, g]_r^c$ we denote the $\mathsf{Pres}$-cylinder $[f', g']_r^c$ (Section 2.6), where $f', g'$ are the restrictions of $f, g$ to the set $X' = \{x \in X : f(x) \leq g(x)\}$. Note also that if $f, g$ are $\mathbb{Z}$-definable, then so is $[f, g]_r^c$.

**Definition 4.36.** Let $C \subseteq M^n$ be a Presburger cell with $0 \in C$. We define the *$\mathbb{Z}$-cover* of $C$, denoted by $\mathcal{Z}(C)$, recursively, as follows. For $C \subseteq M^0$, let $\mathcal{Z}(C) = C$. Let $n > 0$.

- If $C = \Gamma(\alpha)_D$, then $\mathcal{Z}(C) = \Gamma(\alpha)_{\mathcal{Z}(D)}$.
- If $C = [\alpha, \beta]_r^c$ with $\alpha, \beta : D \to M$, or $\alpha = -\infty$ or $\beta = \infty$, then $\mathcal{Z}(C) = [\zeta(\alpha), \zeta(\beta)]_r^c$ with $\zeta(\alpha), \zeta(\beta) : \mathcal{Z}(D) \to M$.

**Note:** The above definition runs along with the following statement, which can be easily proved by induction (and is left to the reader): if $\gamma : C \to M$ is a $\mathsf{Pres}$-linear map, then there is unique $\mathsf{Pres}$-linear extension of $\gamma$ to $\mathcal{Z}(C)$, also denoted by $\gamma$.

We next prove some key properties of $\mathbb{Z}$-covers.

**Lemma 4.37.** *Let $C \subseteq M^n$ be a Presburger cell with $0 \in C$. Then:*

*(1) $\mathcal{Z}(C)$ is $\mathbb{Z}$-definable.*
*(2) $C \subseteq \mathcal{Z}(C)$.*

*(3)* $Sh(\mathcal{Z}(C)) \cap \mathbb{Z}^n = Sh(C) \cap \mathbb{Z}^n$.

*Proof.* Both (1) and (2) are straightforward from the definitions, so we only prove (3). We work by induction on $n$. For $n = 0$, the statement is obvious. Let $n > 0$, $C \subseteq M^n$ and write $E$ for $\mathcal{Z}(C)$. We split into two cases.

**Case I.** $C = \Gamma(\alpha)_D$. Since $\alpha(0) = 0$, by Lemma 2.10, we have

$$Sh(E) \cap \mathbb{Z}^n = \Gamma(\alpha)_{Sh(\mathcal{Z}(D))\cap\mathbb{Z}^{n-1}} = \Gamma(\alpha)_{Sh(D)\cap\mathbb{Z}^{n-1}} = Sh(C) \cap \mathbb{Z}^n,$$

where the middle equality is by Inductive Hypothesis.

**Case II.** $C = [\alpha, \beta]^0_r$. We may assume that not both $\alpha = -\infty$ and $\beta = \infty$, since in that case, $\zeta(\alpha) = -\infty$ and $\zeta(\beta) = \infty$ and the statement follows from Inductive Hypothesis. We may assume $\alpha \neq -\infty$, the other case being handled similarly.

Observe that by definition of cells, shells and $\mathbb{Z}$-covers, we have that

$$Sh(C), Sh(E) \subseteq M^{n-1} \times rM,$$

and, in particular

$$Sh(E) \cap \mathbb{Z}^n \subseteq Sh(\pi(E)) \times r\mathbb{Z}$$

and

$$Sh(C) \cap \mathbb{Z}^n \subseteq Sh(\pi(C)) \times r\mathbb{Z}.$$

Let $x \in \pi(C)$ be so that $|C_x| > 1$. Then there is $t \in M$, such that $(x, t), (x, t + r) \in C$. It follows that $(0, r) \in Sh(C)$.

**Claim.** $(Sh(\pi(C)) \cap \mathbb{Z}^{n-1}) \times \{0\} \subseteq Sh(C)$.

*Proof of Claim.* We may assume $\alpha(0) = 0$. Indeed, for $C' = C - (0, \alpha(0))$, we obtain $Sh(C') = Sh(C)$, and hence we may prove the statement for $C'$ instead of $C$.

Now let $x = x_1 \pm \cdots \pm x_l \in Sh(\pi(C)) \cap \mathbb{Z}^{n-1}$, with $x_1, \ldots, x_l \in \pi(C)$. Then $y = (x_1, \alpha(x_1)) \pm \cdots \pm (x_l, \alpha(x_l)) \in Sh(C)$, and since $\alpha(0) = 0$, we also have $y \in \mathbb{Z}^n$. Therefore, $\alpha(x_1) \pm \cdots \pm \alpha(x_l) \in r\mathbb{Z}$. Since $(0, r) \in Sh(C)$, it follows that

$$(x, 0) \in y + (\{0\} \times r\mathbb{Z}) \subseteq Sh(C),$$

as needed. □

By the claim and again the fact that $(0, r) \in Sh(C)$ it follows that

$$(Sh(\pi(C)) \cap \mathbb{Z}^{n-1}) \times r\mathbb{Z} \subseteq Sh(C) \cap \mathbb{Z}^n.$$

Hence

$$Sh(E) \cap \mathbb{Z}^n \subseteq (Sh(\pi(E)) \cap \mathbb{Z}^{n-1}) \times r\mathbb{Z} = (Sh(\pi(C)) \cap \mathbb{Z}^{n-1}) \times r\mathbb{Z} \subseteq Sh(C) \cap \mathbb{Z}^n,$$

as needed, where the equality is by Inductive Hypothesis. □

We are now ready to conclude this subsection.

*Proof of Proposition 4.32.* Let $m = |A|$ and $\phi(x_1, \ldots, x_m)$ the formula (over $\emptyset$) expressing that $\{x_1, \ldots, x_m\} \sim A$.

Suppose that $A = \{a_1, \ldots, a_m\}$. By Lemma 4.37(2), $A \subseteq Sh(C) \subseteq Sh(\mathcal{Z}(C))$. Since $A$ is finite and $Sh(\mathcal{Z}(C))$ is generated by the $\mathbb{Z}$-definable $\mathcal{Z}(C)$ (Lemma 4.37(1)), there exist formulas $\sigma_1(x), \ldots, \sigma_m(x)$ over $\mathbb{Z}$, such that for each $i \in [m]$ we have $\sigma_i(\mathcal{M}) \subseteq Sh(\mathcal{Z}(C))$ and $\mathcal{M} \vDash \sigma_i(a_i)$. Thus

$$\mathcal{M} \vDash \exists x_1, \ldots, x_m \phi(x_1, \ldots, x_m) \wedge \bigwedge \sigma_i(x_i),$$

and the sentence above has parameters in $\mathbb{Z}$. Therefore

$$\langle \mathbb{Z}, <, + \rangle \models \exists x_1, \dots, x_m \phi(x_1, \dots, x_m) \wedge \bigwedge \sigma_i(x_i).$$

Take $B = \{b_1, \dots, b_m\} \subseteq \mathbb{Z}^n$ witnessing $\phi(x_1, \dots, x_m) \wedge \bigwedge \sigma_i(x_i)$. Then $B \sim A$ and $B \subseteq Sh(\mathcal{Z}(C)) \cap \mathbb{Z}^n$. By Lemma 4.37(3), $B \subseteq Sh(C) \cap \mathbb{Z}^n$, as needed. □

4.7. **Product cells, nesting lines, and purely unbounded sets revisited.** We finally introduce three notions that will be used in Section 5, but their definition can be given in the current setting. We begin with that of a 'product cell', which is in a sense a product of a purely unbounded cell and a bounded linear cell.

**Definition 4.38.** A cell $C \subseteq M^n$ is a *product cell* if $C = J + D$, where

(1) $J$ is a purely unbounded cell and $D$ a bounded cell, and
(2) for every $c \in C$, there are unique $g \in J$ and $d \in D$, with $c = g + d$.

If we write '$C = J + D$ is a product cell', we mean (1)&(2) for these $J$ and $D$ hold.

The following two notions will later on (Sections 5 and Appendices A.1 and A.2) be specialised to semilinear and Presburger sets.

**Definition 4.39.** A *line* $l \subseteq M^n$ is the image $Im(\gamma)$ of a linear map $\gamma : M_{\geq 0} \to M^n$. Such $l$ inherits the order of $M$ via $\gamma$, also denoted by $<$.

**Definition 4.40.** Let $\mathcal{X} = \{X_t\}_{t \in Y}$ be a family of sets in $M^n$, with $Y \subseteq M^m$. A line $l \subseteq Y$ is a *nesting line for* $\mathcal{X}$ if the following hold:z

(1) for all $t, t' \in l$, $t < t' \Rightarrow X_t \subseteq X_{t'}$, and
(2) for every $t \in Y$, there is $t' \in l$, such that $X_t \subseteq X_{t'}$.

A line $l$ containing 0 is called a *nesting direction for* $\mathcal{X}$ if for every $p \in Y$, $p + l$ is a nesting line for $\mathcal{X}$.

A subtle use of purely unbounded cells will be made in Section 5.4, involving the following notion. Recall that if $C \subseteq M^n$ and $x = (x_1, \dots, x_n) \in C$, then $(C_{\pi(x)})_{>x_n}$ (or $(C_{\pi(x)})_{<x_n}$) denotes the set of elements in the fiber $C_{\pi(x)}$ of $C$ that are larger (or smaller) than $x_n$. We identify those lines in $C$ along which both of those sets become arbitrarily large.

**Definition 4.41.** Let $C \subseteq M^n$, $n > 0$, and $l \subseteq C$ a line. We say that $l$ *witnesses that $S$ is purely unbounded* if $\pi(l)$ witnesses that $\pi(S)$ is purely unbounded and one of the following holds:

(1) $C = \Gamma(\alpha)_D$,
(2) for every $N \in \mathbb{N}$, there is $x \in l$, such that

$$(C_{\pi(x)})_{<x_n} \not\subseteq B_N(x_n) \quad \text{and} \quad (C_{\pi(x)})_{>x_n} \not\subseteq B_N(x_n).$$

The above notions will be connected in Proposition 5.17 and Corollary 5.18, which assert that for certain semilinear families, a line witnessing that the parameter set is purely unbounded is also a nesting line for the family. For now we only observe that if a line $l$ witnesses that $C$ is purely unbounded, then $C$ is purely unbounded.

## 5. Zarankiewicz problem for $\langle \mathbb{R}, <, +, \mathbb{Z} \rangle$

*In this section, 'semilinear set/linear map/linear cell' and 'Presburger set/*Pres*linear map/Presburger cell' are taken in $\langle \mathbb{R}, <, + \rangle$ and $\langle \mathbb{Z}, <, + \rangle$, respectively.*

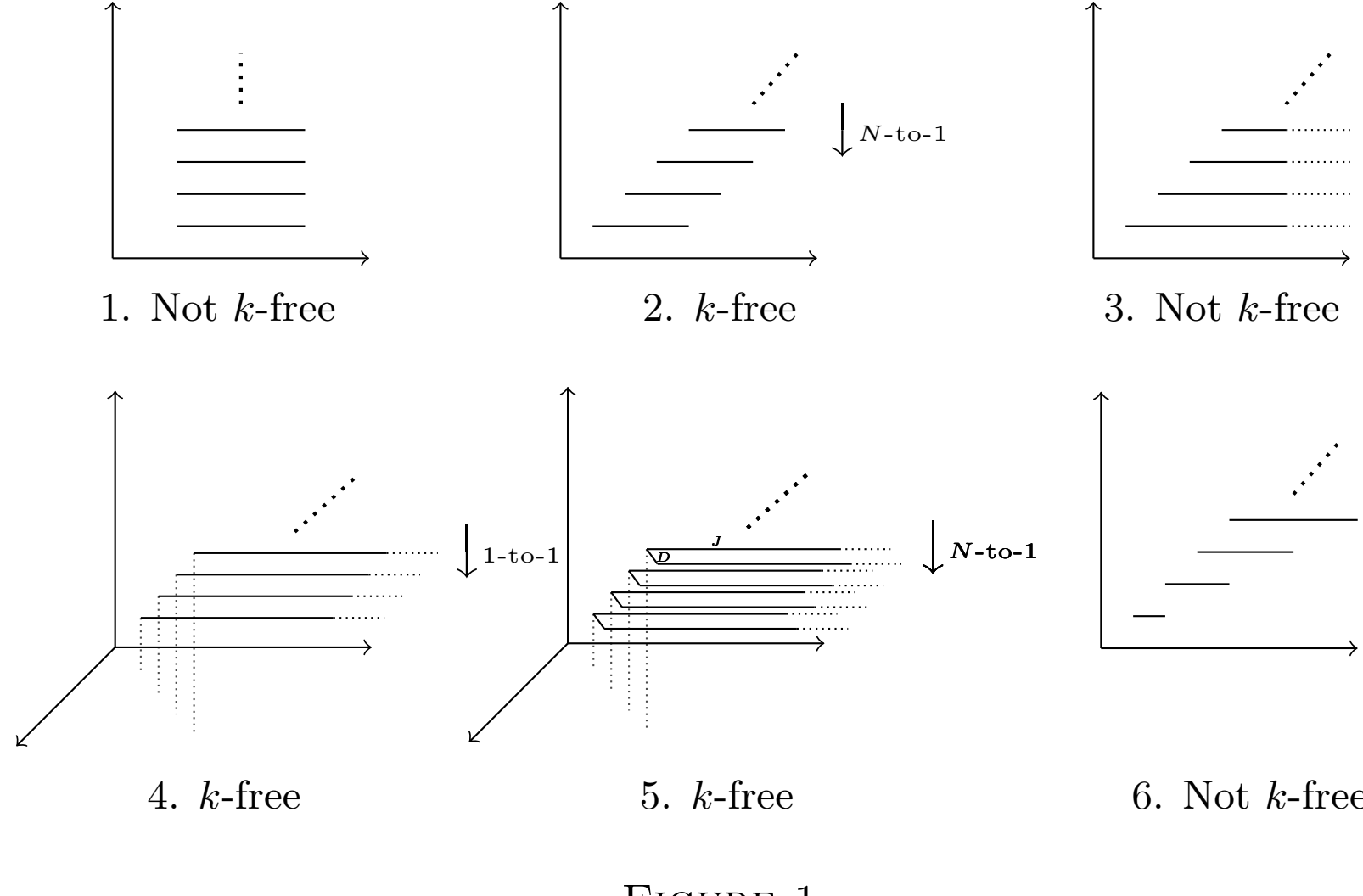

Figure 1

We prove a Zarankiewicz statement in the third and most intriguing geometrical setting of our paper, namely $\langle \mathbb{R}, <, +, \mathbb{Z} \rangle$ (Theorem 5.21 below). The main structure theorem for definable sets in this setting is given by Fact 5.9. A special case of a definable set is that of $E = f(S) + C$ (the 'sum case'), where $S$ is a Presburger set, $C$ a semilinear set, and $f$ a linear map. We first prove $Zar(E)$ for such an $E$ in Section 5.1 (Corollary 5.7), using the reduction to subgroups via $(*)$ mentioned in Section 2.7 (Lemma 5.5). As a preliminary tool, we prove in Proposition A.2 that every linear cell can be partitioned into finitely many product cells (Definition 4.38), which may be a result of its own interest.

In Sections 5.3–5.6, we reduce the general case to Corollary 5.7, after proving Proposition 5.20, which establishes that if a definable set of an appropriate form is $k$-free ($k \in \mathbb{N}$), then it is contained in a $k$-free set of the form $f(S) + C$, as in the sum case. Hence the majority of the content of Sections 5.3–5.5 consists of making sure we may assume that our definable sets are of that appropriate form. Towards proving Proposition 5.20, we establish yet another result in linear o-minimal structures, again of independent interest: for certain semilinear families, a line witnessing that the parameter set is purely unbounded (Definition 4.41) is a nesting line (Definition 4.40, Proposition 5.17 and Corollary 5.18). A sketch of the proof for this general case and hence Theorem 5.21 is given in Section 5.2.

Working in the standard model is essential in this section, due to the fact that we are applying Corollary 4.27 in the proof of Lemma 5.5 below. We also use the fact that linear maps have scalars in $\mathbb{Q}$ (versus any ordered division ring) in Fact 5.2 and the proof of Lemma 5.3.

Before starting with our plan, let us present in Figure 1, pictures of four binary relations in $\mathbb{R}^2$ and two ternary ones in $\mathbb{R}^3$, indicating whether they are $k$-free and which projections are $N$-to-1, for some $N \in \mathbb{N}$. The first five belong to the sum case (Section 5.1), whereas the sixth reduces to the third one that contains it as described above (if the sixth were $k$-free, so would have to be the third).

5.1. **Sum case.** In this subsection, we prove $Zar(E)$ for the 'sum case' $E = f(S)+C$, where $S$ is a Presburger set, $C$ a semilinear set, and $f$ a linear map (Corollary 5.7). The key statement is Lemma 5.4, which finds a uniformly finite-to-1 projection on $E$, assuming further that $C = J + D$ is a product cone and $S, J$ are subgroups of $\langle \mathbb{R}^n, +\rangle$ (Figure 1(5)). We reduce the sum case to this one in the proof of Proposition 5.6, after implementing our Reduction Strategy (3) from Section 2.7 in Lemma 5.5.

We begin with some preliminary lemmas. The first one for $D = \{0\}$ follows from Theorem 2.7, and for $S = \{0\}$ from Lemma 2.13 (for $q = 0$) and Fact 2.20.

**Lemma 5.1.** *Suppose $S, J \leqslant \langle \mathbb{R}^n, +\rangle$ and $D \subseteq \mathbb{R}^n$ is a linear cell, with $J + D$ a linear cell. Suppose further that $S + J + D$ is $k$-free, for some $k \in \mathbb{N}$. Then there is $l \in [r]$ such that both $p_{l\restriction S+J}$ and $p_{l\restriction J+D}$ are injective.*

*Proof.* We may assume that $0 \in J+D$ (by replacing $D$ by $D-a$, for any $a \in J+D$). Suppose the conclusion fails, namely for every $l \in [r]$, we have either $p_{l\restriction S+J}$ or $p_{l\restriction J+D}$ is non-injective. To facilitate notation, suppose the former case happens for $l \in [q]$ and the latter for $l \in [r] \setminus [q]$, for some $q \in \{0, \dots, r\}$. By Fact 2.5, each $U_l^{S+J}$ is infinite, $l \in [q]$. By Lemma 2.13 and Fact 2.20, there are infinite $D_l \subseteq \mathbb{R}^{d_l}$, $l \in [r] \setminus [q]$, such that

$$\{0\} \times D_{q+1} \times \cdots \times D_r \subseteq J + D.$$

Hence

$$\begin{gathered}\left(U_1^{S+J} \times \cdots \times U_q^{S+J} \times \{0\}\right) + \left(\{0\} \times D_{q+1} \times \cdots \times D_r\right) \subseteq \\ S + J + J + D = S + J + D.\end{gathered}$$

In particular, there are infinite $K_i$, $i \in [r]$, with $K_1 \times \cdots \times K_r \subseteq S+J+D$, showing that the latter is not $k$-free, for any $k \in \mathbb{N}$, a contradiction. □

**Fact 5.2.** *Let $q_1, \dots, q_n \in \mathbb{Q}$. Then $q_1\mathbb{Z} + \cdots + q_n\mathbb{Z} = p\mathbb{Z}$, for some $p \in \mathbb{Q}$.*

*Proof.* It is easy to show it for $n = 2$, and then proceed by induction. □

The next lemma might be well-known, but we give a proof, for completeness. Since we are working in the standard model, the scalars of the linear map $f$ below are in $\mathbb{Q}$.

**Lemma 5.3.** *Let $J, f$ be definable in $\langle \mathbb{R}, <, +\rangle$, where $J \leqslant \langle \mathbb{R}^n, +\rangle$ is a subgroup of dimension $< n$, and $f : \mathbb{R}^m \to \mathbb{R}^n$ a linear map with $f(0) = 0$. Let $T \in \mathbb{N}$. Then there is $N \in \mathbb{N}$, such that for every $a \in \mathbb{R}^n$, there are at most $N$ sets of the form $z + J$, $z \in f(\mathbb{Z}^m)$, with $(z + J) \cap B_T(a) \neq \emptyset$.*

*Proof.* Say $0 < \mathsf{dim} J < n$. After permuting coordinates, we may assume that $J$ projects injectively to the first $n - l$ coordinates. Since $\mathsf{dim} J = l$, this implies that for every $z \in f(\mathbb{Z}^m)$, $(z + J) \cap (\{0\} \times \mathbb{R}^{n-l})$ is a singleton.

**Claim.** *It is enough to prove that*

$$(f(\mathbb{Z}^m) + J) \cap (\{0\} \times \mathbb{R}) = \{0\} \times r\mathbb{Z},$$

*for some $r \in \mathbb{R}$.*

*Proof of Claim.* By increasing $T$ if necessary, we only need to prove the conclusion of the lemma for $a \in \mathbb{Z}^n$. Hence, we may assume $a = 0$, and prove that the set $(f(\mathbb{Z}^m) + J) \cap B_T(0)$ is finite. By increasing $T$ further, and the linearity of $f$, it is enough to prove that $(f(\mathbb{Z}^m)+J)\cap(\{0\}\times\mathbb{R}^{n-l})$ is finite. Moreover, it suffices to prove that for each of the last $n-l$ axes, call it $Y_i$, $i=l+1,\dots,n$, we have

$$(f(\mathbb{Z}^m)+J)\cap Y_i \subseteq \{0\}\times r_i\mathbb{Z}\times\{0\}, \tag{1}$$

for some $r_i\in\mathbb{R}$. Indeed, it is not hard to see that this implies

$$(f(\mathbb{Z}^m)+J)\cap(\{0\}\times\mathbb{R}^{n-l})\subseteq\{0\}\times r_{l+1}\mathbb{Z}\times\cdots\times r_n\mathbb{Z},$$

as needed. Without loss of generality, it is enough to prove (1) for $i=n$. □

We now prove the formula of the claim. Let $f=(f_1,\dots,f_n):\mathbb{R}^m\to\mathbb{R}^n$, where each $f_i(z)=f_i^1z_1+\cdots+f_i^mz_m$, with $f_i^j\in\mathbb{Q}$. Let $J$ be generated by $v_1,\dots,v_l$, with $v_i=(v_i^1,\dots,v_i^n)$, for some $v_i^j\in\mathbb{Q}$ (because $J$ is definable in $\langle\mathbb{R},<,+\rangle$). An element of $(f(\mathbb{Z}^m)+J)\cap(\{0\}\times\mathbb{R})$ has form

$$(f_1(z),\dots,f_n(z))+t_1v_1+\cdots+t_lv_l=(0,\dots,0,x),\ \text{ where } t_i,x\in\mathbb{R},$$

which gives rise to the system

$$\begin{pmatrix} f_1(z)\\ \vdots\\ f_{n-1}(z)\end{pmatrix}+\begin{pmatrix} v_1^1 & \dots & v_l^1\\ \vdots & \dots & \vdots\\ v_1^{n-1} & \dots & v_l^{n-1}\end{pmatrix}\begin{pmatrix} t_1\\ \vdots\\ t_l\end{pmatrix}=0,$$

and the equation

$$f_n(z)+t_1v_1^n+\cdots+t_lv_l^n=x.$$

Since all of $v_k^n$ and $f_i^j$ are rationals, and $\mathsf{dim}J=l$, any solution to the system is of the form $t_1,\dots,t_l$, where $t_i=a_1z_1+\cdots+a_mz_m$, for some fixed $a_i\in\mathbb{Q}$. Combining with the equation, we obtain

$$x=b_1z_1+\cdots+b_mz_m,$$

for some fixed $b_i\in\mathbb{Q}$. Therefore, by Fact 5.2, $x\in r\mathbb{Z}$, for some fixed $r\in\mathbb{R}$ (in fact, $r\in\mathbb{Q}$), as desired. □

**Lemma 5.4.** *Suppose $S\leqslant\langle\mathbb{Z}^m,+\rangle$, $J\leqslant\langle\mathbb{R}^n,+\rangle$, and $D\subseteq\mathbb{R}^n$ is a bounded linear cell, with $J+D$ a product cell. Let $f:\mathbb{R}^m\to\mathbb{R}^n$ be a linear map with $f(0)=0$. Suppose that $X=f(S)+J+D$ is $k$-free, for some $k\in\mathbb{N}$. Then there are $l\in[r]$ and $N\in\mathbb{N}$, such that $p_{l\restriction X}$ is $N$-to-1 (illustrated in Figure 1(5)).*

*Proof.* Observe that $f(S)\leqslant\langle\mathbb{R}^n,+\rangle$. By Lemma 5.1, there is $l$ such that both $p_{l\restriction f(S)+J}$ and $p_{l\restriction J+D}$ are injective. We may assume that $f(S)\not\subseteq J$, otherwise we are done. Denote $\pi=p_l$. Let $-D\subseteq B_T(0)$ and $N$ as in Lemma 5.3. Fix $a\in\pi(\mathbb{R}^n)$. We want to prove that the set

$$A=\{x\in f(S)+J+D:\pi(x)=a\}$$

has size at most $N$. If $x=s+g+d$, where $s\in S,g\in J,d\in D$, with $\pi(x)=a$, then

$$\pi(s+g)=a-\pi(d)\in B_T(a).$$

Hence the set $A$ has cardinality at most that of

$$\{w\in f(S)+J:\pi(w)\in B_T(a)\cap(a-\pi(D))\},$$

which in turn, since $\pi_{\restriction f(S)+J}$ is injective, has cardinality at most that of

$$\{y\in\pi(f(S)+J):y\in B_T(a)\cap(a-\pi(D))\},$$

which is a subset of the set

$$C = \pi(f(\mathbb{Z}^m) + \pi(J)) \cap B_T(a) \cap (a - \pi(D)).$$

Observe that $\pi(J) \leqslant \langle \mathbb{R}^{n-1}, +\rangle$ has dimension $< n-1$. Indeed, if $\mathsf{dim}\pi(J) = n-1$, then $\pi(J) = \mathbb{R}^{n-1}$, and since $f(S) \not\subseteq J$, it follows that $\pi_{\restriction f(S)+J}$ is not injective, a contradiction. By Lemma 5.3, there are at most $N$ sets of the form $z + \pi(J)$, $z \in f(\pi(S)) = \pi(f(S))$, that intersect $B_T(a)$. Since $C$ is contained in the union of those sets, it suffices to prove the following claim.

**Claim.** *Let $z \in \pi(f(S))$. Then $|(z + \pi(J)) \cap (a - \pi(D))| \leq 1$.*

*Proof of Claim.* Suppose there are two such elements $y_1 = z + \pi(g_1)$ and $y_2 = z + \pi(g_2)$, $g_1, g_2 \in J$, with

$$y_1 + \pi(d_1) = y_2 + \pi(d_2) = a,$$

for some $d_1, d_2 \in D$. Hence

$$\pi(g_1) + \pi(d_1) = \pi(g_2) + \pi(d_2),$$

Since $\pi_{\restriction J+D}$ is injective,

$$g_1 + d_1 = g_2 + d_2.$$

Since $J + D$ is a product cell, we have $g_1 = g_2$ and $d_1 = d_2$, and hence $x_1 = x_2$. □

This completes the proof of lemma. □

If $E = S+C$ for some Presburger cell $S$ and linear cell $C$, we cannot claim that if $E$ is $k$-free then so is $Sh(E)$ (as we did for $S$ and $C$, Proposition 4.34 and Corollary 2.21, respectively). For example, the shell of the set in Figure 1(2) contains the set in (3). We remedy this as follows.

**Lemma 5.5.** *Suppose $S \subseteq \mathbb{Z}^m$ is a Presburger cell, $J + D$ a product cell, and $f : \mathbb{R}^m \to \mathbb{R}^n$ a linear map, such that $0 \in f(S) \cap J$. Let $E = f(S) + J + D$. If $E$ is $k$-free, then so is $f(Sh(S)) + Sh(J) + D$.*

*Proof.* The idea is to apply Corollary 4.27 to both $\langle \mathbb{Z}, <, +\rangle$ and to $\langle \mathbb{R}, <, +\rangle$, noting that in both structures, grids are of $\mathbb{Z}$-volume. More precisely, it is enough to find, for every finite set $B \subseteq f(Sh(S))+Sh(J)+D$, a set $B' \sim B$ with $B' \subseteq f(S)+J+D$. Given $b \in B$, choose $s(b) \in Sh(S)$, $g(b) \in Sh(J)$, $d(b) \in D$, such that

$$b = f(s(b)) + g(b) + d(b)$$

(they need not be unique). Let $T = \{s(b) : b \in B\}$ and $W = \{g(b) : b \in B\}$. By Corollary 4.27 applied to $\langle \mathbb{Z}, <, +\rangle$, there is $a \in \mathbb{R}^n$, such that $a + T \subseteq S$, and since $f$ is linear, there is also $a' \in \mathbb{R}^n$, such that $a' + f(T) = f(S)$. By Corollary 4.27 applied to $\langle \mathbb{R}, < +\rangle$, there is $c \in \mathbb{R}^n$ such that $c + W \subseteq J$. It follows that the set

$$B' = \{a' + f(s(b)) + c + g(b) + d(b) : b \in B\}$$

satisfies both $B' \subseteq f(S) + J + D$ and $B' = a' + c + B \sim B$, as needed. □

We can now conclude the main result of this subsection.

**Proposition 5.6.** *Let $E = f(S)+C$, where $S \subseteq \mathbb{Z}^m$ is a Presburger cell, $C = J+D$ a product cell, and $f : \mathbb{R}^m \to \mathbb{R}^n$ a linear map with $f(0) = 0$. If $E$ is $k$-free for some $k \in \mathbb{N}$, then some projection $p_{l\restriction E}$, $l \in [r]$, is $N$-to-1. In particular, $Zar(E)$.*

*Proof.* Let $S' = S - a$ and $J' = J - g$, for some $a \in S$ and $g \in J$. Then $E' = f(S') + J' + D$ is a translate of $E$, and it suffices to show the proposition for it. So we may assume $0 \in f(S) \cap J$. By Lemma 5.5, $f(Sh(S)) + Sh(J) + D$ is $k$-free. We note also that $Sh(J)$ is purely unbounded (since $J$ is), and definable in $\langle \mathbb{R}, <, + \rangle$ (as can easily be seen). By Lemmas 5.4 and Lemma 2.1, we are done. □

**Corollary 5.7.** *Let $E = f(S)+C$, where $S \subseteq \mathbb{Z}^n$ is a Presburger set, $C$ a semilinear set, and $f : \mathbb{R}^m \to \mathbb{R}^n$ a linear map. Then $Zar(E)$.*

*Proof.* Since $f$ is linear, it equals $h + d$, where $h : \mathbb{R}^m \to \mathbb{R}^n$ is a linear map with $h(0) = 0$, and $d \in \mathbb{R}^n$. Hence, after replacing $C$ by $d + C$, we may assume that $f(0) = 0$. By Presburger cell decomposition, $S$ is a finite union of Presburger cells. By linear cell decomposition and Proposition A.2, $C$ is a finite union of product cells. Hence $E$ is a finite union of sets of the form $E' = f(S') + C'$ where $S'$ is a Presburger cell and $C'$ a product cell. By Remark 2.2, it suffices to prove $Zar(E')$ for each such set $E'$. But this is Proposition 5.6. □

5.2. **Sketch of proof for the general case (Theorem 5.21).** In the rest of Section 5, we reduce the general case to the sum case, as described in the beginning of the section and outlined further here. We use the notion of a line (Definition 4.39) in both structures $\langle \mathbb{R}, <, + \rangle$ and $\langle \mathbb{Z}, <, + \rangle$. To avoid ambiguity, we call a line in the former structure simply a *line*, and a line in the latter one a *Presburger line.* A line thus has form $L = p + d\mathbb{R}_{\geq 0} \subseteq \mathbb{R}^n$, for some $p \in \mathbb{R}^n$ and $d \in \mathbb{Q}^n$, and a Presburger line $l = p + d\mathbb{N} \subseteq \mathbb{Z}^n$, for some $p, d \in \mathbb{Z}^n$. If $l = p + d\mathbb{N}$ is a Presburger line, we call $L = p + d\mathbb{R}_{\geq 0}$ the *linear extension of* $l$.

*Remark* 5.8. Consider a Presburger cell $S$ and a linear cell $T$, with $S \subseteq T \subseteq \mathbb{R}^m$. Let $l \subseteq S$ be a Presburger line witnessing that $S$ is purely unbounded. Then the linear extension $L$ of $l$ witnesses that $T$ is purely unbounded (note that $L \subseteq T$, by convexity of $T$ ([25, Lemma 3.4]).

The proof of Theorem 5.21, namely that $Zar(X)$ for any definable $X$, spans across the next four subsections. The task carried out in each of them is outlined below.

- **Section 5.3.** We may assume $X = \bigcup_{t \in S} f(t) + Y_t$, where $f$ is a linear map, $S$ a Presburger cell, and $\{Y_t\}_{t \in S}$ a semilinear family with certain properties (Corollary 5.14).
- **Section 5.4.** There is a Presburger line $l$ containing 0, such that for every $p \in S$, $p + l$ witnesses that $S$ is purely unbounded (Proposition 5.15).
- **Section 5.5.** Such an $l$ is a nesting direction for $\{Y_t\}_{t \in S}$ (by Corollary 5.18).
- **Section 5.6.** If $X$ is $k$-free, then so is $f(S) + \bigcup_{t \in S} Y_t$ (Proposition 5.20). Hence we can conclude by the sum case, Corollary 5.7.

5.3. **Structure theorems for definable sets.** Our goal is to prove Corollary 5.14 below. The starting point is the following folklore result.

**Fact 5.9** (Structure Theorem)**.** *Let $X \subseteq \mathbb{R}^n$ be a definable set. Then:*

*(1) If $X \subseteq \mathbb{Z}^n$, then $X$ is a Presburger set.*
*(2) $X = \bigcup_{t \in S} X_t$, for some Presburger set $S \subseteq \mathbb{Z}^m$ and a semilinear family $\{X_t\}_{t \in \mathbb{R}^m}$ of sets $X_t \subseteq \mathbb{R}^n$.*

*Proof.* We provide a sketch of the proof. By [40, Appendix], the theory of $\langle \mathbb{R}, <, +, \mathbb{Z} \rangle$ admits quantifier elimination in the language $\{<, +, -, 0, 1, (\lambda_q)_{q\in\mathbb{Q}}, \lfloor \cdot \rfloor\}$, where $\lambda_q : x \mapsto qx$ denotes scalar multiplication by $\lambda_q$, and $\lfloor \cdot \rfloor : \mathbb{R} \to \mathbb{Z}$ the floor function $\lfloor x \rfloor = \max(\mathbb{Z} \cap (-\infty, x])$. By straightforward inductive arguments on the complexity of formulas, we can then prove the following two results:

(a) Every definable set $S \subseteq \mathbb{R}^m$ with $\mathsf{tdim} S = 0$ is of the form $f(Y)$, where $Y \subseteq \mathbb{Z}^l$ is a Presburger set and $f : \mathbb{R}^l \to \mathbb{R}^m$ is semilinear. Here, for $X \subseteq \mathbb{R}^m$, $\mathsf{tdim} X$ denotes the maximum $k \in \mathbb{N}$ such that some projection of $S$ onto $k$ coordinates has non-empty interior, if $S$ is non-empty, and $\mathsf{tdim}\emptyset = -\infty$. This result in particular shows that every definable subset of $\mathbb{Z}^l$ is a Presburger set, establishing (1).

(b) Every definable set $X \subseteq \mathbb{R}^n$ is a finite union of sets $Y$, each satisfying the following property: there are definable functions $f_1, \dots, f_k : \mathbb{R}^l \to \mathbb{R}$ and $g_1, \dots g_m : \mathbb{R}^l \to \mathbb{R}$, which are given by compositional iterates of basic functions, such that

$$Y = \{x \in \mathbb{R}^l \,:\, \forall i, j,\ f_i(x) = 0,\ g_j(x) > 0\}.$$

Here, by *basic functions* we mean either (a) semilinear functions, (b) the successor function $s : \mathbb{Z} \to \mathbb{Z}$, $x \mapsto x + 1$, or (c) the floor function $\lfloor \cdot \rfloor$.

Now, let $X \subseteq \mathbb{R}^n$ be a definable set. By (b) and [24, Proposition 5.2], there is a semilinear family $\{X'_t\}_{t\in\mathbb{R}^m}$ of subsets of $\mathbb{R}^n$ and a definable set $S' \subseteq \mathbb{R}^m$ of $\mathsf{tdim} S' = 0$, such that $X = \bigcup_{t\in S'} X'_t$. By (a), $S' = f(S)$ for some semilinear function $f$ and a Presburger set $S \subseteq \mathbb{Z}^m$. Thus, setting $X_t = X'_{f(t)}$, $t \in S$, we have

$$X = \bigcup_{t\in f(S)} X'_t = \bigcup_{t\in S} X_t,$$

as needed. ☐

Recall that we use '$\mathsf{dim}$' both in the o-minimal and Presburger settings, each use being clear from the context.

**Lemma 5.10.** *Let $C \subseteq \mathbb{Z}^m$ be Presburger cell. Then there is a linear cell $C' \subseteq \mathbb{R}^m$ with $\mathsf{dim} C = \mathsf{dim} C'$ and $C \subseteq C'$.*

*Proof.* By induction on $m$. For $m = 0$, $\mathbb{R}^0 = \{0\}$ and we let $C' = C$. Let $m > 0$, and $D = \pi(C) \subseteq \mathbb{Z}^{m-1}$. Let $D'$ be provided by Inductive Hypothesis for $D$. If $C = \Gamma(\alpha)$, let $C' = \Gamma(\alpha)_{D'}$, otherwise let $C' = (-\infty, \infty)_{D'}$. The result easily follows. ☐

We only need the next lemma for $T$ a linear cell.

**Lemma 5.11.** *Let $S \subseteq T \subseteq \mathbb{R}^m$, such that $S \subseteq \mathbb{Z}^m$ is a Presburger set and $T$ a semilinear set. Then there are Presburger cells $S_1, \dots, S_l$ and linear cells $T_1, \dots, T_l$, such that*

- *$S = S_1 \cup \cdots \cup S_l$, and*
- *for every $i = 1, \dots, l$, we have $S_i \subseteq T_i \subseteq T$ with $\mathsf{dim} S_i = \mathsf{dim} T_i$.*

*(We do not require $T = T_1 \cup \cdots \cup T_l$.)*

*Proof.* By induction on $r = \mathsf{dim} S$. For $r = 0$, $S$ is a singleton, and we can let $S_1 = T_1 = S$. Let $r > 0$. By Presburger cell decomposition, we may assume that $S$ is a Presburger cell. By Lemma 5.10, there is a semilinear set $U$ with $r = \mathsf{dim} U$ and $S \subseteq U$. By linear cell decomposition, $U \cap T$ is a finite union of linear cells $U_1, \dots, U_l$. By Fact 5.9(1), $S \cap U_j$ is a Presburger set, and by Presburger cell decomposition, it is a finite union of Presburger cells, $S_1, \dots, S_k$. It suffices to prove the statement for each $S_i$, $i = 1, \dots, k$, and $T \cap U$.

If $\mathsf{dim} S_i = r$, then $S_i \cap U_j \subseteq T \cap U_j \subseteq T \cap U$ have the same dimension, as needed.

If $\mathsf{dim} S_i < r$, then apply Inductive Hypothesis to $S_i$ and $T \cap U_j$, to obtain a decomposition of the former into Presburger cells, each contained in a linear cell of the same dimension and contained in $T \cap U_j \subseteq T \cap U$, as needed. □

**Fact 5.12.** *Let $\{X_t\}_{t\in\mathbb{R}^m}$ be a semilinear family of sets $X_t \subseteq \mathbb{R}^n$, such that $\bigcup_{t\in T}\{t\} \times X_t$ is a linear cell. Let $f : \mathbb{R}^m \to \mathbb{R}^l$ be a map, such that $f(T)$ is a linear cell. Then $\bigcup_{t\in f(T)}\{t\} \times X_t$ is a linear cell.*

*Proof.* By an easy induction on $n$, left to the reader. □

**Lemma 5.13.** *Every definable set $X \subseteq \mathbb{R}^n$ is a finite union of sets of the form $\bigcup_{t\in S} X_t$, where:*

- *$\{X_t\}_{t\in\mathbb{R}^m}$ is a semilinear family of sets $X_t \subseteq \mathbb{R}^n$,*
- *$S \subseteq \mathbb{Z}^m$ is a Presburger cell of dimension $m$,*
- *$S$ is contained in a linear cell $T \subseteq \mathbb{R}^m$ of dimension $m$, such that*
$$\bigcup_{t\in T}\{t\} \times X_t$$
*is a linear cell.*

*Proof.* By Fact 5.9, $X = \bigcup_{t\in S} X_t$, for some Presburger set $S \subseteq \mathbb{Z}^m$ and a semilinear family $\mathcal{X} = \{X_t\}_{t\in\mathbb{R}^m}$ of sets $X_t \subseteq \mathbb{R}^n$. By linear cell decomposition, the set $\bigcup_{t\in\mathbb{R}^m}\{t\} \times X_t$ is a finite union of linear cells $C_i$, $i = 1, \dots, l$, each of the form $C_i = \bigcup_{t\in T_i}\{t\} \times C^i_t$, for some linear cell $T_i \subseteq \mathbb{R}^m$ and $C^i_t \subseteq X_t$. We may therefore assume that $X = \bigcup_{t\in S} X_t$, for some Presburger set $S \subseteq \mathbb{Z}^m$ contained in a linear cell $T \subseteq \mathbb{R}^m$, such that $\bigcup_{t\in T}\{t\} \times X_t$ is also a linear cell. By Lemma 5.11, $S$ is a finite union of Presburger cells $S_i$, $i = 1, \dots, l$, such that each of them is contained in a linear cell $T_i \subseteq T$ with $\mathsf{dim} S_i = \mathsf{dim} T_i$. In particular, $\bigcup_{t\in S} X_t = \bigcup_i \bigcup_{t\in S_i} X_t$, and each $\bigcup_{t\in T_i}\{t\} \times X_t$ is a linear cell (Fact 5.12). Hence we may assume that $S$ is a Presburger cell.

Finally, we may assume that $\mathsf{dim} S = \mathsf{dim} T = m$, as follows. If $T$ is a linear $(i_1, \dots, i_m)$-cell, the projection $f : T \to \mathbb{R}^l$ onto the coordinates $j$ with $i_j = 0$ is injective ($l$ is the number of them). For $t \in f(T)$, let $Y_t = X_{f^{-1}(t)}$. Then $\bigcup_{t\in S} X_t = \bigcup_{t\in f(S)} Y_t$, and $\bigcup_{t\in f(T)}\{t\} \times Y_t$ is also a linear cell (Fact 5.12). Finally, $f(S)$ is easily seen to be a Presburger cell, hence we are done. □

**Corollary 5.14.** *Every definable subset of $\mathbb{R}^n$ is a finite union of sets of the form*
$$X = \bigcup_{t\in S} f(t) + Y_t,$$
*where $f : \mathbb{R}^m \to \mathbb{R}^n$ is a linear map, and $S$ a Presburger cell $S \subseteq \mathbb{Z}^m$ contained in a linear cell $T \subseteq \mathbb{R}^m$, such that:*

- *$\{Y_t\}_{t\in\mathbb{R}^m}$ is a semilinear family of sets in $\mathbb{R}^n$*
- *$\bigcup_{t\in T}\{t\} \times Y_t$ is a linear cell*
- *$\bigcap_{t\in T} Y_t \neq \emptyset$*
- *$\mathsf{dim} S = m$.*

*Proof.* By Lemmas 5.13 and A.8. □

5.4. **Lines witnessing purely unbounded Presburger cells.** We strengthen Remark 4.3 as follows.

**Proposition 5.15.** *Let $S \subseteq \mathbb{Z}^n$, $n > 0$, be a Presburger cell with* $\mathsf{dim} S = n$. *Then there is $d \in \mathbb{Z}^n$, such that for every $p \in S$, the Presburger line $p + d\mathbb{N}$ witnesses that $S$ is purely unbounded.*

*Proof.* We work by induction on $n$. For $n = 1$, without loss of generality, $S$ is of the form $[\alpha, \infty)^c_r$ or $(-\infty, \infty)^c_r$. Here, $\pi : \mathbb{Z} \to \mathbb{Z}^0 = \{0\}$, and hence for every $x \in S$, $S_{\pi(x)} = S$. It is then easy to check that for $d = r$ the statement holds.

Now let $n > 1$. Without loss of generality, $S = [\alpha, \beta]^{c_n}_{r_n}$, $[\alpha, \infty)^{c_n}_{r_n}$ or $(-\infty, \infty)^{c_n}_{r_n}$, where $0 \leq c_n < r_n$, and $\alpha, \beta$ are $\mathsf{Pres}$-linear maps. We handle the first two cases as the third one is clear. Suppose that

$$\alpha(x) = \sum_{i=1}^{n-1} a_i \left( \frac{x_i - c_i}{r_i} \right) + \gamma \quad \text{and} \quad \beta(x) = \sum_{i=1}^{n-1} b_i \left( \frac{x_i - c_i}{r_i} \right) + \delta.$$

for some $\gamma, \delta, a_i, b_i, c_i, r_i \in \mathbb{Z}$, with $0 \leq c_i < r_i$. By the definition of Presburger cells, $a_i \neq b_i$. Note also that for every $\hat{p} = (p_1, \ldots, p_{n-1}) \in \pi(S)$ and $\hat{d} = (d_1, \ldots, d_{n-1}) \in \prod_{i \in [r-1]} r_i\mathbb{Z}$, we have

$$\alpha(\hat{p} + \hat{d}) = \alpha(\hat{p}) + \sum_{i=1}^{n-1} a_i \frac{d_i}{r_i} \quad \text{and} \quad \beta(\hat{p} + \hat{d}) = \beta(\hat{p}) + \sum_{i=1}^{n-1} b_i \frac{d_i}{r_i}.$$

By Inductive Hypothesis, there is $\hat{d} = (d_1, \ldots, d_{n-1}) \in \mathbb{Z}^{n-1}$ that verifies the conclusion of the proposition for $\pi(S)$. Note that

$$\sum_{i=1}^{n-1} (b_i - a_i) d_i > 0. \tag{2}$$

Indeed, if not, for large enough $m \in \mathbb{N}$, the expression

$$(\beta - \alpha)(\hat{p} + m\hat{d}) = (\beta - \alpha)(\hat{p}) + \sum_{i=1}^{n-1} (b_i - a_i) \frac{m d_i}{r_i}$$

would become negative, contradicting the fact that $\hat{p} + \mathbb{N}\hat{d} \subseteq \pi(S)$.

We denote:

$$r := r_1 \ldots r_{n-1} \in \mathbb{N}_{>0}.$$

It is then easy to see that $2rr_n\hat{d}$ also satisfies the conclusion of the proposition for $\pi(S)$. Hence, we may assume that $\hat{d} \in 2rr_n\mathbb{Z}^{n-1}$, which we do in Case I below. In particular, $\hat{d} \in rr_n\mathbb{Z}^{n-1}$, which we use in Case II below. For $i = 1, \ldots, n-1$, denote:

$$r_{\neq i} = r_1 \ldots r_{i-1} r_i \ldots r_{n-1}.$$

**Case I: $S = [\alpha, \beta]^{c_n}_{r_n}$.** Recall that $d_i = 2rr_n z_i$, for some $z_i \in \mathbb{Z}$, $i = 1, \ldots, n-1$. Notice that since $d_i$ and $z_i$ have the same sign, (2) also yields

$$\sum_{i=1}^{n-1} (b_i - a_i) r_i r_n z_i > 0, \tag{3}$$

Now let

$$d = (\hat{d}, d_n), \text{ where } d_n = \sum_{i=1}^{n-1} (a_i + b_i) \frac{d_i}{2r_{\neq i}} = \sum_{i=1}^{n-1} (a_i + b_i) r_i r_n z_i.$$

Let $p = (\hat{p}, p_n) \in S$. We verify Definition 4.41(1), namely, that for every $m \in \mathbb{N}$, $p + md \in S$. It suffices to prove:

(a) $p_n + md_n \equiv c_n \operatorname{mod} r_n$, and
(b) $\alpha(\hat{p} + m\hat{d}) < p_n + md_n < \beta(\hat{p} + m\hat{d})$.

For (a), $p_n \equiv c_n \operatorname{mod} r_n$ (since $p \in S$) and

$$md_n = \sum_{i=1}^{n-1} m(a_i + b_i) r_i r_n z_i \equiv 0 \operatorname{mod} r_n.$$

For (b), it suffices to show that

$$\begin{aligned} &p_n + md_n - \alpha(\hat{p} + m\hat{d}) = p_n - \alpha(\hat{p}) + m\sum_{i=1}^{n-1}(b_i - a_i) r_i r_n z_i, \text{ and} \\ &\beta(\hat{p} + m\hat{d}) - (p_n + md_n) = \beta(\hat{p}) - p_n + m\sum_{i=1}^{n-1}(b_i - a_i) r_i r_n z_i, \end{aligned} \tag{4}$$

since then both quantities will be positive, by the fact that $\alpha(\hat{p}) < p_n < \beta(\hat{p})$ and (3). To verify the above equalities (say the former), observe that its left side equals:

$$p_n + md_n - \alpha(\hat{p}) - \sum_{i=1}^{n-1} a_i \frac{md_i}{r_i} = p_n - \alpha(\hat{p}) + m\sum_{i=1}^{n-1}(a_i + b_i) r_i r_n z_i - \sum_{i=1}^{n-1} a_i m 2 r_i r_n z_i,$$

from which the result follows.

Finally, Definition 4.41(2) also follows from (4) since for any $N \in \mathbb{N}$, the quantities in (4) become bigger than $N$ for sufficiently large $m \in \mathbb{N}$.

**Case II:** $S = [\alpha, \infty]^{c_n}_{r_n}$**.** Recall that $d_i = r r_n z_i$, for some $z_i \in \mathbb{Z}$, $i = 1, \dots, n-1$. Let

$$d = (\hat{d}, d_n), \text{ where } d_n = \sum_{i=1}^{n-1}(a_i + (-1)^{s_i}) \frac{d_i}{r_{\neq i}} = \sum_{i=1}^{n-1}(a_i + (-1)^{s_i}) r_i r_n z_i,$$

where $s_i = 0, 1$, so that $(-1)^{s_i} z_i > 0$. Let $p = (\hat{p}, p_n) \in S$. The proofs of Definition 4.41(1) and (2) are similar to those in Case I, after observing that (a) again holds, since

$$md_n = \sum_{i=1}^{n-1} m(a_i + (-1)^{s_i}) r_i r_n z_i \equiv 0 \operatorname{mod} r_n,$$

and, instead of (4),

$$p_n + md_n - \alpha(\hat{p} + m\hat{d}) = p_n + md_n - \alpha(\hat{p}) - \sum_{i=1}^{n-1} a_i \frac{md_i}{r_i} =$$

$$p_n - \alpha(\hat{p}) + m\sum_{i=1}^{n-1}(a_i + (-1)^{s_i}) r_i r_n z_i - \sum_{i=1}^{n-1} a_i m r_i r_n z_i = p_n - \alpha(\hat{p}) + \sum_{i=1}^{n-1} m r_i r_n (-1)^{s_i} z_i,$$

which is again positive and, for any $N \in \mathbb{N}$, it can grow bigger than $N$, for sufficiently large $m \in \mathbb{N}$. □

*Remark* 5.16. A common generalisation of Propositions 4.22 and 5.15 would be that there is a Presburger line $l \subseteq S$, such that for every $N \in \mathbb{N}$, we can find an $N$-internal point on $l$. We will not need such a statement.

5.5. **Nesting directions.** We need the following proposition from Appendix A.2.

**Proposition 5.17.** *Let $\{X_t\}_{t\in\mathbb{R}^m}$ be a semilinear family of sets in $\mathbb{R}^n$, such that $C=\bigcup_{t\in T}\{t\}\times X_t$ is a linear cell, and $\bigcap_{t\in T} X_t\neq\emptyset$. Suppose $L$ is a line witnessing $T$ is purely unbounded. Then $L$ is a nesting line for $\{X_t\}_{t\in T}$.*

*Proof.* By Proposition A.6, applied to $\langle\mathbb{R},<,+\rangle$ as a vector space over $\mathbb{Q}$. □

**Corollary 5.18.** *Let $\{X_t\}_{t\in\mathbb{R}^m}$ be a semilinear family of sets in $\mathbb{R}^n$, such that $\bigcup_{t\in T}\{t\}\times X_t$ is a linear cell, and $\bigcap_{t\in T} X_t\neq\emptyset$. Let $S\subseteq\mathbb{Z}^m$ be a Presburger cell with $S\subseteq T$ and $\mathsf{dim} S=m$. Then every $l=Im(\gamma)$ witnessing $S$ is purely unbounded is a nesting line for $\mathcal{X}=\{X_t\}_{t\in S}$. Moreover,*

$$\bigcup_{t\in S} X_t=\bigcup_{t\in T} X_t.$$

*Proof.* We first check conditions (1)-(2) of Definition 4.40. For (1), by Remark 5.8, the linear extension $L$ of $l$ witnesses $T$ is purely unbounded. By Proposition 5.17, $L$ is a nesting line for $\mathcal{X}'=\{X_t\}_{t\in T}$, so (1)-(2) hold for $L$ and $\mathcal{X}'$. This in particular implies (1) for $l$ and $\mathcal{X}$.

For (2), take any $t\in S\subseteq T$. By (2) for $L$ and $\mathcal{X}'$, there is $x\in L$ with $X_t\subseteq X_x$. Take $t'\in l\subseteq L$, with $x<t'$. Then by (1), $X_t\subseteq X_x\subseteq X_{t'}$, and we are done.

For the 'moreover' clause, let $t\in T$. The proof of the last paragraph applies identically to show that there is $t'\in l$ with $X_t\subseteq X_{t'}$, as needed. □

5.6. **Reduction to the sum case.** Here we conclude the proofs of the reduction and Theorem 5.21 below. We only need the following corollary for $l=m$.

**Corollary 5.19.** *Let $\mathcal{X}=\{X_t\}_{t\in S}$ be a family of sets $X_t\subseteq\mathbb{R}^n$, with $S\subseteq\mathbb{R}^m$, and assume that $\mathcal{X}$ has a nesting direction $L$. Then for every $t_1,\dots,t_l,s_1,\dots,s_m\in S$, there is $a\in S$, such that for every $i,j$, $X_{t_i}\subseteq X_{s_j+a}$.*

*Proof.* For every $i,j$, by Definition 4.40(2) applied to $l=s_j+L$, there is $a_{ij}\in L$, such that $X_{t_i}\subseteq X_{s_j+a_{ij}}$. Let $a=\max a_{ij}$ in the order of $L$. By Definition 4.40(1), applied to all of $s_j+L$, $X_{t_i}\subseteq X_{s_j+a}$, as needed. □

**Proposition 5.20.** *Let $\mathcal{X}=\{X_t\}_{t\in S}$ be a family of sets $X_t\subseteq\mathbb{R}^n$, with $S\subseteq\mathbb{R}^m$, and assume that $\mathcal{X}$ has a nesting direction. Let $f:\mathbb{R}^m\to\mathbb{R}^n$ be a linear map, and $k\in\mathbb{N}$. If $X=\bigcup_{t\in S} f(t)+X_t$ is $k$-free, then so is*

$$Y=f(S)+\bigcup_{t\in S} X_t.$$

*Proof.* It suffices to find, for every finite set $B\subseteq Y$, a set $B'\sim B$ in $X$. Let $B=\{b_1,\dots,b_l\}\subseteq Y$. Suppose that

$$b_i=f(s_i)+x_i\in f(s_i)+X_{t_i},$$

for some $s_i,t_i\in S$ and $x_i\in X_{t_i}$. It need not be $s_i=t_i$, but these two relate as follows: let $a\in S$ as in Corollary 5.19. Then for every $i=1,\dots,l$,

$$X_{t_i}\subseteq X_{s_i+a}.$$

We can therefore conclude as follows. Since $f:\mathbb{R}^m\to\mathbb{R}^k$ is linear, there is $c\in\mathbb{R}^k$, such that for every $i=1,\dots,l$, we have

$$f(a+s_i)=f(a)+f(s_i)+c.$$

We show that $B' = f(a) + c + B \subseteq X$. We have

$$f(a) + c + b_i = f(a) + c + f(s_i) + x_i \in f(a + s_i) + X_{t_i} \subseteq f(a + s_i) + X_{s_i + a} \subseteq X,$$

as needed. □

We can now conclude the main theorem of this section.

**Theorem 5.21.** *Let $X$ be a definable set in $\langle \mathbb{R}, <, +, \mathbb{Z}\rangle$. Then $Zar(X)$.*

*Proof.* Suppose $X$ is $k$-free. By Remark 2.2, we may assume that

$$X = \bigcup_{t\in S} f(t) + Y_t,$$

as in Corollary 5.14. Let $d$ be as in Proposition 5.15, and $l = \mathbb{N}d$. Then for every $p \in S$, $p + l$ witnesses that $S$ is purely unbounded. By Corollary 5.18, $p + l$ is a nesting line for $\{Y_t\}_{t\in S}$, and $\bigcup_{t\in S} Y_t = \bigcup_{t\in T} Y_t$. Hence $l$ is a nesting direction for $\{Y_t\}_{t\in S}$. By Proposition 5.20, $X' = f(S) + \bigcup_{t\in S} Y_t$ is $k$-free. Since $X \subseteq X'$, it suffices to show that $X'$ has linear $Z$-bounds for $\mathcal{C}$. But since $X' = f(S) + C$, with $S$ a Presburger set (in fact, cell) and $C = \bigcup_{t\in S} Y_t = \bigcup_{t\in T} Y_t$ a semilinear set, we are done by Corollary 5.7. □

We conclude this section with some applications of Theorem 5.21. We illustrate the first one with a simple example.

**Example 5.22.** Consider the family $\mathcal{E} = \{E_n\}_{n\in\mathbb{Z}}$ of semilinear sets, such that for every $n \in \mathbb{Z}$,

$$E_n = T_n + D,$$

where $T_n = \{(x, x) \in \mathbb{Z}^2 : 0 \leq x \leq n\}$ and $D = \{(x, 0) : x \in (0, 1)\}$ (as in Figure 1(2)). Each $E_n$ is 2-free, so we know it has linear $Z$-bounds, but Fact 1.3 does not guarantee that they are witnessed by the same $\alpha$ ($\mathcal{E}$ is *not* a subfamily of a semilinear family). However, thanks to Theorem 5.21, and since $\bigcup_{n\in\mathbb{Z}} \mathcal{E}$ is 2-free and definable in $\langle \mathbb{R}, <, +, \mathbb{Z}\rangle$, we can find such a common witness (of course, in this example it is easy to see directly that $\alpha = 1$).

In general, the following holds, by Theorem 5.21.

**Corollary 5.23.** *Let $\mathcal{M} = \langle \mathbb{R}, <, +\rangle$, and $\mathcal{E} = \{E_b\}_{b\in I}$ be a family of semilinear relations*

$$E_b \subseteq \prod_{i\in[r]} \mathbb{R}^{d_i},$$

*with $\bigcup \mathcal{E}$ definable in $\langle \mathbb{R}, <, +, \mathbb{Z}\rangle$ and $k$-free, some $k \in \mathbb{N}$. Then there is $\alpha \in \mathbb{R}_{>0}$, such that for every $b \in I$, $E_b$ has linear $Z$-bounds witnessed by $\alpha$.*

Our second application concerns Zarankiewicz's problem for some further ordered abelian groups.[3]

**Corollary 5.24.** *Let $\mathcal{M}$ be the direct sum of any finite number of copies of $\mathbb{R}$ and $\mathbb{Z}$, equipped with the lexicographic order. Then for every definable $E$, $Zar(E)$ holds.*

---

[3]We thank Joshua Losh for pointing out this application to us.

*Proof.* The structure $\mathcal{M}$ is 'definable' in $\langle \mathbb{R}, <, +, \mathbb{Z} \rangle$, namely, there is an isomorphism between $\mathcal{M}$ and a structure $\mathcal{M}'$ whose universe is the cartesian product of the corresponding copies of $\mathbb{R}$ and $\mathbb{Z}$, and whose atomic relations are all definable in $\langle \mathbb{R}, <, +, \mathbb{Z} \rangle$. Identifying $\mathcal{M}$ with $\mathcal{M}'$, we obtain that $E$ is also definable in $\langle \mathbb{R}, <, +, \mathbb{Z} \rangle$. Since $k$-freeness and having linear $Z$-bounds are preserved under this identification, the result follows. □

## 6. Abstract Zarankiewicz

In this section, we establish Theorem D (6.7), an abstract version of Zarankiewicz's problem, for an arbitrary 'saturated' structure $\mathcal{M}$, an operator $\mathsf{cl} : \mathcal{P}(M) \to \mathcal{P}(M)$, and a class $\mathcal{C}$ of grids satisfying certain properties, (DEF), (UB), and (TIGHT), which we introduce below. In the subsequent sections, we use Theorem D and its consequence Theorem 6.17 to derive several applications. We first recall any necessary model-theoretic background (standard references include [36, 46, 52]) and introduce our key properties.

6.1. **Model-theoretic background.** Let $\mathcal{M} = \langle M, \dots \rangle$ be a structure. Given a set $A \subseteq M$, a *type $p$ over $A$* is a set of formulas $\phi(x)$ with parameters from $A$ such that every finite subset of it has a common realisation in $\mathcal{M}$. We say that $p$ is *complete* if for every formula $\phi(x)$ with parameters from $A$, either $\phi(x)$ or $\neg\phi(x)$ is in $p$. The type of $a \in M^n$ over $A$, denoted by $\mathsf{tp}(a/A)$ is the set of all formulas with parameters from $A$ that $a$ realises, and it is complete. We call $\mathcal{M}$ $\kappa$-saturated if it contains a realisation for every type over a set of parameters of size $< \kappa$. By compactness, every structure $\mathcal{M}$ has an $|M|^+$-saturated extension. We call a set $A \subseteq M$ *small* if $\mathcal{M}$ is $|A|^+$-saturated. We call $\mathcal{M}$ *sufficiently saturated* if it is $|A|^+$-saturated for every parameter set $A$ we consider. Abusing terminology, for simplicity, we drop the word 'sufficiently'. By a *type-definable* set we mean the set of realisations of a type over a small set. In particular, a definable set is also type-definable.

If $a = (a_1, \dots, a_n) \in M^n$, we let $\hat{a} = \{a_1, \dots, a_d\}$, and if $A \subseteq M^n$, we let $\hat{A} = \bigcup\{\hat{a} : a \in A\}$. If $\mathsf{cl} : \mathcal{P}(M) \to \mathcal{P}(M)$ is an operator, by $a \in \mathsf{cl}(A)$ we mean that each $a_i \in \mathsf{cl}(\hat{A})$. If $B \subseteq M^n$ and $A \subseteq M$, we call $B$ $\mathsf{cl}$-*independent over $A$* if for every $b \in B$, $b \notin \mathsf{cl}(A \cup (\hat{B} \setminus \hat{b}))$. Finally, by an $A$-(type-)definable set we mean a set which is $\hat{A}$-(type-)definable.

Note that at this stage we do not assume that $\mathsf{cl}$ has any additional properties (such as being a pregeometry) – this makes our setup applicable to more settings, such as that of stable 1-based theories in Section 7.1 with $\mathsf{cl} = \mathsf{acl}$.

6.2. **Key properties.** We let $\mathcal{C}$ be a class of grids in $\prod_{i \in [r]} M^{d_i}$, and $\mathsf{cl} : \mathcal{P}(M) \to \mathcal{P}(M)$ an operator. Our first property, definability (DEF), expresses the fact that $\mathsf{cl}$ is witnessed by a definable set.

(DEF) Let $(a, b) \in M^{k+m}$ and $A \subseteq M$. If $a \in \mathsf{cl}(Ab)$, then there is an $A$-definable set $X \subseteq M^{k+m}$ with $(a, b) \in X$, and such that for every $(a', b') \in X$, $a' \in \mathsf{cl}(Ab')$.

**Example 6.1.** For any structure $\mathcal{M}$ and $\mathsf{cl} = \mathsf{acl}$, (DEF) holds.

Our second property, uniform bounds (UB), is a uniform finiteness statement relativised to the class $\mathcal{C}$. It imposes a uniform bound on the size of grids from $\mathcal{C}$ intersected with $A$-definable sets contained in $\mathsf{cl}(A)$. In the geometric applications below (Section 8), (UB) forces $\mathcal{C}$ to contain only grids that are 'cl-far apart'; more precisely, '$\mathbb{Z}$-distant' (Definition 1.7), or 'tall' (Definition 8.7), in saturated models of Pres, or semibounded structures, respectively.

(UB) Let $A \subseteq M$ be small, $i \in [r]$ and $\{X_b\}_{b\in I}$ an $A$-definable family of subsets of $M^{d_i}$. Then there is $N \in \mathbb{N}$ such that for every $b \in I$ with $X_b \subseteq \mathsf{cl}(Ab)$, and $B \in \mathcal{C}$, we have $|X_b \cap B_i| \leq N$.

**Example 6.2.** Let $\mathsf{cl} = \mathsf{acl}$, and $\mathcal{C}$ the class of all grids. For any structure whose theory eliminates $\exists^\infty$, (UB) holds. For models of Presburger arithmetic, it fails, as witnessed by the family $\{0 < x < b : b \in M\}$. We show, however, in Proposition 8.2 that it holds with $\mathcal{C}$ the class of all $\mathbb{Z}$-distant grids.

We note that (UB) also trivially holds for any structure $\mathcal{M}$, cl and $\mathcal{C}$ the class of all cl-independent grids, defined next.

**Definition 6.3.** Let $A \subseteq M$. A grid $B$ is called cl-*independent over $A$* if each $B_i$ is cl-independent over $A$.

Finally, tightness (TIGHT) is an abstract version of [4, Lemma 5.5], capturing a 'linearity' property (Example 6.5). It provides some uniform dependence on the tuples of an $\infty$-free $E$.

**Definition 6.4.** Let $E \subseteq \prod_{i\in[r]} M^{d_i}$ and $A \subseteq M$ a set. We say that $E$ is cl-*tight over $A$* if for every $a = (a_1, \ldots, a_r) \in E$, the set $\{a_1, \ldots, a_r\}$ is cl-dependent over $A$.

(TIGHT) Let $A \subseteq M$ be small. Let $E \subseteq \prod_{i\in[r]} M^{d_i}$ be an $A$-type-definable relation that is $\mathcal{C}$-$\infty$-free. Then $E$ is cl-tight over $A$.

**Example 6.5.** We will show that (TIGHT) holds in the stable 1-based context (Proposition 7.7), as well as the weakly locally modular one (Proposition 6.15). It fails for $\mathsf{cl} = \mathsf{acl}$ in the presence of a field: let $E$ be the point-line incidence relation $E = \{(a, b, c, d) \in \mathbb{R}^4 : d = ac + b\}$, definable in the real field. Then $E$ is 2-free (hence $\infty$-free), but for algebraically independent $a, b, c$, the transcendence degree of $\{a, b, c, d\}$, where $d = ac + b$, is 3, and hence neither $\{a, b\} \subseteq \mathsf{acl}(c, d)$ nor $\{c, d\} \subseteq \mathsf{acl}(a, b)$.

*For the rest of Section 6, and unless stated otherwise, we fix a structure $\mathcal{M}$, an operator $\mathsf{cl} : \mathcal{P}(M) \to \mathcal{P}(M)$, and a class $\mathcal{C} \subseteq \prod_{i\in[r]} M^{d_i}$ of grids.*

6.3. **Results.** We sketch the role of our key properties in the proof of Theorem 6.7 below, for a single set $E$. Recall that for every $i = 1, \ldots, r$,

$$p_i : \prod_{j\in[r]} M^{d_j} \to \prod_{j\neq i} M^{d_j}$$

denotes the projection onto all but the $i$-th block of coordinates. So assume $E$ is $\mathcal{C}$-$\infty$-free. By (TIGHT), for every tuple in $E$, one of the blocks of coordinates is in the closure of the others. By (DEF) and compactness, we may assume it is always the $r$-th block of coordinates of tuples in $E$ that is in the closure of the first $r-1$ ones. By (UB), the intersection of the $r$-th component $B_r$ of each grid $B$ with $E$ has a uniform bound $N$, and hence the projection $p_{r\restriction E\cap B}$ is $N$-to-1. Thus we can conclude by Lemma 2.1. Note that the majority of this argument is realised in the proof of Key Lemma 6.6. That proof follows the reasoning in [4, Theorem 5.6], and has its roots in [26, Proposition 3.1].

**Lemma 6.6** (Key Lemma)**.** *Suppose $\mathcal{M}$ is saturated, satisfying* (DEF) *and* (UB)*. Let $A \subseteq M$. Let $\{E_b\}_{b\in I}$ be an $A$-type-definable family of relations $E_b \subseteq \prod_{i\in[r]} M^{d_i}$, with $I \subseteq M^m$. Then there is $\alpha \in \mathbb{R}_{>0}$, such that for every $b \in I$, if $E_b$ is $\mathsf{cl}$-tight over $Ab$, then $E_b$ has linear $Z$-bounds for $\mathcal{C}$ witnessed by $\alpha$.*

*Proof.* Let $b \in I$ such that $E_b$ is $\mathsf{cl}$-tight over $Ab$. Then for every $a = (a_1, \dots, a_r) \in E_b \cap \prod_{i\in[r]} M^{d_i}$, the set $\{a_1, \dots, a_r\}$ is $\mathsf{cl}$-dependent over $Ab$, say $a_i \in \mathsf{cl}(a_{\neq i}bA)$. By (DEF), there is an $A$-definable set $X$ that contains $(a, b)$, such that for every $(a', b') \in X$, we have $a'_i \in \mathsf{cl}(a'_{\neq i}b'A)$. By a standard compactness argument, $E_b$ is contained in a finite union of $A$-definable sets $(E_b)_i$, $i = 1, \dots, r$, such that for every $i = 1, \dots, r$ and $a \in (E_b)_i$, we have $a_i \in \mathsf{cl}(a_{\neq i}bA)$. Thus

$$((E_b)_i)_{a_{\neq i}} \subseteq \mathsf{cl}(a_{\neq i}bA).$$

Now fix some $i \in [r]$. Denote $L = \prod_{j\neq i} M^{d_j}$. Then for each $c \in L$, $(E_b)_c \subseteq M^{d_i}$. Now, by (UB), applied to the family $\{(E_b)_c\}_{b\in I, c\in L}$, we obtain $N_i \in \mathbb{N}$ such that for every $b \in I$ and $c \in L$, with

$$(E_b)_c \subseteq \mathsf{cl}(Abc)$$

and $Y \in \mathcal{C}$, we have $|(E_b)_c \cap Y_i| \leq N_i$. In particular, for every $i = 1, \dots, r$ and $B \in \mathcal{C}$, the map $p_{i\restriction (E_b)_i \cap B}$ is $N_i$-to-1. By Lemma 2.1, $(E_b)_i$ has linear $Z$-bounds for $\mathcal{C}$, witnessed by $N$. Hence $\alpha = \sum_i N_i$ witnesses that $E_b$ has linear $Z$-bounds for $\mathcal{C}$, as needed. ☐

We can now obtain the following abstract version of [4, Theorem 5.6].

**Theorem 6.7.** *Let $\mathcal{M}$ be saturated, $\mathsf{cl}$ an operator, and $\mathcal{C}$ a class of grids. Suppose* (DEF)*,* (UB) *and* (TIGHT)*. Let $\mathcal{E} = \{E_b\}_{b\in I}$ be a type-definable family of relations $E_b \subseteq \prod_{i\in[r]} M^{d_i}$. Then $Zar^\infty(\mathcal{E}, \mathcal{C})$ holds.*

*Proof.* Suppose $\{E_b\}_{b\in I}$ is $A$-type-definable. Let $\alpha$ be as in Key Lemma (6.6). Let $b \in I$ be such that $E_b$ is $\mathcal{C}$-$\infty$-free. By (TIGHT), $E_b$ is $\mathsf{cl}$-tight over $Ab$. Hence, by the Key Lemma, $E_b$ has linear $Z$-bounds for $\mathcal{C}$ witnessed by $\alpha$, as needed. ☐

6.4. **Weakly locally modular pregeometries.** In this subsection, we focus on abstract Zarankiewicz when $\mathsf{cl}$ is a pregeometry satisfying the usual 'independence axioms', as well as 'weak local modularity'. Our Theorem 6.17 below extends [4, Theorem 5.6] in that we do not assume that $\mathcal{M}$ eliminates $\exists^\infty$ or that $\mathsf{cl} = \mathsf{acl}$. This allows us to apply it in Section 7 to the setting of regular types and in Section 8 to the semibounded and Presburger settings (for appropriate $\mathsf{cl}$ and $\mathcal{C}$).

Let us recall the notion of a pregeometry, for any set $M$ and an operator $\mathsf{cl} : \mathcal{P}(M) \to \mathcal{P}(M)$.

**Definition 6.8.** We call $(M, \mathsf{cl})$ a *pregeometry* if for all $A, B \subseteq M$ and $a, b \in M$:

(1) $A \subseteq \mathsf{cl}(A)$
(2) $\mathsf{cl}\big(\mathsf{cl}(A)\big) = \mathsf{cl}(A)$
(3) $\mathsf{cl}(A) = \cup\{\mathsf{cl}(B) : B \subseteq A \text{ finite}\}$
(4) (Exchange) $a \in \mathsf{cl}(bA) \setminus \mathsf{cl}(A) \Rightarrow b \in \mathsf{cl}(aA)$.

**Example 6.9.** If $\mathcal{M}$ is an o-minimal structure or a model of Presburger arithmetic, then $(M, \mathsf{acl})$ is a pregeometry.

We next recall the notion of $\mathsf{cl}$-independence between two sets.

**Definition 6.10.** Given $A, B, C \subseteq M^n$, we call $A$ $\mathsf{cl}$-*independent from* $B$ *over* $C$, denoted $A \mathop{\smile\!\!\!\!\!\!|}^{\mathsf{cl}}_{C} B$, if every finite $A_0 \subseteq A$ that is $\mathsf{cl}$-independent over $C$ is also $\mathsf{cl}$-independent over $BC$.

We record here some types of pregeometries which will be relevant in Section 7.

**Definition 6.11.** We say that a pregeometry $(M, \mathsf{cl})$ is

(1) *trivial* if for all $A \subseteq M$, we have that $\mathsf{cl}(A) = \bigcup_{a \in A} \mathsf{cl}(a)$.
(2) *modular* if for all $A, B \subseteq M$,

$$A \mathop{\smile\!\!\!\!\!\!|}\limits^{\mathsf{cl}}_{\mathsf{cl}(A) \cap \mathsf{cl}(B)} B$$

(3) *locally modular* if for some $x \in M$, the localisation $(M, \mathsf{cl}(-x))$ is modular.

We are interested below in the case when $\mathop{\smile\!\!\!\!\!\!|}^{\mathsf{cl}}$ satisfies a list of 'independence axioms'. Some of them can be stated for any ternary relation $\mathop{\smile\!\!\!\!\!\!|} \subseteq \mathcal{P}(M)^3$, as follows (for all $a, a', b, b' \in M^n$ and $C \subseteq M$):

- (Extension) If $a \mathop{\smile\!\!\!\!\!\!|}_C b$, then for all $d \in M$ there is some $a' \in M$ such that $\mathsf{tp}(a'/Cb) = \mathsf{tp}(a/Cb)$ and $a' \mathop{\smile\!\!\!\!\!\!|}_C bd$.
- (Monotonicity) $aa' \mathop{\smile\!\!\!\!\!\!|}_C bb' \;\Rightarrow\; a \mathop{\smile\!\!\!\!\!\!|}_C b$.
- (Symmetry) $a \mathop{\smile\!\!\!\!\!\!|}_C b \;\Leftrightarrow\; b \mathop{\smile\!\!\!\!\!\!|}_C a$.
- (Transitivity) $a \mathop{\smile\!\!\!\!\!\!|}_C bb' \;\Leftrightarrow\; a \mathop{\smile\!\!\!\!\!\!|}_C b$ and $a \mathop{\smile\!\!\!\!\!\!|}_{Cb} b'$.

**Definition 6.12.** We say that $\mathop{\smile\!\!\!\!\!\!|}^{\mathsf{cl}}$ *satisfies the independence axioms* if, for all $a, a', b, b' \in M^n$ and all $C \subseteq M$:

(1) Extension, Monotonicity, Symmetry, Transitivity hold
(2) $a \mathop{\smile\!\!\!\!\!\!|}_C b \;\Leftrightarrow\; \mathsf{cl}(aC) \mathop{\smile\!\!\!\!\!\!|}_C \mathsf{cl}(bC)$
(3) (Non-degeneracy) If $a \mathop{\smile\!\!\!\!\!\!|}_C b$ and $d \in \mathsf{cl}(aC) \cap \mathsf{cl}(bC)$ then $d \in \mathsf{cl}(C)$.

**Example 6.13.** It is well-known, and easy to check, that for any structure $\mathcal{M}$, with $(M, \mathsf{acl})$ a pregeometry, $\mathop{\smile\!\!\!\!\!\!|}^{\mathsf{acl}}$ satisfies the independence axioms.

Finally, we recall the following 'linearity property' from [7].

(WLM) For all small $A, B \subseteq M$ there is some small $C \subseteq M$ such that

$$C \mathop{\smile\!\!\!\!\!\!|}\limits^{\mathsf{cl}}_{\emptyset} AB \quad \text{and} \quad A \mathop{\smile\!\!\!\!\!\!|}\limits^{\mathsf{cl}}_{\mathsf{cl}(AC) \cap \mathsf{cl}(BC)} B$$

**Example 6.14.** An o-minimal structure is linear if and only if (WLM) holds for $\mathsf{acl}$ ([7, Propositions 4.8 & 6.9 and Theorem 6.10(4)]). We will show in Lemma 8.6 that a semibounded structure satisfies (WLM) for the 'short closure operator' $\mathsf{scl}$.

The following proposition is an abstract version of [4, Lemma 5.5].

**Proposition 6.15.** *Assume* $\mathop{\smile\!\!\!\!\!\!|}^{\mathsf{cl}}$ *satisfies the independence axioms and* (WLM). *Let* $\mathcal{C}$ *be a class of grids that contains, for all small* $A \subseteq M$, *all grids that are* $\mathsf{cl}$-*independent over* $A$. *Then* (TIGHT) *holds.*

*Proof.* Assume that $E \subseteq \prod_{i \in [r]} M^{d_i}$ is $A$-type-definable and contains no infinite grid in $\mathcal{C}$. By the assumption on $\mathcal{C}$, $E$ contains no infinite $\mathsf{cl}$-independent (over $A$) grid. We want to prove that $E$ is $\mathsf{cl}$-tight over $A$. The proof is word-by-word the same with that of [4, Lemma 5.5], after replacing $\mathsf{acl}$ there by $\mathsf{cl}$ (and $b$ by $A$). The proof begins assuming, for a contradiction, the existence of an $\mathsf{acl}$-independent (over $b$) tuple $a \in E$ ($\mathsf{cl}$-independent over $A$ here), and continues to find an infinite grid contained in $E$. As noted in the second to last paragraph of that proof, the infinite grid is also $\mathsf{acl}$-independent over $b$ ($\mathsf{cl}$-independent over $A$ here). A contradiction.

Note that the $A$-definability of $E$ in [4, Lemma 5.5] is not used in its proof, except for the last paragraph, to ensure that two tuples with the same type over $A$ are either both in $E$ or both not in $E$ – this is also true for our $A$-type-definable $E$. Finally, the finiteness of the tuple $b$ is not used at all there , nor is the assumption that $\mathcal{M}$ eliminates $\exists^\infty$. □

*Remark* 6.16. We can replace (TIGHT) in the conclusion of Proposition 6.15 by the stronger 'if an $A$-type-definable $E$ contains no infinite grid $\mathsf{cl}$-independent over $A$, then it is $\mathsf{cl}$-tight over $A$', as the proof shows.

We can conclude Zarankiewicz's problem in the weakly locally modular setting.

**Theorem 6.17.** *Let $\mathcal{M}$ be saturated, $\mathsf{cl}$ an operator, such that $\mathop{\perp\!\!\!\perp}^{\mathsf{cl}}$ satisfies the independence axioms and* (WLM)*, and $\mathcal{C}$ a class of grids containing every grid that is $\mathsf{cl}$-independent over some small set. Suppose* (DEF) *and* (UB)*. Let $\{E_b\}_{b\in I}$ be a type-definable family of relations $E_b \subseteq \prod_{i\in[r]} M^{d_i}$. Then $Zar^\infty(\mathcal{E},\mathcal{C})$ holds.*

*Proof.* By Proposition 6.15, (TIGHT) holds. We conclude by Theorem 6.7. □

*Remark* 6.18. We can replace $Zar^\infty_\alpha(E_b,\mathcal{C})$ in the conclusion of Theorem 6.17 (inside $Zar^\infty(\mathcal{E},\mathcal{C})$) by the stronger 'if $E_b$ does not contain any infinite grid $\mathsf{cl}$-independent over $Ab$, then it has linear $Z$-bounds for $\mathcal{C}$ witnessed by $\alpha$'. Indeed, by Remark 6.16, such an $E_b$ is $\mathsf{cl}$-tight over $Ab$, and the proof of the Theorem 6.7 still goes through.

6.5. **A general remark.** In order to utilise the material of this section in what follows, we need the following remark.

*Remark* 6.19. Let $\mathcal{M}$ be a structure, and $\alpha \in \mathbb{R}_{>0}$.

(1) Assume $\mathcal{M}$ is saturated, and let $E$ be a type-definable relation. Then, by compactness, $E$ is $k$-free for some $k \in \mathbb{N}$ if and only if it is $\infty$-free. In particular, $Zar_\alpha(E)$ is equivalent to $Zar^\infty_\alpha(E)$.
(2) Assume $\mathcal{M}$ is saturated. Then, by compactness, if $Zar^\infty(\mathcal{E})$ (or $Zar(\mathcal{E})$) holds for every definable family $\mathcal{E}$, then so does it for every type-definable relation $\mathcal{E}$.
(3) Let $\mathcal{N}$ be a saturated elementary extension of $\mathcal{M}$, $E$ a type-definable relation, and $E^{\mathcal{N}}$ the interpretation of $E$ in $\mathcal{N}$. Then for a fixed $k \in \mathbb{N}$, $E$ is $k$-free if and only if $E^{\mathcal{N}}$ is. In particular, if $\mathcal{E}$ is a type-definable family of relations, then $Zar(\mathcal{E}^{\mathcal{N}})$ implies $Zar(\mathcal{E})$.

As a note, only $\aleph_1$-saturation is needed in (1) and (2), and $|T|^+$-saturation in (3).

## 7. Model-theoretic versions of Zarankiewicz's problem

In this section, we use Theorems 6.7 and 6.17 to obtain two new model-theoretic versions of Zarankiewicz's problem: in a model of a *stable* 1*-based theory* (Section 7.1), and in the induced structure on a *locally modular regular type* in a stable theory (Section 7.3). Along the way, we deduce Zarankiewicz's problem for *ab initio* Hrushovski constructions (Section 7.2), which are neither pregeometric nor 1-based.

We assume some familiarity with the basic notions of stability theory (as can be found, for example, in [46] or [52]), although we recall the main definitions as we go along. An $\mathcal{L}$-theory $T$ is called *stable* if for all $\mathcal{L}$-formulas $\phi(x,y)$, with $|x|=|y|$, there is $k\in\mathbb{N}$, such that for all $\mathcal{M}\vDash T$, there is no sequence $(a_i : i<k)$ such that

$$\mathcal{M}\vDash \phi(\bar{a}_i,\bar{a}_j) \text{ if and only if } i<j,$$

for all $i,j<k$. Given an $\mathcal{L}$-theory $T$, the theory $T^{\mathsf{eq}}$ is a theory in an expanded language by formulas for quotients of $\emptyset$-definable equivalence relations in $T$ (satisfying appropriate axioms) and the canonical projection functions. In particular, if $\mathcal{M}\models T$, then there is a model $\mathcal{M}^{\mathsf{eq}}\models T^{\mathsf{eq}}$ such that $M$ and every $\emptyset$-definable subset of $M$ is $\emptyset$-definable in $\mathcal{M}^{\mathsf{eq}}$. Since Zarankiewicz statements are preserved when passing from $\mathcal{M}^{\mathsf{eq}}$ to $\mathcal{M}$, we often work in $T^{\mathsf{eq}}$, and write $\mathsf{acl}$ for $\mathsf{acl}^{\mathsf{eq}}$.

By $\perp\!\!\!\!\perp$ we denote *non-forking* independence. It is well-known (for example, [52, Theorem 8.5.5]) that in stable theories it satisfies monotonicity, symmetry, transitivity from Section 6.4 (and a certain converse holds – [9, Fact 10.4.5])

7.1. **Stable 1-based Zarankiewicz.** The property of 1-basedness can be viewed as an analogue of local modularity, outside pregeometric contexts.

**Definition 7.1.** We say that a stable theory $T$ is

(1) *1-based* if for all $\mathcal{M}\vDash T$ and small sets $A,B\subseteq M$, we have $A\perp\!\!\!\!\perp_{A\cap B} B$.
(2) *without the finite cover property* ($\mathsf{nfcp}$) if $T^{\mathsf{eq}}$ eliminates $\exists^\infty$.

Our main theorem of this subsection is the following.

**Theorem 7.2.** *Let $T$ be a stable 1-based theory with* $\mathsf{nfcp}$*, and $\mathcal{M}\models T$. Let $\mathcal{E}=\{E_b\}_{b\in I}$ be a type-definable family of relations $E_b\subseteq\prod_{i\in[r]}M^{d_i}$. Then $Zar(\mathcal{E})$ (and $Zar^\infty(\mathcal{E})$, if $\mathcal{M}$ is saturated) holds.*

A non-parametric version of Theorem 7.2, with $E$ binary, was proved by Evans [26]. It was phrased there as follows.

**Fact 7.3** ([26, Proposition 3.1])**.** *Let $T$ be a complete, stable* 1*-based theory with* $\mathsf{nfcp}$*. Then, in any saturated $\mathcal{M}\vDash T$, any type-definable pseudoplane in $T$ is sparse.*

The notion of *sparseness* is identical to that of having linear $Z$-bounds for the class of all grids. A *pseudoplane* $E$ is a binary relation such that

(1) for every $(a,b)\in E$, both $E_a$ and $E_b$ are infinite, and
(2) for every $a\neq a'$ in the domain of $E$, $E_a\cap E_{a'}$ is finite
(3) for every $b\neq b'$ in the range of $E$, $E^b\cap E^{b'}$ is finite.

So, under (1), and recast in our terminology, we have:

$$E \text{ is 2-free } \Rightarrow\ E \text{ is a pseudoplane} \Rightarrow\ E \text{ is } \infty\text{-free}.$$

Thus, Theorem 7.2 strengthens Fact 7.3 in that $r\geq 2$, and $E$ is $k$-free for some $k\in\mathbb{N}$ (and $\infty$-free, for $\mathcal{M}$ saturated). Moreover, Theorem 7.2 is a parametric statement.

To apply Theorem 6.7, with $\mathsf{cl}=\mathsf{acl}$ in $\mathcal{M}\models T$, and $\mathcal{C}$ the class of all grids, we need to ensure (DEF), (UB) and (TIGHT). The former two properties being straightforward, we focus on (TIGHT) (Proposition 7.7).

We recall some further facts from stability theory. Let $\mathcal{M}$ be a structure.

**Definition 7.4.** Let $p$ be a complete type over $\mathcal{M}$. We say that $p$ is *stationary* if for every saturated elementary extension $\mathcal{N} \succcurlyeq \mathcal{M}$, $p$ has a unique non-forking extension to a type over $\mathcal{N}$.

In particular, if $p$ is a stationary type over $A$, then for any $B \supseteq A$ inside any elementary extension of $\mathcal{M}$, $p$ has a unique non-forking extension to a type over $B$, which we denote $p|_B$.

**Definition 7.5.** Let $p$ and $q$ be stationary types over $\mathcal{M}$. We define the *Morley product* $p \otimes q$ of $p$ and $q$ as follows. Take a saturated elementary extension $\mathcal{N} \succcurlyeq \mathcal{M}$, a tuple $b$ in a possibly larger elementary extension of $\mathcal{N}$, with $b \vDash q|_N$ and $a \vDash p|_{Nb}$. Then $p \otimes q = \mathsf{tp}(a, b/M)$.

It follows from the definitions that whenever $a \mathop{\smile\!\!\!\!\!\!\mid}_D b$, $p = \mathsf{tp}(a/D)$ and $q = \mathsf{tp}(b/D)$, then $\mathsf{tp}(a, b/D) = p \otimes q$. It is also well-known that the Morley product is associative, that is $p \otimes (q \otimes r) = (p \otimes q) \otimes r$. Given this, if $p$ is a stationary type, we define $p^{(1)} := p$ and $p^{(n)} := p \otimes p^{(n-1)}$ for $n > 1$. We define also $p^{(\omega)} := \bigcup_{n\in\omega} p^{(n)}$. A sequence realising $p^{(\omega)}$ is called a *Morley sequence in p*. It easily follows from [46, Lemma I.2.28] that if $\mathcal{M}$ is a model of a stable theory, then $\otimes$ is commutative.

**Definition 7.6.** A grid $B$ is called a *Morley grid* if each $B_i$ is a subset of a Morley sequence in some stationary type $s_i$. If all $s_i = p$, we call $B$ a *Morley grid in p*. We write $\mathcal{C}_R$ ($\mathcal{C}_{R,p}$) for the class of all Morley grids (in $p$).

**Proposition 7.7.** *Let $T = T^{\mathsf{eq}}$ be a stable 1-based theory, $\mathcal{M} \models T$ saturated, $\mathsf{cl} = \mathsf{acl}$ in $\mathcal{M}$, and $\mathcal{C}$ the class of all grids. Then* (TIGHT) *holds. In fact, every $\mathcal{C}_R$-$\infty$-free $A$-type-definable relation $E \subseteq \prod_{i\in[r]} M^{d_i}$ is $\mathsf{acl}$-tight over $A$.*

*Proof.* Let $E$ be the set of realisations of a type $\pi(x_1, \ldots, x_r)$ over $A$. To ease notation, assume $A = \emptyset$. Suppose that there is some realisation $a = (a_1, \ldots, a_r) \vDash \pi(x)$ such that $a_i \notin \mathsf{acl}(a_{\neq i})$ for all $i \in [r]$. We want to show that $E$ is not $\infty$-free.

Since $T$ is 1-based, we know that:

$$\mathsf{acl}(a_i) \mathop{\smile\!\!\!\!\!\!\mid}_{\mathsf{acl}(a_i)\cap\mathsf{acl}(a_{\neq i})} \mathsf{acl}(a_{\neq i}),$$

for all $i \in [r]$. Let $D = \bigcap_{i\in[r]} \mathsf{acl}(a_{\neq i})$, and observe that for all $i \in [r]$ we have that $\mathsf{acl}(a_i) \cap \mathsf{acl}(a_{\neq i}) \subseteq D \subseteq \mathsf{acl}(a_{\neq i})$, so, by monotonicity and transitivity, we have

$$a_i \mathop{\smile\!\!\!\!\!\!\mid}_D a_{\neq i},$$

for all $i \in [r]$.

Let $s_i(x_i) = \mathsf{tp}(a_i/D)$, and observe that $D = \mathsf{acl}(D)$, since it is the intersection of algebraically closed sets. Since $D$ is algebraically closed in $\mathcal{M}^{\mathsf{eq}}$ and $T$ is stable, each $s_i(x_i)$ is stationary, $i \in [r]$ ([52, Corollary 8.5.3]). Thus, their Morley product is defined.

Note that since $(a_1, \ldots, a_r)$ is a tuple which is forking-independent over $D$, we obtain

$$\mathsf{tp}(a_1, \ldots, a_r/D) = s_1(x) \otimes \cdots \otimes s_r(x).$$

Now, take a Morley sequence in $s_1(x) \otimes \cdots \otimes s_r(x)$, say:

$$B := ((a_1^t, \ldots, a_r^t) : t \in \omega) \vDash (s_1(x) \otimes \cdots s_r(x))^{\otimes\omega}$$

Each coordinate $B_i := (a_i^t : t \in \omega)$ is a Morley sequence in $s_i(x)$, and since by assumption these types are non-algebraic, each $B_i$ is infinite. Finally, since $T$ is stable, $\otimes$ is commutative, so for all $(a_1^{t_1}, \dots, a_r^{t_r}) \in \prod_{i \in [r]} B_i$ we have that:

$$(a_1^{t_1}, \dots, a_r^{t_r}) \models \mathsf{tp}(a_1, \dots, a_r/D).$$

In particular, $\prod_{i \in [r]} B_i \subseteq E$, implying that $E$ is not $\mathcal{C}_R$-$\infty$-free. $\square$

We can now conclude the proof of our theorem.

*Proof of Theorem 7.2.* Without loss of generality, we may assume $\mathcal{M} = \mathcal{M}^{\mathsf{eq}}$. (DEF) is by Example 6.1, and (UB) by $\mathsf{nfcp}$. Assume $\mathcal{M}$ is saturated. By Proposition 7.7, (TIGHT) also holds. Hence we conclude by Theorem 6.7 and Remark 6.19(3). $\square$

*Remark* 7.8. As our proof of Theorem 7.2 shows, in $Zar(\mathcal{E})$ and $Zar^{\infty}(\mathcal{E})$, we may replace the assumption that $E_b$ is $k$-free (or $\infty$-free), by the weaker one that $E_b$ is $\mathcal{C}_R$-$k$-free (or $\mathcal{C}_R$-$\infty$-free, respectively).

**Example 7.9.** Any complete theory of $R$-modules (for $R$ a commutative, unital ring) is stable 1-based (folklore, see for example [28, Example 3.17]), and has $\mathsf{nfcp}$ ([3, Remark in Section 4.6]). It thus satisfies the assumptions of Theorem 7.2.

7.2. ***Ab initio* Hrushovski constructions.** We now illustrate an application of Theorem 7.2, showing that the pregeometric or 1-basedness assumptions are not necessary in order to obtain the Zarankiewicz statements (and only the lack of a definable infinite field may be). The definition of the *ab initio Hrushovski construction* structure $\mathcal{M}_0$ can be found in [26, Section 1] and has its origins in [31]. It is well-known that $\mathcal{M}_0$ is stable but not 1-based. What is important here is that by [26, Theorem 3.3], $\mathcal{M}_0$ is a reduct of a stable 1-based (and also trivial, see [46, Definition IV.2.1]) structure. It is also well-known that $\mathsf{acl}$ in $\mathcal{M}_0$ is given by the so-called *self-sufficient closure* (that is, by the pre-dimension function) and does not satisfy the exchange property. See [26] for more details on the above facts.

**Corollary 7.10.** *Let $\{E_b\}_{b \in I}$ be a type-definable family of relations $E_b \subseteq \prod_{i \in [r]} M_0^{d_i}$. Then $Zar(\mathcal{E})$ (and $Zar^{\infty}(\mathcal{E})$, if $\mathcal{M}_0$ is saturated) holds.*

*Proof.* By [26, Section 3.2], $\mathcal{M}_0$ has $\mathsf{nfcp}$. By [26, Theorem 3.3], $\mathcal{M}_0$ is a reduct of a stable 1-based structure. Since $Zar(E_b)$ is preserved under taking reducts, the result follows from Theorem 7.2. $\square$

7.3. **Regular types.** Let $\mathcal{M}$ be a model of a stable theory. Let $p$ be a non-algebraic type, and $U = p(M)$ its set of realisations in $\mathcal{M}$. We call $U$ the *locus* of $p$. Let $\mathsf{cl}_p : \mathcal{P}(U) \to \mathcal{P}(U)$ be the operator defined by:

$$\mathsf{cl}_p(B) = \{b \in U : b \not\perp B\},$$

for all $B \subseteq U$. It is clear that $\mathsf{cl}_p$ satisfies Definition 6.8(1)–(3).

**Definition 7.11.** Let $p$ be a stationary non-algebraic type. Then $p$ is *regular* if $\mathsf{cl}_p$ is a pregeometry.

*Remark* 7.12. The definition of a regular type given above is equivalent to the usual one encountered in geometric stability theory (by [46, Lemma I.4.5.1] and forking symmetry in stable theories). See also [46, Remark VII.1.1]. Regular types generalise types in strongly minimal theories (and, more generally, rank-1 types in stable theories).

**Definition 7.13.** Let $p$ be a regular type. We say that $p$ is *locally modular* (respectively, *trivial*) if $(X, \mathsf{cl}_p)$ is locally modular (respectively, trivial) (Definition 6.11).[4]

**Fact 7.14** ([46, Corollary VII.2.5])**.** *Let $T$ be a stable* 1*-based theory. Then every regular type is locally modular.*

The converse is not true. Indeed, there are non-1-based stable theories, whose all regular types are locally modular (for example, the 'free pseudoplane', [46, Example VII.2.10 and Proposition IV.1.7]). Although Theorem 7.2 does not apply to them, Theorem 7.16 below does apply to the locus of any regular type therein.

Another interesting example of a locally modular regular type in a stable non-1-based theory arises from the *heat variety*, in the theory of $\mathsf{DCF}_{0,2}$ of differentially closed fields of characteristic 0 with two commuting derivations.

**Example 7.15** ([50, Theorem 7.4(2)])**.** Let $\mathcal{M}$ be a saturated model of $\mathsf{DCF}_{0,2}$ and $K \preccurlyeq \mathcal{M}$ a small elementary submodel of $\mathcal{M}$. Let $G$ denote *heat variety*, that is, the group defined by $\delta_1 y = \delta_2^2 y$, and let $p$ be the generic type of $G$ over $K$. Then $p$ is locally modular.

The main result of this subsection is the following application of Theorem 6.17, in the locus $U$ of a regular type $p$. Recall, $\mathcal{C}_{R,p}$ is the class of all Morley grids in $p$.

**Theorem 7.16.** *Let $T$ be a stable theory, $\mathcal{M} \models T$, $p$ a locally modular regular type, $U = p(\mathcal{M})$ its locus in $\mathcal{M}$, and $\mathsf{cl}_p$ the forking closure operator on $U$. Let $\mathcal{U}$ be the induced structure on $U$ by $\mathcal{M}$, and $\mathcal{C} = \mathcal{C}_{R,p}$. Let $\mathcal{E} = \{E_b\}_{b \in I}$ be a type-definable (in $\mathcal{U}$) family of relations $E_b \subseteq \prod_{i\in[r]} U^{d_i}$. Then $Zar(\mathcal{E}, \mathcal{C})$ (and $Zar^{\infty}(\mathcal{E}, \mathcal{C})$, if $\mathcal{M}$ is saturated) holds.*

*Proof.* First, note that $\mathop{\smile\!\!\!\!\!\!|}^{\mathsf{cl}_p} = \mathop{\smile\!\!\!\!\!\!|}$ on $\mathcal{P}(U)^3$. Moreover, it satisfies the independence axioms from Definition 6.12. Indeed, it is well-known that the properties in (1) hold for forking independence in stable theories, and (2)–(3) are true of any independence relation obtained from a pregeometry (we do not know this for (1)).

As usual, we assume that $\mathcal{M}$ is saturated, and the non-saturated case will follow from Remark 6.19(3). It follows then that $\mathcal{U}$ is also saturated. We need only verify the assumptions of Theorem 6.17 (for $\mathcal{U}$, $\mathsf{cl}_p$ and $\mathcal{C}_{R,p}$).

(DEF) Since in stable theories forking and dividing coincide, forking is witnessed by a formula, so an easy argument shows (DEF).

(UB) By Example 6.2.

(WLM) By local modularity of $(M, \mathsf{cl}_p)$. □

*Remark* 7.17. The Zarankiewicz statements for arbitrary relations do not always hold in models of $\mathsf{DCF}_{0,n}$, since the latter define infinite fields and thus Szemerédi-Trotter phenomena occur. Our machinery, thus, allows us to uncover *local* versions of Zarankiewicz's statements for restricted classes of grids, concentrating entirely on the locus of a regular type, such as the generic type of the heat variety.

[4]By [46, Remark IV.2.2] this is consistent with the notion of triviality used in [26].

## 8. Geometric settings revisited

In this section, we use Theorem 6.17 to obtain Theorem B′ (Theorem 8.4), as well as strengthen Theorem A′ in the saturated setting (Theorem 8.9). We note that we cannot reduce Presburger Zarankiewicz to the saturated setting (as was the case for example with the linear Zarankiewicz [4, Theorem 5.6]), since the class $\mathcal{C}_{\mathbb{Z}}$ of $\mathbb{Z}$-distant grids (Definition 1.7) we consider is smaller than that of all grids, resulting to $\mathcal{C}_{\mathbb{Z}}$-$k$-freeness (for a fixed $k \in \mathbb{N}$) not being a first-order property.

### 8.1. Parametric Zarankiewicz for models of Pres.

Here we establish Theorem B′. As mentioned in Example 6.2, models of $\mathsf{Pres}$ do not satisfy (ub) for the class of all grids. We remedy this problem by showing that saturated models satisfy it for the class of all $\mathbb{Z}$-distant grids (Definition 1.7).

*In the rest of this subsection, $\mathcal{M} \models \mathsf{Pres}$ is saturated, $\mathsf{cl} = \mathsf{acl}$ and $\mathcal{C}_{\mathbb{Z}}$ is the class of all $\mathbb{Z}$-distant grids.*

We verify properties (ub) an (wlm).

**Uniform bounds** (ub). It suffices to prove Proposition 8.2 below. For the notions of a copy of $\mathbb{Z}^n$ and $\mathbb{Z}$-volume, see Definition 4.26. Note that if $\alpha : X \subseteq M^n \to M$ is a linear map and $x, y \in X$ with $x - y \in \mathbb{Z}^n$, then $\alpha(x) - \alpha(y) \in \mathbb{Z}$. As usual, $\pi : M^l \to M^{l-1}$ denotes the projection onto the first $l-1$ coordinates.

**Lemma 8.1.** *Let $X \subseteq M^{m+l} = M^n$ be a Presburger cell, and $b \in M^m$. Suppose that $X_b$ is finite and $\{b\} \times \pi(X_b)$ has $\mathbb{Z}$-volume. Then $X_b$ has $\mathbb{Z}$-volume.*

*Proof.* For $l = 0$, it is obvious, so let $l > 0$. Fix $w \in \pi(X_b)$. We split into two cases.

**Case I.** $X = \Gamma(\alpha)$. Pick any $x \in \pi(X_b)$. Since $x - w \in \mathbb{Z}^{l-1}$, we have $(b, x) - (b, w) \in \mathbb{Z}^{m+l-1}$, and hence $\alpha(b, x) - \alpha(b, w) \in \mathbb{Z}$. Therefore

$$\{b\} \times X_b \subseteq (b, w, \alpha(b, w)) + \mathbb{Z}^n,$$

as needed.

**Case II.** $X$ is a cylinder. Since $X_b$ is finite, $X = [\alpha, \beta]_r^c$ for some linear maps $\alpha, \beta$. Denote $\alpha_b = \alpha(b, -)$ and $\beta_b = \beta(b, -)$. Then $\alpha_b(w) - \beta_b(w) \in \mathbb{Z}$. Pick any $x \in \pi(X)$. Then

$$\beta_b(x) - \alpha_b(x) = \beta_b(x) - \beta_b(w) + \beta_b(w) - \alpha_b(w) + \alpha_b(w) - \alpha_b(x) \in \mathbb{Z}.$$

Now take any $z \in (X_b)_x$. We have

$$z - \alpha_b(w) = z - \alpha_b(x) + \alpha_b(x) - \alpha_b(w) \leq \beta_b(x) - \alpha_b(x) + \alpha_b(x) - \alpha_b(w) \in \mathbb{Z}$$

and hence

$$\{b\} \times X_b \subseteq (b, w, \alpha(b, w)) + \mathbb{Z}^n,$$

as needed. □

**Proposition 8.2.** *Let $\{X_b\}_{b \in I}$ be a definable family of subsets of $M^l$. Then there is $N \in \mathbb{N}$ such that for every $b \in I$ with $X_b$ finite, and every $\mathbb{Z}$-distant set $Y \subseteq M^l$, we have $|X_b \cap Y| \leq N$.*

*Proof.* By Presburger cell decomposition, $X = \bigcup_{b\in I}\{b\} \times X_b$ is a finite union of Presburger cells, each of the form $C = \bigcup_{b\in J}\{b\} \times C_b$, for some $J \subseteq I$. It is easy to see that if the claim holds for each such $C$, then it holds for $X$. Thus, we may assume $X$ is a Presburger cell. The proposition will follow from the next claim.

**Claim.** *For every $b \in I$, if $X_b$ is finite, then $X_b$ has $\mathbb{Z}$-volume.*

*Proof of Claim.* We work by induction on $l$. For $l = 0$, the claim holds trivially. Let $l > 0$, and assume that the claim holds for $X' = \bigcup_{b\in I}\{b\} \times \pi(X_b)$. Therefore, if $X_b$ is finite (and hence $\pi(X_b)$ also is), then $\pi(X_b)$ has $\mathbb{Z}$-volume. By Lemma 8.1, we have that $X_b \cap (D \times M)$ has $\mathbb{Z}$-volume, as needed. □

Since a $\mathbb{Z}$-distant set can only contain one element from each copy of $\mathbb{Z}^l$, we conclude the proposition with $N$ the number of Presburger cells in the cell decomposition of the first paragraph. □

**Weak local modularity** (WLM). We adapt the proof of [7, Proposition 6.9], where linear o-minimal structures are shown to be weakly locally modular. Our argument is essentially the one given there, with a few differences that we point out. First, we need the following easy fact about definable functions in Presburger arithmetic.

**Fact 8.3.** *Let $a, b \in M$ and suppose that $b \in \mathsf{dcl}(a) \setminus \mathsf{dcl}(\emptyset)$. Then, there is a $\emptyset$-definable linear function $f$ defined on a* Pres*-interval $C$ containing $a$, such that $f(a) = b$, and $f$ is strictly monotone on $C$.*

Next let us isolate the following property from [7], for a theory $T$:

(‡) For any singletons $a, b$ and tuple $c$ from a saturated model $\mathcal{N} \models T$, if $a \in \mathsf{acl}(b, c)$, then there exists a tuple $u$ from $M$ and a singleton $d \in \mathsf{acl}(c, u)$, such that $u \downarrow abc$ and $a \in \mathsf{acl}(b, d, u)$.

In [7, Theorem 4.3], the authors show that if $T$ is a pregeometric theory and eliminates $\exists^\infty$, then (‡) is equivalent to weak local modularity. It is an easy observation that the proof of [7, Theorem 4.3] does not, in fact, require the assumption that $T$ eliminates $\exists^\infty$, and it suffices to assume that $T$ is pregeometric.[5] Thus for Presburger arithmetic, (‡) suffices to prove (WLM).

In the proof of [7, Proposition 6.9], it is shown that (‡) holds for the theory $T$ of linear o-minimal structures. The proof of (‡) for $T = \mathsf{Pres}$ is word-for-word the same with that one, after replacing the use of 'continuity' in the o-minimal context, with 'linearity' in models $\mathcal{M}$ of Pres, and the group interval $(b_1, b_2)$ with $M$, and using Fact 8.3 to obtain the interdefinability mentioned in the second paragraph of the proof. The use of [7, Proposition 6.8] is replaced by the fact that $M$ is dcl-independent from $P(\mathcal{M})$ over $P(\mathcal{M})$, for any $P \subseteq M$ (which is trivial).

**Theorem 8.4.** *Let $\mathcal{M} \models \mathsf{Pres}$, $\mathcal{E} = \{E_b\}_{b\in I}$ a type-definable family of relations $E_b \subseteq \prod_{i\in[r]} M^{d_i}$, and $\mathcal{C}_\mathbb{Z}$ the class of all $\mathbb{Z}$-distant grids. Then $Zar(\mathcal{E}, \mathcal{C}_\mathbb{Z})$ (and $Zar^\infty(\mathcal{E}, \mathcal{C}_\mathbb{Z})$, if $\mathcal{M}$ is saturated) holds.*

*Proof.* Assume first that $\mathcal{M}$ is saturated. We have proved (UB) an (WLM), whereas (DEF) is by Example 6.1. We can thus apply Theorem 6.17. The non-saturated case follows from Remark 8.10. □

[5] We thank A. Berenstein and E. Vassiliev for confirming this.

8.2. **Semibounded Zarankiewicz.** Here we use Theorem 6.17 to point out an alternative proof to (a stronger version of) Theorem A′ in the saturated setting.

*For the rest of this subsection, $\mathcal{M} = \langle M, <, +, \dots \rangle$ denotes a saturated semibounded o-minimal expansion of an ordered group, and $R \subseteq M$ a fixed short interval; that is $R$ is bounded and a field with domain $R$ is definable in $\mathcal{M}$ (see Section 3).*

To apply Theorem 6.17, we need an operator that satisfies (WLM), (DEF) and (UB) for a suitable class of grids (Definition 8.7). As mentioned in Example 6.14, (WLM) for $\mathsf{acl}$ characterises linearity among o-minimal structures, and hence fails in our $\mathcal{M}$. However, we show that (WLM) holds for the *short closure operator* $\mathsf{scl}$ from [21], defined as follows. For $A \subseteq M$, let

$$\mathsf{scl}(A) = \{a \in M : \text{there is an } A\text{-definable short interval that contains } a\}.$$

By [1, Fact 4.4],

$$\mathsf{scl}(A) = \mathsf{dcl}(AR).$$

**Fact 8.5.** *If $B$ is an $A$-definable short set, then $B \subseteq \mathsf{scl}(A)$.*

*Proof.* Clearly, the projection of $B$ onto each coordinate is short, so a finite union of $A$-definable short intervals and points. So if $b \in B$, then every coordinate of $b$ is contained in an $A$-definable short interval or is in $\mathsf{dcl}(A)$. Hence $b \in \mathsf{scl}(A)$. ☐

Let $\Lambda$ be as in Section 3.1. Denote by $\mathcal{M}_{lin} = \langle M, <, +, 0, \{\lambda\}_{\lambda \in \Lambda} \rangle$ the *linear reduct* of $\mathcal{M}$, and by $\mathsf{dcl}_{lin}$ the definable closure in $\mathcal{M}_{lin}$. Since $\mathsf{dcl}_{lin}$ is weakly locally modular ([7, Proposition 6.9]), and weak local modularity is preserved under taking localisations, our result will follow from the following proposition. This strategy is inspired by ([7, Theorem 4.3]), which shows that in dense pairs $\langle \mathcal{N}, P \rangle$ of linear o-minimal structures, the localisations to $P$ of the full $\mathsf{acl}$ and of the $\mathsf{acl}$ in the base structure $\mathcal{N}$ coincide.

**Lemma 8.6.** *We have*

$$\mathsf{dcl}(-R) = \mathsf{dcl}_{lin}(-\mathsf{dcl}(R)).$$

*Hence,* $\mathsf{scl}$ *is weakly locally modular.*

*Proof.* ($\supseteq$). $\mathsf{dcl}_{lin}(A\,\mathsf{dcl}(R)) \subseteq \mathsf{dcl}(A\,\mathsf{dcl}(R)) \subseteq \mathsf{dcl}(AR)$.

($\subseteq$). If $d \in \mathsf{dcl}(AR)$, then there are a tuple $c \subseteq A$ and an $R$-definable function such that $d = f(c)$. By the full version of the Structure Theorem ([21, Theorem 3.8]), we may assume that $f : C \to M^n$, where $C = B + \sum_{i=1}^{k} v_i t_i | J_i$ is an $R$-definable cone containing $c$, and such that $f$ is 'almost linear with respect to $C$'. That is,

$$c = b + \sum_{i=1}^{k} v_i t_i \quad \text{and} \quad d = f(b) + \sum_{i=1}^{k} \mu_i t_i,$$

for some $b \in B$, $t_i \in J_i$, and $\mu_i \in \Lambda^n$. By [1, Lemma 4.10], $B$ is $R$-definable, hence $f(b) \in \mathsf{dcl}(R)$. Moreover, by the first displayed formula, and since cones are normalized, every $t_i \in \mathsf{dcl}_{lin}(cb) \subseteq \mathsf{dcl}_{lin}(A\,\mathsf{dcl}(R))$. Hence, by the second displayed formula, $d \in \mathsf{dcl}_{lin}(A\,\mathsf{dcl}(R))$, as needed. ☐

**Definition 8.7.** A grid $B$ is *tall* if for every $i \in [r]$ and $x, y \in B_i$, $x - y$ is tall.

**Lemma 8.8.** (UB) *holds for* $\mathsf{cl} = \mathsf{scl}$ *and* $\mathcal{C}$ *the class of all tall grids.*

*Proof.* Let $\{X_b\}_{b\in I}$ be as in (UB). Note that $X_b \subseteq \mathsf{scl}(Ab)$ implies that $X_b$ is short (by [21, Corollary 5.10] the 'long dimension' of $X_b$ is 0, and by [21, Lemma 3.6], it is short). Now, by o-minimality, there is $N$ such that every $X_b \subseteq M^d$ is a finite union of $N \in \mathbb{N}$ definably connected sets, each of course short, when $X_b$ is. So it suffices to prove that every definable, definably connected short set $X \subseteq M^d$ cannot contain two elements $x, y$ whose difference is tall. By taking the projection on each of the $d$ coordinates, we see that $X$ is contained in the product of short intervals. So for any $a, b \in X$ and $i = 1, \dots, d$, $a_i - b_i$ is short, as needed. □

We are now ready to state the main result of this section. Let us note that as a corollary to Theorem A′, we obtain that for the family $\{E_b\}_{b\in I}$ mentioned there, if $\mathcal{M}$ is saturated, we can choose $m_b$ to be the same element $m$ for each $b \in I$. Indeed, every $m(E_b)$ is short, and hence we can let $m$ be any tall element. The advantage of the next theorem is that we can enlarge the class $\mathcal{C}$ of grids even more to the class of all tall grids (not just those that are $m$-distant for some fixed tall $m$).

**Theorem 8.9.** *Let $\mathcal{M} = \langle M, <, +, \dots\rangle$ be a saturated semibounded o-minimal expansion of an ordered group, and $\mathcal{C}$ the class of all tall grids. Let $\mathcal{E} = \{E_b\}_{b\in I}$ be a definable family of relations $E_b \subseteq \prod_{i\in[r]} M^{d_i}$. Then $Zar^\infty(\mathcal{E}, \mathcal{C})$.*

*Proof.* For $\mathcal{C}$ and $\mathsf{cl} = \mathsf{scl}$, we have proved (UB) and (WLM) (Lemmas 8.8 and 8.6), whereas (DEF) holds by [21, Remark 5.3]. Apply Theorem 6.17. □

We follow-up Remark 6.19(1)&(2), extending them to more classes $\mathcal{C}$.

*Remark* 8.10. If $\mathcal{M}$ is a saturated model of $\mathsf{Pres}$, or a semibounded o-minimal structure, and $\mathcal{C}$ the class of all $\mathbb{Z}$-distant grids, or $m$-distant grids, $m \in \mathbb{N}$, respectively, then Remark 6.19(1)&(2) still hold after replacing $k$-free, $\infty$-free, $Zar_\alpha(E)$, and $Zar^\infty(\mathcal{E})$ by $\mathcal{C}$-$k$-free, $\mathcal{C}$-$\infty$-free, $Zar_\alpha(E, \mathcal{C})$, and $Zar^\infty(\mathcal{E}, \mathcal{C})$, respectively.

## 9. Open questions

In this last section, we summarise the results of the paper and list some open questions. In Table 1 below, we list horizontally the various aspects of Zarankiewicz's problem in our settings: whether it holds for a single set or parametrically, $Zar$ or $Zar^\infty$, for a definable or type-definable relation/family, in the saturated or any model of a given theory. We note that in the saturated setting, the definable $Zar$ (single or parametric) implies the type-definable $Zar^\infty$, by Remark 6.19(1)&(2). Moreover, in any model $\mathcal{M}$, the type-definable $Zar$ (single or parametric) in a saturated elementary extension of $\mathcal{M}$ implies it for $\mathcal{M}$, by Remark 6.19(3).

| setting | $\mathcal{C}$ | single | par | $Zar$ | $Zar^\infty$ | def | tp-def | sat | any | resource |
|---|---|---|---|---|---|---|---|---|---|---|
| **linear** | all | | • | | • | • | | | • | [4, Cor 5.11] |
| | | | • | | • | | • | • | | Rmk 6.19(1)&(2) |
| | | | • | • | | | • | | • | Rmk 6.19(3) |
| | | | • | | • | | • | | • | ? |
| **semibounded** | $m$-distant | | • | | • | • | | | • | Thm A′ |
| | tall | | • | | • | | • | • | | Thm 8.9 |
| | $m$-distant | | • | • | | | • | | • | Rmk 8.10 |
| | $m$-distant | | • | | • | | • | | • | ? |
| **Presburger** | all | • | | | • | • | | | • | fails (Rmk 4.28) |
| | | • | | • | | | • | | • | Thm B, Rmk 8.10 |
| | | | • | • | | • | | • | | fails (Section 1.3) |
| | $\mathbb{Z}$-distant | | • | | • | | • | • | | Thm 8.4 |
| $\langle \mathbb{R}, <, +, \mathbb{Z} \rangle$ | all | • | | • | | • | | — | | Thm C |
| **abstract** | $\supseteq$ cl-indep | | • | | • | | • | • | | Thm D |
| **1-based** | all | | • | • | | | • | | • | Thm E |
| | | | • | | • | | • | • | | Rmk 6.19 |
| **regular** $p$ | Morley in $p$ | | • | • | | | • | | • | Thm F |
| | | | • | | • | | • | • | | Rmk 6.19 |

Table 1. Summary of results

(For example, in any linear o-minimal structure, every definable family $\mathcal{E}$ satisfies $Zar^\infty(\mathcal{E})$ [4])

We now list various questions in different settings, often referring the reader to the corresponding literature for any undefined terminology.

9.1. **Semibounded and Presburger settings.** We ask whether we can strengthen Theorem A″ so that we recover fields on intervals of any given length.

**Question 9.1.** Let $\mathcal{M} = \langle M, <, +, \dots \rangle$ be an o-minimal expansion of an ordered group, and $R = (-\gamma, \gamma)$ an open interval. Are the following equivalent?

(1) There is no definable field with domain $R$.
(2) Let $1 \leq \beta < \frac{4}{3}$ in $\mathbb{R}_{>0}$. For every binary definable $E \subseteq M^{d_1} \times M^{d_2}$, there are $\alpha \in \mathbb{R}_{>0}$ and $m \in R$, such that if $E$ is $k$-free, for some $k \in \mathbb{N}$, then for every $m$-distant $n$-grid $B \subseteq M^{d_1} \times M^{d_2}$, we have

$$|E \cap B| \leq \alpha n^\beta.$$

The idea of restricting our class of grids to sufficiently distant ones in the semibounded Zarankiewicz can be adopted to other statements as well, such as those in [5, Introduction]. We only deal with [5, Theorem 1.4] here. Let $\mathcal{S}$ be a family of subsets of $M^n$, and $m \in M_{\geq 0}$. The *$m$-distant shatter function* $\pi_{\mathcal{S},m} : \mathbb{N} \to \mathbb{N}$ is

$$\pi_{\mathcal{S},t}(t) = \max\{|S \cap A| : A \subseteq M^n \ m\text{-distant}, |A| = t\}.$$

For $m = 0$, we obtain the usual shatter function $\pi_{\mathcal{S}}$.

**Question 9.2.** Let $\mathcal{M} = \langle M, <, +, \dots \rangle$ be an o-minimal expansion of an ordered group. Are the following equivalent?

(1) $\mathcal{M}$ is semibounded.
(2) for every set system $\mathcal{S}$ definable in $\mathcal{M}$, there is $m \in M_{\geq 0}$ such that $\pi_{\mathcal{S},m}(t)$ is asymptotic to a polynomial.

We could also adopt this question to $\mathsf{Pres}$.

**Question 9.3.** Let $\mathcal{M} \models \mathsf{Pres}$, and $\mathcal{S}$ definable. Is $\pi_{\mathcal{S}}$ asymptotic to a polynomial?

Although $Zar(\mathcal{E})$ need not hold for definable families in models of $\mathsf{Pres}$ (Section 1.3), the following parametric version seems plausible.

**Question 9.4.** Let $\mathcal{M} \models \mathsf{Pres}$, and $\mathcal{E} = \{E_b\}_{b\in I}$ a definable family of relations $E_b \subseteq \prod_{i\in[r]} M^{d_i}$. Fix $k \in \mathbb{N}$. Then is there $\alpha \in \mathbb{R}_{>0}$ (depending on $k$) such that for every $b \in E$, if $E$ is $k$-free, then $E$ has linear $Z$-bounds for the class of all grids witnessed by $\alpha$?

9.2. **Zarankiewicz's problem in more settings.** There is a variety of settings where the Zarankiewicz statements could be sought. We list here some.

**Question 9.5.** Do we obtain a Zarankiewicz statement for sets definable in (a) valued vector spaces (asked by A. Chernikov), (b) ordered abelian groups (M. Vicaria), (c) $\langle \mathbb{R}, <, +, P\rangle$ where $P$ is some multiplicative group, such as $P = \mathbb{Z}[\frac{1}{2}]$, (d) more real-integer structures $\langle \mathbb{R}, <, +, Z\rangle$ in the sense of [38], where for example $Z$ could be a principal ideal domain, (e) a non-standard model of the theory of $\langle \mathbb{R}, <, +, \mathbb{Z}\rangle$?

**Question 9.6.** What other contexts can the abstract Zarankiewicz (Theorems 6.7 and 6.17) be applied to? In particular, can Theorem 7.16 be proved for a broader class $\mathcal{C}$ than $\mathcal{C}_{R,p}$, yielding even applications relevant to the heat variety?

9.3. **Converse Zarankiewicz's problem.** So far, the only known way to violate the Zarankiewicz statements is the use of the Szemerédi-Trotter theorem ([38, Proposition 4.2.1]). A strong question we could ask is the following (suggested also by M. Aschenbrenner).

**Question 9.7.** Let $\mathcal{M}$ be a structure that does not interpret an infinite field, and $E$ a definable relation. Does $Zar(E)$ hold?

More specifically, we can ask for a converse to Theorem B$'$.

**Question 9.8.** Let $\mathcal{M} = \langle \mathbb{Z}, <, +, \dots\rangle$. Are the following equivalent?

(1) The multiplication $\cdot$ is not definable in $\mathcal{M}$.
(2) For every binary definable $E \subseteq M^{d_1} \times M^{d_2}$, $Zar(E)$ holds.

We can also explore a converse Zarankiewicz's problem in the stable context, where the linearity condition is captured by 1-basedness.

**Question 9.9.** Let $T$ be a stable theory with $\mathsf{nfcp}$, and assume that all types have finite $U$-rank. Are the following equivalent?

(1) $T$ is 1-based.
(2) For every $\mathcal{M} \models T$ and binary definable $E \subseteq M^{d_1} \times M^{d_2}$, $Zar(E)$ holds.

Without the assumption of finite $U$-rank, the question admits a negative answer, as manifested by the ab initio Hrushovski constructions (Section 7.2): (1) is not always preserved under taking reducts, but (2) is. We could also replace (1) by $T$ being not $n$-ample, for some $n \in \mathbb{N}$ (see [47] for the notion of ampleness). Finally, we could replace the assumption of being stable with finite $U$-rank by $T$ being topological, in the sense of the recent [8], where 1-basedness is explored in that setting.

A relevant question is the following.

**Question 9.10.** Let $\mathcal{M}$ be a structure interpretable in another structure $\mathcal{N}$. Suppose that every definable relation $E$ in $\mathcal{N}$ satisfies $Zar(E)$. Does every definable relation $E'$ in $\mathcal{M}$ also satisfy $Zar(E')$?

9.4. **Tame expansions.** In [7], the authors show that in dense pairs of linear o-minimal structures $\langle \mathcal{M}, \mathcal{P} \rangle$, $\mathsf{dcl}$ in the pair coincides with $\mathsf{dcl}$ in $\mathcal{M}$, and hence the parametric Zarankiewicz statement holds in this setting as well, by [4, Theorem 5.6].

**Question 9.11.** Do Zarankiewicz statements hold in dense pairs of semibounded o-minimal structures?

**Question 9.12.** What is a suitable notion of lovely pairs of models of $\mathsf{Pres}$, and do Zarankiewicz statements hold in them?

**Question 9.13.** Do Zarankiewicz statements hold in $\langle \mathbb{N}, <, +, 2^{\mathbb{Z}} \rangle$, $\langle \mathbb{R}, <, +, 2^{\mathbb{Z}} \rangle$?

Hieronymi-Walsberg [30] introduced a tetrachotomy for tame expansions of the real ordered group and it is plausible to ask whether the Zarankiewicz statements can capture combinatorially some of the categories of the structures they consider.

**Question 9.14.** Let $\mathcal{M} = \langle \mathbb{R}, <, +, \dots \rangle$ be an expansion of the real ordered group of type $A$. For what classes of grids $\mathcal{C}$, are the following equivalent?

(1) $\mathcal{M}$ is not of field-type (for example, the expansion of $\langle \mathbb{R}, <, + \rangle$ by all subsets of all $\mathbb{Z}^n$).
(2) For every binary definable $E \subseteq M^{d_1} \times M^{d_2}$, $Zar(E, \mathcal{C})$ holds.

9.5. **Local Zarankiewicz.** This paper deals with the *global* linear Zarankiewicz problem – namely, we assume the $k$-freeness condition on the ambient relation $E$, versus the *local* version, where we would assume it only on the intersections $E \cap B$ whose size we want to bound. The local linear Zarankiewicz problem has been dealt with in [4] and the Presburger one in [11, Remark 2.24], [54, Theorem 6.9], where they obtain 'almost linear' Zarankiewicz bounds.

**Question 9.15.** Do we obtain almost linear Zarankiewicz bounds in the local Zarankiewicz problem in the settings (a) (for sufficiently distant grids) and (c)-(e) from the abstract of this paper, as well in ordered abelian groups?

Observe that our Reduction Strategy (Section 2.7) would not be applicable here since assuming $k$-freeness on each $E \cap B$ might not yield enough information for $E$.

9.6. **Recovering other structure from Zarankiewicz statements.** In private communicatoin, J. Pila asked whether there is a way to tweak the Zarankiewicz statements so that in condition (1) of Theorem A$''$ we recover some structure other than the multiplication, such as exp. While it is still unclear what a suitable statement would be, we risk the following question.

**Question 9.16.** Let $\mathcal{M} = \langle \mathbb{R}, <, +, \dots \rangle$ be an o-minimal expansion of the real ordered group. For what classes $\mathcal{C}$ of grids, are the following equivalent?

(1) The exponential exp is not definable in $\mathcal{M}$.
(2) For every binary definable $E \subseteq M^{d_1} \times M^{d_2}$, $Zar(E, \mathcal{C})$ holds.

We could also change (1) to $\mathcal{M}$ having strictly more structure than the real field.

9.7. **Varying the Zarankiewicz bounds.** As mentioned in the introduction, [12] proves a version of Elekes-Szabó in the o-minimal and other settings, with different Zarankiewicz bounds and recognising algebraic groups instead of fields. We can ask whether we can recover *global* algebraic groups, following the idea of Theorem A$''$ in considering sufficiently distant grids (compare with [12, Corollary 6.21]). (Similar questions can be asked regarding [12, Theorem 6.4, Corollaries 6.19–6.20]).

**Question 9.17** (with J. Dobrowolski)**.** Assume $r \geq 3$. Let $\mathcal{M} = \langle M, <, \ldots \rangle$ be an o-minimal structure, and assume that there is an ordered group $\mathcal{R} = \langle M, <, + \rangle$. (We could further assume that $\mathcal{M}$ is a reduct of an o-minimal expansion of $\mathcal{R}$.) Let $E \subseteq M^r$ be definable, such that every $p_{l \restriction E}$, $l \in [r]$, is finite-to-1. Then at least one of the following holds (with the notion of $m$-distant taken in $\mathcal{R}$).

(1) There are $m \in M_{>0}$ and $\alpha \in \mathbb{R}_{>0}$, such that for any finite $m$-distant grid $A = A_1 \times \cdots \times A_r \subseteq M^r$, with $|A_i| = n$, $i \in [r]$, we have

$$|E \cap (A_1 \times \cdots \times A_r)| \leq \alpha n^{r-1-\gamma},$$

where $\gamma = \frac{1}{3}$ if $r \geq 4$, and $\gamma = \frac{1}{6}$ if $r = 3$.

(2) There exist open sets $U_i \subseteq M$, $i \in [r]$ and homeomorphisms $\pi_i : U_i \to M$ such that for all $x_i \in U_i$, $i \in [r]$,

$$\pi_1(x_1) + \cdots + \pi_r(x_r) = 0 \Leftrightarrow E(x_1, \ldots, x_r).$$

9.8. **Definably complete ordered groups.** In this setting, in Section 4, we introduced the notions of $\mathbf{Z}$-cells, purely $\mathbf{Z}$-unbounded sets, and $N$-internal points, $N \in \mathbf{Z}$, for a fixed copy $\mathbf{Z}$ of the integers in our structure $\mathcal{M}$.

**Question 9.18.** Let $\mathcal{M} = \langle M, <, + \rangle$ be a definably complete ordered group.

(1) Let $E \subseteq M^n$ be a definable set. Does $Zar(E)$ hold (without the assumption in Corollary 4.30 of $\mathbf{Z}$ being cofinal in $M$)? Does Corollary 4.33 hold after replacing 'Presburger cell' by '$\mathbf{Z}$-cell', and '$\mathbb{N}$' by '$\mathbf{N}$'?
(2) Do Propositions A.2 and A.6 hold after replacing 'linear cell' by '$\mathbf{Z}$-cell'? Do they hold after also replacing 'purely unbounded' by 'purely $\mathbf{Z}$-unbounded'?
(3) (asked by J. Losh) Suppose $<$ is dense. Does every $\mathbf{Z}$-cell $C$ containing 0 satisfy the conclusion of Lemma 2.13, and hence $(*)_\infty$ and $Zar(C)$ as in Proposition 2.19?
(4) Under what further conditions on $\mathcal{M}$ does a $\mathbf{Z}$-cell decomposition theorem holds?

9.9. ***N*-internal points in linear o-minimal structures.** In Remark 4.24, we noticed that purely unbounded cells in an ordered vector space $\mathcal{M}$ contain $N$-internal points, for $N \in M$. The following more refined statement remains open.

**Question 9.19.** Let $\mathcal{M}$ be an ordered vector space over an ordered division ring. Let $C$ be a purely $\mathbf{Z}$-unbounded cell. Then for every $N \in \mathbf{N}$, $C$ contains an $N$-internal point.

9.10. **Applications.** Zarankiewicz's problem has been motivated by many applications in discrete geometry and more (see, for example, [27, 19]).

**Question 9.20.** What would be some applications of the local/global versions of Zarankiewicz's problem in any of our settings, and specifically of the semibounded one (perhaps also in view of Remark 3.14).

## Appendix A. Lemmas for semilinear sets and families (with Pablo Andújar Guerrero)

*In this appendix,* $\mathcal{M} = \langle M, <, +, \{x \mapsto \lambda x\}_{\lambda \in \Lambda} \rangle$ *is an ordered vector space over an ordered division ring* $\Lambda$*. 'Semilinear set/linear map/cell' is taken in* $\mathcal{M}$*.* We recall (Section 2.5) that in this setting, a linear map $\alpha : X \subseteq M^n \to M$ extends to a linear map on $M^n$, which we also denote by $\alpha$.

We prove Propositions A.2, A.6 and A.8, which in Section 5 are applied to $\mathcal{M} = \langle \mathbb{R}, <, + \rangle$, as a vector space over $\mathbb{Q}$. It is possible that they can be further extended to the setting of definably complete ordered groups (Section 4, Question 9.18(2)).

A.1. **Decomposition of semilinear sets into product cells.**

**Fact A.1.** *Let $C = (\alpha, \beta)_D \subseteq M^n$ be a linear cell, $n > 0$, $\alpha, \beta : D \to M$. Then $C$ is purely unbounded if and only if $D$ is purely unbounded and $\alpha - \beta$ is not constant.*

*Proof.* By a straightforward induction on $n$, left to the reader. $\Box$

The notion of a 'product cell' was introduced in Definition 4.38.

**Proposition A.2.** *Every linear cell $C \subseteq M^n$ can be partitioned into a finite union of product cells.*

*Proof.* By induction on $n$. For $n = 0$, $C = M^0$ is purely unbounded, and for $n = 1$, $C$ is either purely unbounded or bounded. Let $n > 1$. By Inductive Hypothesis, $\pi(C)$ is a finite union of product cells, and since the restriction of $C$ to each of them is a cell, we may assume that $\pi(C) = J + D$ is a product cell. We prove that

$$C = K + E,$$

where $K$ is a purely unbounded cell and $E$ a bounded cell, and leave the verification of Definition 4.38(2) to the reader. We only consider the cases $C = \Gamma(\alpha), (\alpha, \infty), (\alpha, \beta)$, where $\alpha, \beta : \pi(C) \to M$, the cases $(-\infty, \beta)$ and $(-\infty, \infty)$ being similar. Let $a = \alpha(0)$. We only provide full details in Case III.

**Case I.** $C = \Gamma(\alpha)$. We can take

$$K = \Gamma(\alpha)_J \quad \text{and} \quad E = \Gamma(\alpha - a)_D.$$

**Case II.** $C = (\alpha, \infty)$. We can take

$$K = (\alpha, \infty)_J \quad \text{and} \quad E = \Gamma(\alpha - a)_D.$$

**Convention**: In what follows, if $\gamma : J + D \to M$ is a linear map and $d \in \overline{D}$, then $g \mapsto \gamma(g + d) : J \to M$ is constant if and only if for any $d' \in \overline{D}$, $g \mapsto \gamma(g + d') : J \to M$ is constant. We say in this case that '$\gamma_J$ is constant'. Similarly for $\gamma_D$.

**Case III.** $C = (\alpha, \beta)$. We first handle two subcases:

**Subcase IIIa.** $(\beta - \alpha)_J$ non-constant and $(\beta - \alpha)_D$ constant. Then $C$ consists of all elements $x$ of the form

$$x = (g + d, \alpha(g + d)) + (0, y) = (g, \alpha(g) + y) + (d, \alpha(d) - a),$$

where $g \in J$, $d \in D$ and $y \in (0, \beta(g + d_0) - \alpha(g + d_0))$, for some/any $d_0 \in D$ (since $(\beta - \alpha)_D$ is constant). Hence we can let

$$K = (\alpha, \varepsilon)_J \quad \text{and} \quad E = \Gamma(\alpha - a)_D,$$

where $\varepsilon(g) = \alpha(g) + \beta(g + d_0) - \alpha(g + d_0)$.

**Subcase IIIb.** $(\beta - \alpha)_J$ constant. Then $C$ consists of all elements $x$ of the form

$$x = (g + d, \alpha(g + d)) + (0, y) = (g, \alpha(g) - a) + (d, \alpha(d) + y),$$

where $g \in J$, $d \in D$ and $y \in (0, \beta(g_0 + d) - \alpha(g_0 + d))$, for some/any $g_0 \in J$ (since $(\beta - \alpha)_J$ is constant). Hence we can let

$$K = \Gamma(\alpha - a)_J \quad \text{and} \quad E = (\alpha, \varepsilon)_D,$$

where $\varepsilon(d) = \alpha(d) + \beta(g_0 + d) - \alpha(g_0 + d)$.

In general, we partition $C$ into cells that fall in some of the previous cases.

**Claim.** *There is a linear map $\gamma : \pi(C) \to M$, such that $C = C_1 \cup C_2 \cup C_3$ is the union of three linear cells (in fact, disjoint), where*

*(1) $C_1 = (\alpha, \gamma)$, with $(\gamma - \alpha)_J$ non-constant and $(\gamma - \alpha)_D$ constant,*
*(2) $C_2 = \Gamma(\gamma)$,*
*(3) $C_3 = (\gamma, \beta)$, with $(\beta - \gamma)_J$ constant.*

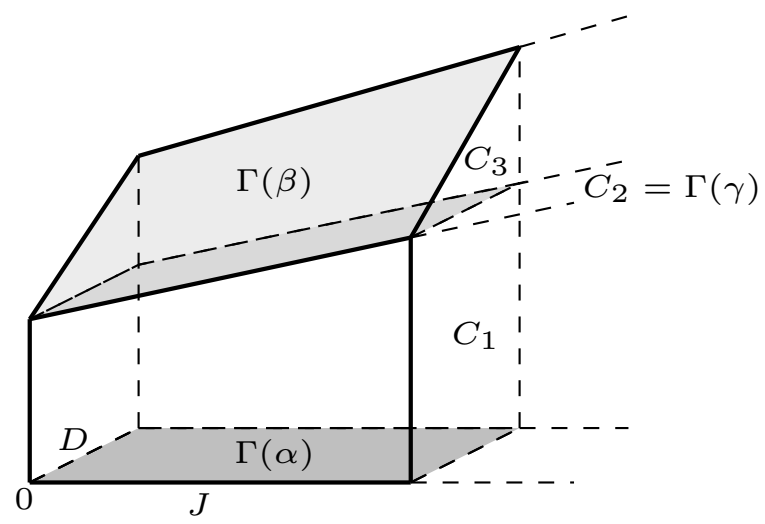

*Proof of Claim.* We may assume that both $(\beta - \alpha)_J$ and $(\beta - \alpha)_D$ are non-constant, otherwise we are directly in cases (1) or (3). Let $k = \inf(\beta - \alpha)_{\pi(C)}$, and suppose that $\beta(g_0 + d_0) - \alpha(g_0 + d_0) = k$, for some $g_0 + d_0 \in \overline{\pi(C)}$, with $g_0 \in \overline{J}$ and $d_0 \in \overline{D}$ (note that $\overline{J} + \overline{D} = \overline{J + D}$).

**Subclaim.** *For every $g \in J$ with $g + d_0 \in \pi(C)$, and $d \in D$,*

$$0 < \beta(g + d_0) - \alpha(g + d_0) < \beta(g + d) - \alpha(g + d).$$

*Proof of Subclaim.* The first $<$ is because $(\beta - \alpha)_J$ is non-constant and hence

$$\beta(g + d_0) - \alpha(g + d_0) > \beta(g_0 + d_0) - \alpha(g_0 + d_0) \geq 0.$$

The second $<$ is by linearity, since

$$\beta(g + d) - \alpha(g + d) - (\beta(g + d_0) - \alpha(g + d_0)) =$$

$$\beta(g_0 + d) - \alpha(g_0 + d) - (\beta(g_0 + d_0) - \alpha(g_0 + d_0)) > 0,$$

where the last inequality is because $(\beta - \alpha)_D$ is non-constant. □

Define now $\gamma : \pi(C) \to M$ via

$$\gamma(g + d) = \beta(g + d_0) - \alpha(g + d_0) + \alpha(g + d).$$

It follows from the subclaim that $\alpha < \gamma < \beta$. Hence for $C_1 = (\alpha, \gamma)$, $C_2 = \Gamma(\gamma)$ and $C_3 = (\gamma, \beta)$, we obtain $C = C_1 \cup C_2 \cup C_3$. We now check the extra properties for $C_1$ and $C_3$.

- $(\gamma - \alpha)_J$ **non-constant:** for every $d \in D$, we have

$$\gamma(g + d) - \alpha(g + d) = \beta(g + d_0) - \alpha(g + d_0),$$

which is non-constant on $J$.

- $(\gamma - \alpha)_D$ **constant:** for every $g \in J$, we have
$$\gamma(g+d) - \alpha(g+d) = \beta(g+d_0) - \alpha(g+d_0),$$
which is fixed.
- $(\beta - \gamma)_J$ **constant:** for every $d \in D$, we have
$$\beta(g+d) - \gamma(g+d) = \beta(g+d) - \beta(g+d_0) + \alpha(g+d_0) - \alpha(g+d),$$
but both $\beta(g+d) - \beta(g+d_0)$ and $\alpha(g+d_0) - \alpha(g+d)$ are constant as $g$ varies, by linearity.

This ends the proof of the claim. □

Let $C_1, C_2, C_3$ be as in the claim. We are then done after observing that each of them is a product cell, by Cases IIIa, I, and IIIb, respectively. □

*Remark* A.3. A bounded set $D \subseteq M^n$ is a product cone $\{0\} + D$, since singletons are purely unbounded (Remark 4.3).

A.2. **Nesting lines for semilinear families.** In Definitions 4.39, 4.40 we defined the notions of a line, nesting line and nesting direction. A line in the current setting has the form $L = p + dM_{\geq 0}$, for some $p \in M^n$ and $d \in \Lambda^n$.

**Fact A.4.** *Let $\alpha : M^n \to M$ be a linear map, and $S_1, S_2 \subseteq M^n$ two definable sets. Suppose that $\mathsf{dim}(S_1 + S_2) = n$. If both $\alpha_{\restriction S_1}$ and $\alpha_{\restriction S_2}$ are constant, then so is $\alpha$.*

*Proof.* Suppose $\alpha_{\restriction S_1}$ and $\alpha_{\restriction S_2}$ are constant. Using linearity, it is easy to see that $\alpha_{\restriction S_1+S_2}$ is also constant. Hence, we need to prove that if $\mathsf{dim} S = n$ and $\alpha_{\restriction S}$ is constant, then so is $\alpha$. Take an open box contained in $S$. Then for every $x$ in it and $\varepsilon$ sufficiently small, we have $\alpha(x+\varepsilon) - \alpha(x) = 0$. By linearity, for every $y \in M^n$ and $\varepsilon \in M^n$ sufficiently small, $\alpha(y+\varepsilon) - \alpha(y) = 0$. So $\alpha$ is locally constant, and since it is continuous, it is constant (by a standard connectedness argument). □

We first prove a claim for semilinear families of unary sets.

**Claim A.5.** *Let $\{X_t\}_{t\in M^m}$ be a semilinear family of sets in $M$, $m > 0$, such that $C = \bigcup_{t\in T}\{t\} \times X_t$ is a linear cell, and $\bigcap_{t\in T} X_t \neq \emptyset$. Suppose $L$ is a line witnessing $T$ is purely unbounded. Suppose also that $C$ is a graph $\Gamma(\alpha)$ or a cylinder $(\alpha, \beta)$, $(\alpha, \infty)$, $(-\infty, \beta)$, or $(-\infty, \infty)$. Then, for $\alpha$ and $\beta$ appearing in the definition of $C$, we have:*

*(1) If $\alpha_{\restriction L}$ (respectively, $\beta_{\restriction L}$) is constant, then so is $\alpha$ (respectively, $\beta$).*
*(2) If $\alpha_{\restriction L}$ (respectively, $\beta_{\restriction L}$) is not constant, then it is strictly decreasing (respectively, strictly increasing).*

*Proof.* If $C = \Gamma(\alpha)$, since $\bigcap_{t\in T} X_t \neq \emptyset$, $\alpha$ has to be constant, and the result follows. So suppose that $C$ is a cylinder. The proof is similar in the various cases, and we handle here only that of $C = (\alpha, \infty)$, by induction on $m$.

For $m = 1$, $T$ is either a point, in which case the claim is obvious, or $\mathsf{dim} L = 1$, and the statement follows from Fact A.4 with $S_1 = L$ and $S_2 = \{0\}$. Now let $m > 1$. We may assume $T$ is a cylinder. Indeed, if $T$ is a graph, the projection $f : T \to M^{m-1}$ onto the first $m-1$ coordinates is injective. For $t \in f(T)$, let $Y_t = X_{f^{-1}(t)}$. Then clearly $\bigcap_{t\in T} Y_t \neq \emptyset$, $f(L)$ is a line witnessing $f(T)$ is purely unbounded, and $f_{\restriction L}$ is order-preserving. It is easy to check that $C' = \bigcup_{t\in f(T)}\{t\} \times Y_t$ is also a linear cell $\Gamma(\alpha')$, $(\alpha', \beta')$, $(\alpha', \infty)$, $(-\infty, \beta')$, or $(-\infty, \infty)$, where $\alpha', \beta' : f(T) \to M$ with $\alpha'(t) = \alpha(f^{-1}(t))$, $\beta'(t) = \beta(f^{-1}(t))$. Hence the conclusion for $C$ follows from that for $C'$, by Inductive Hypothesis.

(1). Assume that $\alpha_{\restriction L}$ is constant. Let $T = (f,g)_{\pi(T)}$, where $f, g$ can also be $-\infty, \infty$, respectively. We can find a linear map $\delta : \pi(L) \to M$ whose graph equals $L$. Indeed, if $L = p + dM_{\geq 0}$, for some $p \in M^n$, $d \in \Lambda^n$, let $\delta : \pi(L) \to M$, with

$$\delta(p_1 + xd_1, \dots, p_{m-1} + xd_{m-1}) = p_m + xd_m.$$

Since $\alpha_{\restriction L}$ is constant, so is the map $s \mapsto \alpha(s, \delta(s)) : \pi(L) \to M$.

We first show that $\alpha_{\restriction (f,g)_{\pi(L)}}$ is constant. If not, then by linearity of $\alpha$, for every $s \in \pi(L)$, $\alpha(s, -)_{\restriction T_s}$ is strictly monotone and of the same behavior. Since $L$ witnesses $T$ is purely unbounded, both $\delta(s) - f(s)$ and $g(s) - \delta(s)$ can become arbitrarily large (or are infinite), by varying $s \in \pi(L)$. If each $\alpha(s, -)_{\restriction T_s}$, $s \in \pi(L)$, is strictly increasing (respectively, decreasing), then by linearity $\alpha(s, x - \delta(s))$ (respectively, $\alpha(s, \delta(s) - x)$) can become arbitrarily large, as $x \in T_s$ approaches $g(s)$ or $f(s)$, respectively. Since $\alpha(s, \delta(s))$ is constant, this yields arbitrarily large values $\alpha(s, x)$ (or $-\alpha(s, x)$), with $x \in T_s$, contradicting $\bigcap_{t \in T} X_t \neq \emptyset$.

We now define a linear map $h : \pi(T) \to M$, as follows.

$$h = \begin{cases} \frac{f+g}{2} & \text{if } f \neq -\infty,\ g \neq \infty \\ f + 1 & \text{If } f \neq -\infty,\ g = \infty \\ g - 1 & \text{if } f = -\infty,\ g \neq \infty \\ 0 & \text{if } f = -\infty,\ g = \infty. \end{cases}$$

Clearly, $f < h < g$. Define also $\alpha' : \pi(T) \to M$,

$$\alpha'(x) = \alpha(x, h(x)).$$

Since $\alpha_{\restriction (f,g)_{\pi(L)}}$ is constant, so is $\alpha'_{\restriction \pi(L)}$. Finally, for $s \in \pi(T)$, let $Y_s = X_{(s, h(s))}$. Then $\bigcup_{s \in \pi(T)} \{s\} \times Y_s = (\alpha', \infty)$ is a linear cell, $\bigcap_{s \in \pi(T)} Y_s \neq \emptyset$, and $\pi(L)$ witnesses that $\pi(T)$ is purely unbounded. By Inductive Hypothesis, we obtain that $\alpha'_{\restriction \pi(T)}$ is constant. That is, $\alpha_{\restriction \Gamma(h)}$ is constant. By Fact A.4, for $S_1 = (f,g)_{\pi(L)}$ (or just any $T_s$, $s \in \pi(L)$) and $S_2 = \Gamma(h)$, we obtain that $\alpha$ is constant.

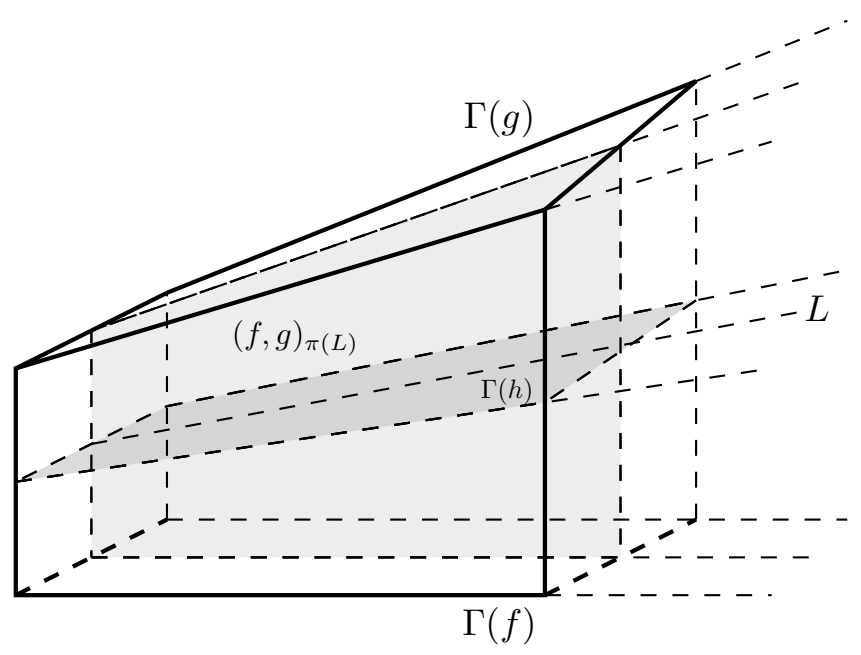

(2). Assume that $\alpha_{\restriction L}$ is not constant. If it is strictly increasing, we contradict $\bigcap_{t \in T} X_t \neq \emptyset$. Therefore, it is strictly decreasing. □

**Proposition A.6.** *Let $\{X_t\}_{t \in M^m}$ be a semilinear family of sets in $M^n$, $m > 0$, such that $C = \bigcup_{t \in T} \{t\} \times X_t$ is a linear cell, and $\bigcap_{t \in T} X_t \neq \emptyset$. Suppose $L$ is a line witnessing $T$ is purely unbounded. Then $L$ is a nesting line for $\{X_t\}_{t \in T}$.*

*Proof.* Let $C = \Gamma(\alpha)$, $(\alpha, \beta)$, $(\alpha, \infty)$, $(-\infty, \beta)$, or $(-\infty, \infty)$. We perform induction on $n$. For $n = 0$, all of $\alpha_{\restriction L}, \beta_{\restriction L}, \alpha, \beta$ are constant equal to 0, so we are done. Let $n > 0$. We establish the two properties from Definition 4.40:

(1) For all $t, t' \in L$, $t < t' \Rightarrow X_t \subseteq X_{t'}$.
(2) For all $t \in T$, there is $t' \in L$, $X_t \subseteq X_{t'}$.

Recall that $\pi : M^n \to M^{n-1}$ denotes the projection onto the first $n-1$ coordinates. Fix some $x = (x_1, x_2) \in \bigcap_{t \in T} X_t$, with $x_1 = \pi(x)$. Then $x_1 \in \bigcap_{t \in T} \pi(X_t)$. Therefore $\{\pi(X_t)\}_{t \in T}$ satisfies the conditions of the proposition . By Inductive Hypothesis, $L$ is a nesting line for $\{\pi(X_t)\}_{t \in T}$. Therefore:

(a) for every $t < t'$ in $L$, $\pi(X_t) \subseteq \pi(X_{t'})$, and
(b) for every $t \in T$, there is $t' \in L$, such that $\pi(X_t) \subseteq \pi(X_{t'})$.

We continue the proof by considering two cases, in which we prove (1) and (2) as follows (onwards, for any $t \in T$ and $s \in M^{n-1}$, $X_{t,s}$ denotes the fiber $(X_t)_s = X_{(s,t)}$):

- to see (1), we let $t < t'$ in $L$, and prove that for every $s \in \pi(X_t) \subseteq \pi(X_{t'})$ (which holds by (a)), $X_{t,s} \subseteq X_{t',s}$.
- to see (2), we let $t \in T$, and find $t' \in L$, such that $\pi(X_t) \subseteq \pi(X_{t'})$ and for every $s \in \pi(X_t)$, $X_{t,s} \subseteq X_{t',s}$.

These are sufficient, since in both cases we obtain $X_t \subseteq X_{t'}$.

**Case (I):** $C = \Gamma(\alpha)$. Since $x_2 \in \bigcap_{t \in T} X_{t,x_1}$ and each $X_{t,x_1} = \{(\alpha(t, x_1)\}$, we have that $\alpha(-, x_1)_{\restriction T}$ is constant (with image $\{x_2\}$). Since $\alpha$ is a linear map, it follows that for every $s \in M^{n-1}$, $\alpha(-, s)_{\restriction T}$ is constant. Hence for every $s \in M^{n-1}$, and $t, t' \in T$ with $s \in \pi(X_t) \cap \pi(X_{t'})$, $\alpha(t, s) = \alpha(t', s)$.

To see (1): let $t < t'$ in $L$. For every $s \in \pi(X_t) \subseteq \pi(X_{t'})$, we have

$$X_{t,s} = \{(s, \alpha(t, s))\} = \{(s, \alpha(t', s))\} = X_{t',s},$$

as needed.

To see (2): let $t \in T$, and take $t' \in L$ as in (b). Then for each $s \in \pi(X_t) \subseteq \pi(X_{t'})$, we again have

$$X_{t,s} = \{(s, \alpha(t, s))\} = \{(s, \alpha(t', s))\} = X_{t',s},$$

as needed.

**Case (II):** $C$ is a cylinder $(\alpha, \beta)$, $(\alpha, \infty)$, $(-\infty, \beta)$, or $(-\infty, \infty)$. We only handle the case $C = (\alpha, \beta)$, as the others are similar. Since $x_2 \in \bigcap_{t \in T} X_{t,x_1}$, the family $\{X_{t,x_1}\}_{t \in T}$ satisfies the conditions of Claim A.5. There are three subcases – after stating them, we prove (1) and (2) simultaneously.

**Subcase IIa.** Both $\alpha(-, x_1)_{\restriction L}$ and $\beta(-, x_1)_{\restriction L}$ are constant. By Claim A.5, both $\alpha(-, x_1)_{\restriction T}$ and $\beta(-, x_1)_{\restriction T}$ are constant. By linearity, we obtain that for every $s \in M^{n-1}$, $\alpha(-, s)_{\restriction T}$ and $\beta(-, s)_{\restriction T}$ are constant.

**Subcase IIb.** $\alpha(-, x_1)_{\restriction L}$ is strictly decreasing and $\beta(-, x_1)_{\restriction L}$ is constant (or the former is constant and the latter strictly increasing, which is analogous). By Claim A.5, $\beta(-, x_1)_{\restriction T}$ is also constant. By linearity, we obtain that for every $s \in M^{n-1}$, $\alpha(-, s)_{\restriction L}$ is strictly decreasing and $\beta(-, s)_{\restriction T}$ is constant.

**Subcase IIc.** $\alpha(-, x_1)_{\restriction L}$ is strictly decreasing and $\beta(-, x_1)_{\restriction L}$ is strictly increasing. Then for every $s \in M^{n-1}$, $\alpha(-, s)_{\restriction L}$ is strictly decreasing and $\beta(-, s)_{\restriction L}$ is strictly increasing.

To see (1): let $t < t'$ in $L$. For every $s \in \pi(X_t) \subseteq \pi(X_{t'})$, we have

$$X_{t,s} = (\alpha(t,s), \beta(t,s)) \subseteq (\alpha(t',s), \beta(t',s)) = X_{t',s},$$

as needed.

To see (2): let $t \in T$. By (b), there is $t' \in L$, such that $\pi(X_t) \subseteq \pi(X_{t'})$. Take any $t'' > t'$ in $L$ such that $\alpha(t'', x_1) \leq \alpha(t, x_1)$ and $\beta(t'', x_1) \geq \beta(t, x_1)$ (which exists in all subcases). Now, by linearity, for any $s \in M^{n-1}$, we have that $\alpha(t'', s) \leq \alpha(t, s)$ and $\beta(t'', s) \geq \beta(t, s)$. Observe also that still $\pi(X_t) \subseteq \pi(X_{t''})$, by $(a)$. Hence again

$$X_{t,s} = (\alpha(t,s), \beta(t,s)) \subseteq (\alpha(t'',s), \beta(t'',s)) = X_{t'',s},$$

as needed. □

*Remark* A.7. It is not very hard to prove that for every purely unbounded linear cell $T$, there is a line $L$ containing 0, such that for every $p \in T$, $p + L$ witnesses that $T$ is purely unbounded. We omit the proof as it is a (much easier) version of the proof of Proposition 5.15 for Presburger cells. Therefore, every semilinear family as in Proposition A.6 has a nesting direction.

### A.3. **Reconfiguring semilinear families.**

**Lemma A.8.** *Let $\{X_t\}_{t\in T}$ be a semilinear family of non-empty sets $X_t \subseteq M^n$, with $T \subseteq M^m$, such that $X = \bigcup_{t\in T}\{t\} \times X_t$ is a linear cell. Then there is a linear map $f : T \subseteq M^m \to M^n$ with $f(t) \in X_t$, such that for $Y_t = X_t - f(t)$, $t \in T$, the set*

$$Y = \bigcup_{t\in T}\{t\} \times Y_t$$

*is a linear cell. In particular, $0 \in \bigcap_{t\in T} Y_t \neq \emptyset$.*

*Proof.* The proof is an extension of the proof of definable choice in o-minimal structures ([17, Proposition 6.1.2]). Indeed, we need to choose $f(t) \in X_t$ so that, in addition, $\bigcup_{t\in T}\{t\} \times (X_t - f(t))$ remains a linear cell. We work by induction on $n$.

For $n = 0$, each $X_t = M^0 = \{0\}$, and we can let $Y_t = X_t$. For $n > 1$, the cell $X$ is either the graph of a function or a cylinder. Let $X'$ denote the projection of $X$ onto the first $m + n - 1$ coordinates, and $\pi : M^n \to M^{n-1}$ the projection map onto the first $n - 1$ coordinates. By Inductive Hypothesis, there is a linear map $g : T \to M^{n-1}$ with $g(t) \in \pi(X_t)$, such that for $G_t = \pi(X_t) - g(t)$, the set

$$G = \bigcup_{t\in T}\{t\} \times G_t$$

is a linear cell.

From the cylinder case, we only handle the subcase where $X$ is of the form $X = (\alpha, \beta)_{X'}$, with $\alpha, \beta : X \to M$ linear maps, as the other cases are similar . In particular, for every $t \in T$,

$$X_t = (\alpha(t, -), \beta(t, -))_{\pi(X_t)}.$$

Define the linear maps $h : T \to M$ and $f : T \to M^n$, via

$$h(t) = \frac{\alpha(t, g(t)) + \beta(t, g(t))}{2} \quad \text{and} \quad f(t) = (g(t), h(t)) \in X_t.$$

For the graph case, $X = \Gamma(\alpha)$, define the linear maps $h : T \to M$ and $f : T \to M^n$,

$$h(t) = \alpha(t, g(t)) \quad \text{and} \quad f(t) = (g(t), h(t)) \in X_t.$$

The verifications that $f$ is as needed are left to the reader. □

**Note:** The above proof goes through in any o-minimal structure $\mathcal{M}$, yielding the same lemma after replacing 'semilinear', 'linear cell' and 'linear map' by 'definable in $\mathcal{M}$', 'cell' and 'continuous', respectively.

School of Mathematics, University of Leeds, Leeds LS2 9JT, United Kingdom
*Email address*: `p.eleftheriou@leeds.ac.uk`

Department of Mathematics, University of Maryland, College Park
*Email address*: `aris@umd.edu`